\documentclass{amsart}
\usepackage{amsmath,amssymb, amsbsy}
\usepackage{subfigure}
\usepackage{graphpap,latexsym,epsf}
\usepackage{color,psfrag}
\usepackage[dvips]{graphicx}
\newcommand{\R}{{\mathbb R}}
\newcommand{\N}{{\mathbb N}}

\newcommand{\SN}{{\mathbb S}^{N-1}}

\newcommand{\e }{\varepsilon}
\newcommand{\Di}{{\mathcal D}^{1,2}(\R^N)}
\newcommand{\alchi}{\raisebox{1.7pt}{$\chi$}}

\renewcommand{\geq }{\geqslant}

\renewcommand{\leq }{\leqslant}

\newenvironment{pf}{\noindent{\sc Proof}.\enspace}{\hfill\qed\medskip}
\newenvironment{pfn}[1]{\noindent{\bf Proof of
    {#1}.\enspace}}{\hfill\qed\medskip}
\newtheorem{Theorem}{Theorem}[section]
\newtheorem{Corollary}[Theorem]{Corollary}
\newtheorem{Lemma}[Theorem]{Lemma}
\newtheorem{Proposition}[Theorem]{Proposition}
\theoremstyle{definition} 
\newtheorem{Example}[Theorem]{Example}

\newtheorem{remark}[Theorem]{Remark}
\begin{document}

\title[Schr\"odinger equations with many-particle and cylindrical
potentials]{On the behavior at collisions of solutions to
  Schr\"odinger equations with many-particle and cylindrical potentials}

\author[Veronica Felli \and Alberto Ferrero \and Susanna
Terracini]{Veronica Felli \and Alberto Ferrero \and Susanna Terracini}
\address{\hbox{\parbox{5.7in}{\medskip\noindent{Universit\`a di Milano
        Bicocca,\\
        Dipartimento di Ma\-t\-ema\-ti\-ca e Applicazioni, \\
        Via Cozzi
        53, 20125 Milano, Italy. \\[3pt]
        \em{E-mail addresses: }{\tt veronica.felli@unimib.it,
          alberto.ferrero@unimib.it, susanna.terracini@unimib.it}.}}}}

\date{}


 \begin{abstract}
   \noindent The asymptotic behavior of solutions to Schr\"odinger
   equations with singular homogeneous potentials is investigated.
   Through an Almgren type monotonicity formula and separation of
   variables, we describe the exact asymptotics near the singularity
   of solutions to at most critical semilinear elliptic equations with
   cylindrical and quantum multi-body singular
   potentials. Furthermore, by an iterative Brezis-Kato procedure,
   pointwise upper estimate are derived.
 \end{abstract}

\maketitle

\section{Introduction}\label{intro}

The purpose of the present paper is to describe the behavior of
solutions to a class of  Schr\"odinger equations with singular homogeneous
potentials including cylindrical and quantum multi-body ones.

The interaction between $M$ particles of coordinates $y^1,\dots, y^M$
in $\R^k$ is described in classical mechanics by potentials of the
form
\begin{equation*}
 V(y^1,\dots,y^M)=\sum_{\underset{j<m}{j,m=1}}^MV_{j,m}(y^j-y^m)
\end{equation*}
where $V_{j,m}(y)\to 0$ as $|y|\to +\infty$, see
\cite{HS}. From the mathematical point of view,
a particular interest arises in the case of inverse square
potentials $V_{j,m}(y)=\frac{\lambda_j \lambda_m}{|y|^2}$, since they
have the same order of homogeneity as the laplacian thus making the
corresponding Schr\"odinger operator invariant by scaling.
Schr\"odinger equations with the  resulting $M$-body potential
\begin{equation} \label{iniziale}
 V(y^1,\dots,y^M)=\sum_{\underset{j<m}{j,m=1}}^M
\frac{\lambda_j \lambda_m}{|y^j-y^m|^2},\quad \lambda_j,\lambda_m\in\R,
\end{equation}
have been studied by several authors; we mention in particular
\cite{HHLT} where many-particle Hardy inequalities are proved and
\cite{CSW} where the existence of ground state solutions for
semilinear Schr\"odinger equations with potentials of type
(\ref{iniziale}) is investigated. It is worth  pointing out that
hamiltonians with singular potentials having the same homogeneity
as the operator arise in relativistic quantum mechanics, see \cite{LT}.

There is a natural relation between 2-particle potentials
(\ref{iniziale}) and cylindrical potentials, whose singular set is
some $k$-codimensional subspace of the configuration space. Indeed, in
the special case $M=2$,
after the change of variables in $\R^{2k}$
\begin{equation}\label{eq:24}
z^1=\frac{1}{\sqrt 2}(y^1-y^2), \quad
z^2=\frac{1}{\sqrt 2}(y^1+y^2),
\end{equation}
the potential $V(y^1,y^2)=\frac{\lambda_1\lambda_2}{|y_1-y_2|^2}$
takes the form
\begin{equation}\label{speciale}
\frac{\lambda_1\lambda_2}{2|z^1|^2}.
\end{equation}
Elliptic equations with cylindrical inverse square potentials arise in
several fields of applications, e.g. in the search for solitary waves
with no vanishing angular momentum of nonlinear evolution equations of
Schr\"odinger and Klein-Gordon type, see \cite{BBR}.  In the recent
literature, many papers have been devoted to the study of semilinear
elliptic equations with cylindrical potentials; we mention among
others \cite{BBR,BR,BT,MFS,Musina}.  We point out that cylindrical
type (and a fortiori many-particle) potentials give rise to
substantially major difficulties with respect to the case of an
isolated singularity, because in the cylindrical/many-particle case 
 separation of variables (radial and angular) does not
actually ``eliminate'' the singularity, being the angular part of
the operator also singular.

We consider both linear and semilinear Schr\"odinger equations
with singular homogeneous potentials belonging to a class
 including as particular cases both
\eqref{iniziale} and \eqref{speciale}.
For every $3\leq k\leq N$, let us define the sets
$$
 {\mathcal
A}_k:=\big\{J\subseteq\{1,2,\dots,N\}\text{ such that }\#J=k\big\}
$$
and
$$
{\mathcal B}_k:=\{(J_1,J_2)\in {\mathcal A}_k\times {\mathcal A}_k
 \text{ such that }
J_1\cap J_2=\emptyset \text{ and } J_1<J_2 \}
$$
where $\# J$ stands for the cardinality of $J$ and $J_1<J_2$
stands for the ``alphabetic ordering'' for multi-indices (see
the list of notations at the end of this section).

In the sequel, for every $x=(x_1,x_2,\dots,x_N)\in \R^N$ and
$J\in{\mathcal A}_k$, we denote as $x_J$ the $k$-uple $(x_i)_{i\in
J}$ so that $|x_J|^2=\sum_{i\in J}x_i^2$. In a similar way, for
 any $x\in \R^N\setminus \{0\}$ and $J\in{\mathcal A}_k$ we write
$\theta_J=\frac{x_J}{|x|}$.
 Moreover we denote
\begin{align} \label{sigma}
\Sigma:=&\{ (\theta_1,\dots,\theta_{N})\in {\mathbb
  S}^{N-1}: \theta_J=0 \text{ for some }J\in{\mathcal A}_k\}\\
\notag &  \cup \{(\theta_1,\dots,\theta_N)\in \SN:
\theta_{J_1}=\theta_{J_2} \text{ for some } (J_1,J_2)\in \mathcal
B_k\}
\end{align}
and
\begin{equation}\label{eq:sigmatilde}
  \widetilde\Sigma=\{x\in\R^N\setminus\{0\}:x/|x|\in \Sigma\}\cup
  \{0\}.
\end{equation}
The potentials we are going to consider are of the type
\begin{equation}\label{eq:V(x)}
V(x)=\sum_{J\in{\mathcal
A}_k}\frac{\alpha_J}{|x_J|^2}+\sum_{(J_1,J_2)\in \mathcal B_k}
\frac{\alpha_{J_1 J_2}}{|x_{J_1}-x_{J_2}|^2},\quad
\text{for all }x\in \R^N\setminus   \widetilde\Sigma,
\end{equation}
where $\alpha_J,\alpha_{J_1J_2}\in\R$.  We notice that
${\mathcal B}_k$ is empty whenever $k>\frac{N}2$; in such a case
we consider potentials $V$ with only the cylindrical part, i.e.
with only the first summation at right hand side of~(\ref{eq:V(x)}).

Letting, for all $\theta \in\SN\setminus\Sigma$,
\begin{equation}\label{eq:ateta}
a(\theta)=\sum_{J\in{\mathcal
A}_k}\frac{\alpha_J}{|\theta_J|^2}+\sum_{(J_1,J_2)\in \mathcal
B_k}\frac{\alpha_{J_1\, J_2}}{|\theta_{J_1}-\theta_{J_2}|^2} \not
\equiv 0,
\end{equation}
we can  write the potential $V$ in (\ref{eq:V(x)}) as
$$
V(x)=\frac{a(\frac{x}{|x|})}{|x|^2}
$$
and the associated hamiltonian as
$$
\mathcal L_{a}=-\Delta-\frac{a\big(\frac{x}{|x|}\big)}{|x|^2}.
$$
As a natural setting to study the properties of operators
${\mathcal L}_{a}$, we introduce the functional space $\Di$ defined
as the completion of $C^\infty_{\rm c}(\R^N)$ with respect to the
Dirichlet norm
$$
\|u\|_{\Di}:=\bigg(\int_{\R^N}|\nabla u(x)|^2\,dx\bigg)^{\!\!1/2}.
$$
The potential $V$ in (\ref{eq:V(x)}) satisfies a Hardy type
inequality. Indeed, it was proved in \cite{Mazja} (see also
\cite{BT} and \cite{SSW}) that the
following Hardy's inequality for cylindrically singular potentials
holds:
\begin{equation}\label{eq:ssw}
\bigg(\frac{k-2}2\bigg)^{\!\!2}\int_{\R^N}\frac{|u(x)|^2}{|x_J|^2}\,dx
\leq \int_{\R^N}|\nabla u(x)|^2\,dx
\end{equation}
for all $u\in \Di$ and $J\in{\mathcal A}_k$, being the constant
$\big(\frac{k-2}2\big)^{\!2}$ optimal. Using a change of variables
of type (\ref{eq:24}), from
\eqref{eq:ssw} it follows the ``two-particle Hardy inequality'':
\begin{equation}\label{eq:MPHardy}
\frac{(k-2)^2}2
\int_{\R^N}\frac{|u(x)|^2}{|x_{J_1}-x_{J_2}|^2}\,dx \leq
\int_{\R^N}|\nabla u(x)|^2\,dx
\end{equation}
for all $u\in \Di$ and $(J_1,J_2)\in \mathcal B_k$, being the constant
$\frac{(k-2)^2}2$ optimal. From (\ref{eq:ssw}) and (\ref{eq:MPHardy})
we deduce that the potential $V$ in (\ref{eq:V(x)}) satisfies the
following ``many-particle Hardy inequality'':
\begin{equation}\label{eq:sswJ}
\bigg(\frac{k-2}2\bigg)^{\!\!2} \int_{\R^N}V(x)|u(x)|^2\,dx \leq
\bigg(\sum_{J\in{\mathcal A}_k}\alpha_J^+ +\sum_{(J_1,J_2)\in
\mathcal B_k} \alpha_{J_1 J_2}^+ \bigg)\int_{\R^N}|\nabla
u(x)|^2\,dx
\end{equation}
for all $u\in \Di$, where $\alpha_J^+=\max\{\alpha_J,0\}$ and
$\alpha_{J_1J_2}^+=\max\{\alpha_{J_1J_2},0\}$. We refer to \cite{HHLT}
for a deep analysis of many-particle Hardy inequalities and
related best constants.

In order to
discuss the positivity properties of the Schr\"odinger operator
${\mathcal L}_{a}$ in $\Di$, we consider the best constant in the
Hardy-type inequality (\ref{eq:sswJ}), i.e.
\begin{equation}
\label{eq:bound}
\Lambda(a):=\sup_{u\in\Di\setminus\{0\}}\dfrac{{\displaystyle
{\int_{\R^N}{{|x|^{-2}}{a(x/|x|)}\,u^2(x)\,dx}}}}
{{\displaystyle{\int_{\R^N}{|\nabla u(x)|^2\,dx}}}}.
\end{equation}
By (\ref{eq:sswJ}),  $\Lambda(a)\leq
\frac4{(k-2)^2}(\sum_{J\in{\mathcal A}_k}\alpha_J^+ +
\sum_{(J_1,J_2)\in {\mathcal B}_k} \alpha_{J_1 J_2}^+)$. It is
easy to verify that the quadratic form associated to ${\mathcal
L}_{a}$ is positive definite in $\Di$ if and only if
\begin{equation}\label{eq:lambdamin1}
\Lambda(a)<1.
\end{equation}
The relation between the value $\Lambda(a)$ and the first eigenvalue
of the angular component of the operator on the unit $(N-1)$-dimensional
sphere $\mathbb S^{N-1}$ is discussed in Lemma \ref{l:pos}. More
precisely, Lemma \ref{l:pos} ensures that the quadratic form
associated to ${\mathcal L}_{a}$ is positive definite if and only if
$$
\mu_1(a)>-\bigg(\frac{N-2}2\bigg)^{\!\!2},
$$
where  $\mu_1(a)$ is the first eigenvalue of
 the operator $L_{a}:=-\Delta_{\mathbb S^{N-1}}- a$ on the sphere $\mathbb S^{N-1}$.
The spectrum of the angular operator  $L_{a}$ is
discrete and consists in  a nondecreasing sequence of  eigenvalues
\[
\mu_1(a)\leq\mu_2(a)\leq\cdots\leq\mu_k(a)\leq\cdots
\]
diverging to $+\infty$, see Lemma \ref{l:spe}.

We study nonlinear equations obtained as perturbations
of the operator ${\mathcal L}_{a}$ in a bounded domain $\Omega\subset \R^N$
containing the origin. More
precisely, we deal with semilinear equations of the type
\begin{equation} \label{u}
\mathcal L_{a} u=h(x)\, u+f(x,u), \quad \text{in } \Omega \, .
\end{equation}
We assume that the linear perturbing potential $h$ is negligible with
respect to the potential $V$ near the collision singular set
$\widetilde \Sigma$ defined in (\ref{eq:sigmatilde}), in the sense
that there exist $C_h>0$ and $\varepsilon>0$ such that, for a.e. $x\in
\Omega\setminus \widetilde\Sigma$,
\begin{equation}\label{H1-2}\tag{$\bf H$}
h\in W^{1,\infty}_{\rm loc}\big(\Omega\setminus
  \widetilde\Sigma \big) \ \text{and}\
|h(x)|+|\nabla h(x)\cdot x|\leq C_h \bigg(\sum_{J\in{\mathcal
A}_k}|x_{J}|^{-2+\varepsilon}+\sum_{(J_1,J_2)\in \mathcal B_k}
|x_{J_1}-x_{J_2}|^{-2+\varepsilon}\bigg).
\end{equation}
We notice that it is not restrictive to assume $\e\in(0,1)$ in
(\ref{H1-2}).

As far as the nonlinear perturbation is concerned, we assume that $f$
satisfies
\begin{equation}\label{F}\tag{$\bf F$}
\left\{\!\!
\begin{array}{l}
f\in C^0(\Omega\times \R),\quad F\in C^1(\Omega\times \R),
\quad s\mapsto f(x,s)\in C^1(\R)\text{ for a.e. }x\in\Omega, \\[5pt]
|f(x,s)s|+|f'_s(x,s)s^2|+|\nabla_x F(x,s)\cdot x|\leq C_f(|s|^2+|s|^{2^*})
\quad\text{for a.e. $x\in\Omega$ and all $s\in\R$},
\end{array}
\right.
\end{equation}
 where $F(x,s)=\int_0^s f(x,t)\,dt$,
$2^*=\frac{2N}{N-2}$ is the critical Sobolev exponent, $C_f>0$ is a
constant independent of $x\in\Omega$ and $s\in\R$, $\nabla_x F$
denotes the gradient of $F$ with respect to the $x$ variable, and
$f'_s(x,s)=\frac{\partial f}{\partial s}(x,s)$.

We say
that a function $u\in H^1(\Omega)$ is a $H^1(\Omega)$-weak
solution to (\ref{u}) if, for all $w\in H^1_0(\Omega)$,
\begin{equation*}
{\mathcal Q}_{a}^{\Omega}(u,w)
=\int_{\Omega} h(x)
  u(x)w(x)\,dx+\int_\Omega f(x,u(x))w(x)\, dx,
\end{equation*}
where ${\mathcal Q}_{a}^{\Omega}:H^1(\Omega)\times
H^1(\Omega)\to\R$ is defined by
\[
{\mathcal Q}_{a}^{\Omega}(u,w):= \int_{\Omega}
\nabla u(x)\cdot\nabla w(x)\,dx -\int_{\Omega}
\frac{a(x/|x|)}{|x|^2}\,u(x)w(x)\,dx.
\]

Schr\"odinger equations with inverse square homogeneous singular
potentials can be regarded as critical from the mathematical point of
view, as they do not belong to the Kato class.  A rich literature
deals with such critical equations, both in the case of one isolated
pole, see e.g. \cite{egnell, FG, GP, Jan, SM, terracini96}, and in
that of multiple singularities, see \cite{esteban, chen, duyckaerts,
  FMT1, FT}.  The analysis of fundamental spectral properties such as
essential self-adjointness and positivity carried out in
\cite{FMT1,FMT3} for Schr\"odinger operators with isolated inverse
square singularities, highlighted how the asymptotic behavior of
solutions to associated elliptic equations near the singularity plays
a crucial role.  A precise evaluation of the asymptotics of solutions
turned out to be an important tool also to establish existence of
“ground states” for nonlinear Schr\"odinger equations with
multi-singular Hardy potentials (see \cite{FT}) and of solutions to
nonlinear systems of Schr\"odinger equations with Hardy potentials
\cite{AFP}.  A first result about the study of the asymptotic behavior
of solutions near isolated singularities is contained in \cite{FS3},
where H\"older continuity of solutions to degenerate elliptic
equations with singular weights has been established thus allowing the
evaluation of the exact asymptotic behavior of solutions to
Schr\"odinger equations with Hardy potentials near the pole. An
extension to the case of Schr\"odinger equations with dipole-type
potentials (namely purely angular multiples of inverse square
potentials) has been obtained in \cite{FMT2} by separation of
variables and comparison principles, and later generalized to
Schr\"odinger equations with singular homogeneous electromagnetic
potentials of Aharonov-Bohm type \cite{FFT} by the Almgren
monotonicity formula.  Comparison and maximum principles play a
crucial role also in \cite{pinchover94}, where the existence of the
limit at the singularity of any quotient of two positive solutions to
Fuchsian type elliptic equations is proved.  We mention that
related asymptotic expansions near singularities were obtained in
\cite{mazzeo91, mazzeo91-2} for elliptic equations on manifolds with
conical singularities by Mellin transform methods (see also
\cite{lesch}); we refer to \cite{addendumFFT} for a comparison between
such results and asymptotics via Almgren monotonicity methods. It is
also worth citing \cite{buslaev_levin}, where some asymptotic formulas
are heuristically obtained for the three-body one-dimensional problem.

Due to the presence of  multiple collisions,
one should expect that solutions to equations (\ref{u})
behave singularly at the origin: our purpose is to
describe the rate and the shape of the singularity of
solutions, by relating them to the eigenvalues and the eigenfunctions of
the angular operator  $L_{a}$  on the sphere ${\mathbb S}^{N-1}$.

The following theorem provides a classification of the behavior of
any solution $u$ to (\ref{u}) near the singularity based on the limit as
$r\to 0^+$ of the \emph{Almgren's frequency} function (see \cite{almgren,GL})
\begin{equation}\label{eq:14}
  {\mathcal N}_{u,h,f}(r)=\frac{r\int_{B_r} \big(|\nabla u(x)|^2
-\frac{a(x/|x|)}{|x|^2}u^2(x) -
    h(x)u^2(x)-f(x,u(x)) \big) \, dx}{\int_{\partial B_r}|u(x)|^2 \, dS},
  \end{equation}
where, for any $r>0$, $B_r$ denotes the ball $\{x\in\R^N: |x|<r\}$.

\begin{Theorem}\label{Main-h}
  Let $u\not\equiv 0$ be a nontrivial weak $H^1(\Omega)$-solution to
  (\ref{u}) in a bounded open set $\Omega\subset\R^N$ containing $0$,
  $N\geq k\geq 3$, with $a$ satisfying (\ref{eq:ateta}) and
  (\ref{eq:lambdamin1}), $h$ satisfying \eqref{H1-2}, and $f$
  satisfying \eqref{F}.  Then, letting $\mathcal N_{u,h,f}(r)$ as in
  \eqref{eq:14}, there exists $k_0\in \N$, $k_0\geq 1$, such that
\begin{equation}\label{lim-N}
  \lim_{r\to 0^+} \mathcal
  N_{u,h,f}(r)=-\frac{N-2}2+\sqrt{\left(\frac{N-2}2\right)^{\!\!2}+\mu_{k_0}(a)}.
\end{equation}
Furthermore, if $\gamma$ denotes the limit in (\ref{lim-N}), $m\geq 1$
is the multiplicity of the eigenvalue $\mu_{k_0}(a)$ and
$\{\psi_i:j_0\leq i\leq j_0+m-1\}$ ($j_0\leq k_0\leq j_0+m-1$) is an
$L^2(\SN)$-orthonormal basis for the eigenspace associated to
$\mu_{k_0}(a)$, then
\begin{equation} \label{convergence}
\lambda^{-\gamma} u(\lambda x)\to |x|^\gamma\sum_{i=j_0}^{j_0+m-1}
\beta_i \psi_i\bigg(\frac{x}{|x|}\bigg) \qquad \text{in } H^1(B_1) \quad
\text{as } \lambda\to 0^+
\end{equation}
where
\begin{align}\label{eq:38}
  \beta_i= \int_{{\mathbb S}^{N-1}}\!\bigg[ R^{-\gamma}u(R\theta)+
  \int_{0}^R\frac{ h(s\theta)u(s\theta)
+f\big(s\theta,u(s\theta)\big)}{2\gamma+N-2}
  \bigg(s^{1-\gamma}-\frac{s^{\gamma+N-1}}{R^{2\gamma+N-2}}\bigg)ds
  \bigg]\psi_{i}(\theta)\,dS(\theta),
 \end{align}
 for all $R>0$ such that $\overline{B_{R}}=
\{x\in\R^N:|x|\leq R\}\subset\Omega$
and $(\beta_{j_0},\beta_{j_0+1},\dots,\beta_{j_0+m-1})\neq(0,0,\dots,0)$.
\end{Theorem}
Due to the homogeneity of the potentials, Schr\"odinger operators
${\mathcal L}_{a}$ are invariant by the Kelvin transform,
\begin{equation*}
\tilde u(x)=|x|^{-(N-2)}u\bigg(\frac{x}{|x|^2}\bigg),
\end{equation*}
which is an isomorphism of $\Di$. Indeed, if $u\in H^1(\Omega)$
weakly solves (\ref{u}) in a bounded open set $\Omega$ containing $0$,
then its Kelvin's transform $\tilde u$ weakly solves (\ref{u}) with
$h$ replaced by $|x|^{-4}h(\frac x{|x|^2})$ 
and $f(x,\cdot)$ replaced by $|x|^{-N-2}f\big(\frac{x}{|x|^2},|x|^{N-2}\cdot\big)$
in the external domain
$\widetilde \Omega=\big\{x\in\R^N: x/{|x|^2}\in\Omega\big\}$.  
Therefore, under suitable decay conditions on $h$ at $\infty$ and proper 
subcriticality assumptions on $f$, the asymptotic behavior at infinity 
of solutions to 
(\ref{u}) in external domains can be easily deduced from Theorem 
\ref{Main-h}  and the Kelvin transform (see \cite[Theorems 1.4 and 1.6]{FFT}).

A major breakthrough in the description of the singularity of
solutions at zero can be done by evaluating the behavior of
eigenfunctions $\psi_i$; indeed such eigenfunctions solve an elliptic
equation on $\SN$ exhibiting itself a potential which is singular on
$\Sigma$. After a stereographic projection of $\SN$ onto $\R^{N-1}$,
the equation satisfied by each $\psi_i$ takes a form which is similar
to (\ref{u}) in a lowered dimension with a potential whose singular
set is $(N-1-k)-$dimensional and to which we can apply the above
theorem to deduce a precise asymptotics in terms of eigenvalues and
eigenfunctions of an operator on ${\mathbb S}^{N-2}$; the procedure
can be iterated $(N-k)-$times until we come to an equation with a
potential with isolated singularities whose corresponding angular
operator is no more singular. A detailed analysis of the asymptotic
behavior of eigenfunctions is performed in section
\ref{sec:asympt-behav-eigenf}.

A pointwise upper estimate on the behavior of solutions can be derived
by a Brezis-Kato type iteration argument, see \cite{BK}. More
precisely, we can estimate the solutions by terms of the first
eigenvalue and eigenfunction of the angular potential $\hat a$ obtained 
by summing up only the positive contributions of $a$, i.e.
\begin{equation}\label{eq:atetahat}
\hat a(\theta)=\sum_{J\in{\mathcal
A}_k}\frac{\alpha_J^+}{|\theta_J|^2}+\sum_{(J_1,J_2)\in \mathcal
B_k}\frac{\alpha_{J_1\, J_2}^+}{|\theta_{J_1}-\theta_{J_2}|^2}.
\end{equation}
Under the assumption 
\begin{equation}\label{eq:53}
  \Lambda(\hat a)
  =
  \sup_{u\in\Di\setminus\{0\}}\dfrac{\int_{\R^N}{{|x|^{-2}}{\hat a(x/|x|)}\,u^2(x)\,dx}}
  {\int_{\R^N}{|\nabla u(x)|^2\,dx}}
  <1,
\end{equation}
by Lemma \ref{l:pos} the number 
\begin{equation}\label{def-sigma-hat}
\hat\sigma=-\frac{N-2}{2} +
\sqrt{\left(\frac{N-2}2\right)^{\!\!2}+\mu_1(\hat a)}
\end{equation}
is well defined. We denote as $\hat \psi_1\in H^1(\SN)$, $\|\hat
\psi_1\|_{L^2(\SN)}=1$,  the first positive $L^2-$normalized eigenfuntion of the
eigenvalue problem
$L_a \psi=\mu_1(\hat a)\psi$ in   $\SN$.

\begin{Theorem}\label{t:BK}
  Let $u$ be a weak $H^1(\Omega)$-solution to (\ref{u}) in a bounded
  open set $\Omega\subset\R^N$ containing $0$, $N\geq k\geq 3$, with
  $a$ satisfying (\ref{eq:ateta}) and $\hat a$ as in
  (\ref{eq:atetahat}) satisfying (\ref{eq:53}). If $h$ satisfies
  \eqref{H1-2} and $f$ satisfies \eqref{F}, then for any
  $\Omega'\Subset \Omega$  there exists $C>0$ such that 
$$
|u(x)|\leq C |x|^{\hat\sigma} \hat\psi_1\Big(\frac{x}{|x|}\Big)\quad\text{for a.e. }
x\in \Omega'.
$$
\end{Theorem}
In particular, if all $\alpha_J,\alpha_{J_1J_2}$ are positive, then
$\hat a\equiv a$ and the above theorem ensures that all solutions are
pointwise bounded by $|x|^\sigma\psi_1(x/|x|)$ where $\sigma=
-\frac{N-2}{2} +
\big[(\frac{N-2}2)^{2}+\mu_1( a)\big]^{1/2}$.
On the other hand, if all $\alpha_J,\alpha_{J_1J_2}$ are negative, then
$\hat a\equiv 0$ and the above theorem implies that all solutions are
bounded.

The paper is organized as follows. In section
\ref{sec:hardys-cylindr-ineq} we prove some Hardy-type inequalities
with singular potentials of type (\ref{eq:V(x)}) and discuss the
relation between the positivity of the quadratic form associated to
${\mathcal L}_{a}$ and the first eigenvalue of the angular operator on
the sphere ${\mathbb S}^{N-1}$.  In section \ref{sec:pohoz-type-ident}
we derive a Pohozaev-type identity for solutions to~(\ref{u}) through
a suitable approximating procedure which allows getting rid of the
singularity of the angular potential. In Section
\ref{sec:brezis-kato-type} we deduce a Brezis-Kato estimate to prove
an a-priori super-critical summability of solutions to (\ref{u}) which
allows us to include the critical growth case in the Almgren type
monotonicity formula which is obtained in Section~\ref{sec:monot-prop}
and which is used in section \ref{sec:blow-up-argument} together with
a blow-up method to prove Theorem \ref{Main-h}.  Section
\ref{sec:asympt-behav-eigenf} is devoted to the study of the
asymptotic behavior of the eigenfunctions of the angular
operator. Section \ref{sec:pointw-estim-under} contains some
Brezis-Kato estimates in weighted Sobolev spaces which allow proving
Theorem \ref{t:BK}.  A final appendix contains a Pohozaev-type
identity for semilinear elliptic equations with an anisotropic 
inverse-square potential with a bounded angular coefficient.

\medskip
\noindent
{\bf Notation. } We list below some notation used throughout the
paper.\par
\begin{itemize}
\item[-] For all $r>0$, $B_r$ denotes the ball $\{x\in\R^N: |x|<r\}$
  in $\R^N$ with center at $0$ and radius $r$.
\item[-] For all $r>0$, $\overline{B_{r}}= \{x\in\R^N:|x|\leq r\}$
denotes the closure of $B_r$.
\item[-] $dS$ denotes the volume
element on the spheres $\partial B_r$, $r>0$.
\item[-] If
$J_1=\{j_{1,1},\dots, j_{1,k}\}$ and $J_2=\{j_{2,1},\dots,
j_{2,k}\}$ are two multi-indices of $k$ elements, by $J_1<J_2$ we
mean that there exists $n\in \{1,\dots,k\}$ such that
$j_{1,i}=j_{2,i}$ for any $1\leq i\leq n-1$ and $j_{1,n}<j_{2,n}$.
\item[-] For all $t\in\R$,  $t^+=t_+:=\max\{t,0\}$ (respectively
$t^-=t_-:=\max\{-t,0\}$) denotes the positive (respectively negative)
part of $t$.
\item[-] $S=\inf_{v\in\mathcal D^{1,2}(\R^N)\setminus\{0\}}\|\nabla
  v\|_{L^2}^2 \|v\|_{L^{2^*}}^{-2}$ denotes the best constant in the
  classical Sobolev's embedding.
\end{itemize}

\section{Hardy type inequalities}\label{sec:hardys-cylindr-ineq}

The following Hardy's inequality on the unit sphere holds.
\begin{Lemma}\label{l:hardysphere}
  Let $a$ as in (\ref{eq:ateta}). For every $\psi\in H^1({\mathbb
    S}^{N-1})$ there holds
\begin{align*}
\bigg(\frac{k-2}2\bigg)^{\!\!2}\!\! & \int_{{\mathbb
    S}^{N-1}}\!\!a(\theta)|\psi(\theta)|^2\,dS
\\
&    \leq
 \bigg(\sum_{J\in{\mathcal A}_k}\alpha_J^+ + \sum_{(J_1,J_2)\in \mathcal
 B_k} \alpha_{J_1\, J_2}^+ \bigg)\bigg[\int_{{\mathbb
    S}^{N-1}}\!|\nabla_{{\mathbb
    S}^{N-1}}\psi(\theta)|^2\,dS+\bigg(\frac{N-2}2\bigg)^{\!\!2}
\!\!\int_{{\mathbb
    S}^{N-1}}|\psi(\theta)|^2\,dS\bigg].
\end{align*}
\end{Lemma}
\begin{pf}
Let $\psi\in H^1({\mathbb S}^{N-1})$ and $\phi\in C^\infty_{\rm
c}(0,+\infty)$. Rewriting inequality (\ref{eq:sswJ}) for
$u(x)=\phi(r)\psi(\theta)$, $r=|x|$, $\theta=\frac x{|x|}$, we
obtain that
\begin{align*}
  \bigg(\frac{k-2}2\bigg)^2&\bigg(\int_0^{+\infty
  }\frac{r^{N-1}}{r^2}\,\phi^2(r)\,dr\bigg)\bigg( \int_{{\mathbb
      S}^{N-1}}a(\theta)|\psi(\theta)|^2\,dS\bigg)\\
&\leq \bigg(\sum_{J\in{\mathcal A}_k}\alpha_J^+ +
\sum_{(J_1,J_2)\in \mathcal
 B_k} \alpha_{J_1\, J_2}^+\bigg)\bigg(\int_0^{+\infty
  }r^{N-1}|\phi'(r)|^2\,dr\bigg)\bigg( \int_{{\mathbb
      S}^{N-1}}|\psi(\theta)|^2\,dS\bigg)\\
&\quad+\bigg(\sum_{J\in{\mathcal A}_k}\alpha_J^+ +
\sum_{(J_1,J_2)\in \mathcal
 B_k} \alpha_{J_1\, J_2}^+\bigg)\bigg(\int_0^{+\infty
  }\frac{r^{N-1}}{r^2}\,\phi^2(r)\,dr\bigg)\bigg( \int_{{\mathbb
      S}^{N-1}}|\nabla_{{\mathbb
      S}^{N-1}} \psi(\theta)|^2\,dS\bigg),
\end{align*}
and hence, by optimality of the classical Hardy constant,
\begin{multline*}
  \bigg(\frac{k-2}2\bigg)^2\bigg( \int_{{\mathbb
      S}^{N-1}}a(\theta)|\psi(\theta)|^2\,dS\bigg)\\
\leq \bigg(\sum_{J\in{\mathcal A}_k}\alpha_J^+ +
\sum_{(J_1,J_2)\in \mathcal
 B_k} \alpha_{J_1\, J_2}^+\bigg)\bigg[
\bigg(\int_{{\mathbb
      S}^{N-1}}|\psi(\theta)|^2\,dS\bigg)
\inf_{\phi\in C^\infty_{\rm c}(0,+\infty)}
\frac{\int_0^{+\infty
  }r^{N-1}|\phi'(r)|^2\,dr}{\int_0^{+\infty
  }r^{N-3}\,\phi^2(r)\,dr}
\\
+ \int_{{\mathbb
      S}^{N-1}}|\nabla_{{\mathbb
      S}^{N-1}} \psi(\theta)|^2\,dS\bigg]\\
=\bigg(\sum_{J\in{\mathcal A}_k}\alpha_J^+ + \sum_{(J_1,J_2)\in
\mathcal
 B_k} \alpha_{J_1\, J_2}^+\bigg)\bigg[\bigg(\frac{N-2}2\bigg)^2 \int_{{\mathbb
    S}^{N-1}}|\psi(\theta)|^2\,dS+\int_{{\mathbb
    S}^{N-1}}|\nabla_{{\mathbb
    S}^{N-1}}\psi(\theta)|^2\,dS\bigg].
\end{multline*}
The proof is thereby complete.
\end{pf}

\noindent Let us consider the following class of angular potentials
\begin{equation}\label{eq:classF}
{\mathcal F}:=\bigg\{f\in L^\infty_{\rm loc}(\mathbb S^{N-1}\setminus \Sigma):
\frac{|f(\theta)|}{\sum_{J\in{\mathcal A}_k}|\theta_J|^{-2}+
\sum_{(J_1,J_2)\in \mathcal
 B_k} |\theta_{J_1}-\theta_{J_2}|^{-2}} \in
L^\infty(\mathbb S^{N-1})\bigg\}.
\end{equation}
From Lemma \ref{l:hardysphere} we have that, for every $f\in{\mathcal F}$,
 the supremum
\begin{equation} \label{supremum}
\Lambda(f):=\sup_{\psi\in H^1(\mathbb
  S^{N-1})\setminus\{0\}}\frac{\int_{\mathbb
    S^{N-1}}f(\theta)\,\psi^2(\theta)\,dS(\theta)}{\int_{\mathbb
    S^{N-1}}|\nabla_{\mathbb
    S^{N-1}}\psi(\theta)|^2\,dS(\theta)+\big(\frac{N-2}{2}\big)^2\int_{\mathbb
    S^{N-1}} \psi^2(\theta)\,dS(\theta)}
\end{equation}
is finite.
On the other hand, arguing as in the proof of \cite[Lemma
1.1]{terracini96}, we can easily verify that
\begin{equation}\label{eq:1}
\Lambda(f)=\sup_{u\in\Di\setminus\{0\}}\dfrac{{\displaystyle
{\int_{\R^N}{{|x|^{-2}}{f(x/|x|)}\,u^2(x)\,dx}}}}
{{\displaystyle{\int_{\R^N}{|\nabla u(x)|^2\,dx}}}}.
\end{equation}
Furthermore, it is easy to verify that
$$
\Lambda(f)\geq 0$$
and
$$
\Lambda(f)= 0 \quad\text{if and only if}\quad
 f\leq 0 \text{ a.e. in }\mathbb S^{N-1}.
$$
For every $f\in{\mathcal F}$ satisfying $\Lambda(f)<1$, we can
perform a complete spectral analysis of the angular
Schr\"odinger operator $-\Delta_{\mathbb S^{N-1}}- f$ on the sphere.

\begin{Lemma}\label{l:spe}
  Let $f\in{\mathcal F}$ satisfying $\Lambda(f)<1$.  Then the spectrum of the
  operator
$$
L_f:=-\Delta_{\mathbb S^{N-1}}- f
$$
on $\mathbb S^{N-1}$ consists in a diverging sequence 
$\mu_1(f)\leq\mu_2(f)\leq\cdots\leq\mu_k(f)\leq\cdots$ of real
eigenvalues with finite multiplicity
the first of which admits the variational characterization
\begin{equation}\label{firsteig}
  \mu_1(f)=\min_{\psi\in H^1(\mathbb
    S^{N-1})\setminus\{0\}}\frac{\int_{\mathbb S^{N-1}}
    \big[\big|\nabla_{\mathbb S^{N-1}}\psi(\theta)\big|^2
    - f(\theta)|\psi(\theta)|^2\big]\,dS(\theta)}{\int_{\mathbb S^{N-1}}
    |\psi(\theta)|^2\,dS(\theta)}.
\end{equation}
Moreover $\mu_1(f)$ is simple and its associated  eigenfunctions
do not change sign in $\mathbb S^{N-1}$.
\end{Lemma}

\begin{pf} By Lemma \ref{l:hardysphere} and assumption $\Lambda(f)<1$,
  the operator $T:L^{2}({\mathbb S}^{N-1})\to L^{2}({\mathbb
    S}^{N-1})$ defined as
\[
Th=u\quad\text{if and only if}\quad
-\Delta_{\mathbb S^{N-1}}u- f u+{\textstyle{\big(\frac{N-2}2\big)^2}}u=h
\]
is well-defined, symmetric, and compact. The conclusion  follows
from classical spectral theory. In particular, we point out that 
the simplicity of the first eigenvalue follows from the fact that,
since $k>1$,  
the singular set $\Sigma$ does not disconnect the sphere.~\end{pf}

\noindent For all $f\in{\mathcal F}$, let us consider the quadratic
form associated to the Schr\"odinger operator
${\mathcal L}_f$, i.e.
$$
Q_f(u):=\int_{\R^N} |\nabla  u(x)|^2dx-
\int_{\R^N}\frac{f(x/|x|)\,u^2(x)}{|x|^2}\,dx.
$$
The problem of positivity of $Q_f$ is solved in the following lemma.
\begin{Lemma}\label{l:pos}
Let $f\in{\mathcal F}$. The following conditions are equivalent:
\begin{align*} {\rm i)}\quad&Q_f \text{ is positive definite, i.e. }
  \inf_{u\in\Di\setminus\{0\}}\frac{Q_f(u)}
  {\int_{\R^N}|\nabla u(x)|^2\,dx}>0; \\
  {\rm ii)} \quad&\Lambda(f)<1;\\
  {\rm iii)}\quad&
  \mu_1(f)>-\big({\textstyle{\frac{N-2}2}}\big)^{\!2}\text{ where $\mu_1(f)$
    is defined in (\ref{firsteig})}.
\end{align*}
\end{Lemma}
\begin{pf}
  The equivalence between i) and ii) follows from the definition of
  $\Lambda(f)$, see (\ref{eq:1}). On the other hand, arguing as in
  \cite[Proposition 1.3 and Lemma 1.1]{terracini96} (see also
  \cite[Lemmas 1.1 and 2.1]{FFT}) one can obtain equivalence between
  i) and iii).
\end{pf}

Henceforward, we shall assume that (\ref{eq:lambdamin1}) holds, so
that the quadratic form associated to the operator ${\mathcal L}_a$ is
positive definite.

\begin{Example}\label{ex}
Let us consider cylindrical potentials, i.e. the particular case in which
\begin{equation} \label{alpha-J}
\alpha_J=
\begin{cases}
\alpha,&\text{if }J=\bar J=\{1,2,\dots,k\},\\
0,&\text{if }J\not=\{1,2,\dots,k\},
\end{cases}
\qquad\text{for some $\alpha\in\R$}
\end{equation}
and
\begin{equation} \label{alpha-J1-J2}
\alpha_{J_1\, J_2}=0 \qquad \text{ for any } (J_1,J_2)\in
{\mathcal B}_k,
\end{equation}
so that $a(\theta)=\alpha/|\theta_{\bar J}|^2$. Then, from
the optimality of the constant $\big(\frac{k-2}2\big)^2$ in
(\ref{eq:ssw}), it follows that
$\Lambda(a)=\alpha^+\big(\frac2{k-2}\big)^2$ and
(\ref{eq:lambdamin1}) reads as $\alpha<\big(\frac{k-2}2\big)^2$.
Moreover there holds
\begin{equation}\label{eq:2}
\mu_1(a)=-\frac{(k-2)(N-k)}{2}
-\alpha+(N-k)\sqrt{\left(\frac{k-2}{2}\right)^{\!\!2}-\alpha}.
\end{equation}
In order to verify (\ref{eq:2}), let us set
\begin{equation} \label{defgamma'}
\gamma'=-\frac{k-2}{2}+\sqrt{\left(\frac{k-2}{2}\right)^{\!\!2}-\alpha}
\end{equation}
and consider the function $u(x)=|x_{\bar J}|^{\gamma'}
=\big(\sum_{i=1}^kx_i^2\big)^{\gamma'/2}\in H^1_{\rm loc}(\R^N)$.
Then $u$ solves the equation
\begin{equation} \label{1steig} -\Delta
  u(x)-\frac{\alpha}{|x_{\bar J}|^2}\,u(x)=0 \quad \text{in }
  \{x\in\R^N:x_{\bar J}\neq0\}.
\end{equation}
The function $u$ may be rewritten as $u(x)=|x|^{\gamma'} \psi\big(
\frac{x}{|x|}\big)$ once we define
$\psi(\theta)=|\theta_{\bar J}|^{\gamma'}$ for any $\theta\in
\SN\setminus \Sigma$. Since $u$ solves \eqref{1steig}, we
obtain
\begin{equation*}
  -\gamma'(\gamma'+N-2)r^{\gamma'-2}
  \psi(\theta)-r^{\gamma'-2}\Delta_{\SN}\psi(\theta)=
  r^{\gamma'-2} a(\theta)\psi(\theta), \quad \text{for any } r>0
  \text{ and } \theta\in\SN\setminus \Sigma.
\end{equation*}
This yields
$$
-\Delta_{\SN} \psi(\theta)-
a(\theta)\psi(\theta)=\gamma'(\gamma'+N-2)\psi(\theta), \quad {\rm
in \ } \SN.
$$
This shows that $\psi$ is a positive eigenfunction of the operator
$L_a$ and hence by Lemma \ref{l:spe} the corresponding eigenvalue
must coincide with $\mu_1(a)$, i.e. $\gamma'(\gamma'+N-2)=\mu_1(a)$.
 (\ref{eq:2}) follows by~\eqref{defgamma'}.
\end{Example}

\begin{Example}\label{ex:2}
Let us also consider two-body potentials, i.e. the case in which $N\geq 2k$,
\begin{equation*}
\alpha_J=0 \qquad \text{ for any } J\in \mathcal A_k
\end{equation*}
and
\begin{equation*} 
\alpha_{J_1\, J_2}=\
\begin{cases}
\alpha,&\text{if } J_1=\bar J_1=\{1,2,\dots,k\} \text{  and }
J_2=\bar J_2=\{k+1,k+2,\dots,2k\}
,\\
0,&\text{if }(J_1,J_2)\not=(\bar  J_1,\bar J_2),
\end{cases}
\end{equation*}
so that $a(\theta)=\alpha/|\theta_{\bar
J_1}-\theta_{\bar J_2}|^2$.
The optimality of the constant
$\frac{(k-2)^2}2$ in inequality \eqref{eq:MPHardy} implies that
$\Lambda(a)=\alpha^+ \frac{2}{(k-2)^2}$ and condition
\eqref{eq:lambdamin1} reads as $\alpha<\frac{(k-2)^2}2$. Moreover
we have
\begin{equation} \label{eq2-bis}
\mu_1(a)=-\frac{(k-2)(N-k)}{2}
-\frac{\alpha}2+(N-k)\sqrt{\left(\frac{k-2}{2}\right)^{\!\!2}-\frac{\alpha}2}.
\end{equation}
In order to prove \eqref{eq2-bis} we put
\begin{equation} \label{gamma-second}
\gamma''=-\frac{k-2}2+\sqrt{\left(\frac{k-2}{2}\right)^{\!\!2}-\frac{\alpha}2}
\end{equation}
and we define $u(x)=|x_{\bar J_1}-x_{\bar
J_2}|^{\gamma''}\in H^1_{\rm loc}(\R^N)$. Then $u$ solves the
equation
\begin{equation} \label{1steig-bis}
-\Delta u(x)-\frac{\alpha}{|x_{\bar J_1}-x_{\bar
J_2}|^2} u(x)=0 \quad \text{ in } \{x\in\R^N:
x_{\bar J_1}\neq x_{\bar J_2} \}.
\end{equation}
Proceeding as in Example \ref{ex}, by \eqref{gamma-second} and
\eqref{1steig-bis} we conclude that $\psi(\theta)=|\theta_{\bar
J_1}-\theta_{\bar J_2}|^{\gamma''}$ is an eigenfunction of
$\mu_1(a)$ and that $\mu_1(a)$ is given by \eqref{eq2-bis}.
\end{Example}

We extend to singular potentials of the form (\ref{eq:V(x)}) the
Hardy type inequality with boundary terms proved by Wang and Zhu
in \cite{wz}.

\begin{Lemma}\label{l:hardyboundary}
  Let $a$ be as in (\ref{eq:ateta}) and assume that
  (\ref{eq:lambdamin1}) holds. Then
  \begin{multline}\label{eq:10bi}
    \int_{B_r} \bigg( |\nabla u(x)|^2-\frac{a(\frac{x}{|x|})}{|x|^2}
|u(x)|^2\bigg)\,dx+
    \frac{N-2}{2r}\int_{\partial B_r}|u(x)|^2\,dS\\
\geq \left(\mu_1(a)+\bigg(\frac{N-2}{2}\bigg)^{\!\!2}\right)
      \int_{B_r}\frac{|u(x)|^2}{|x|^2}\,dx
  \end{multline}
for all $r>0$ and $u\in H^1(B_r)$.
  \end{Lemma}
\begin{pf}
  By scaling, it is enough to prove the inequality for $r=1$.  Let
  $u\in C^\infty(B_1)\cap H^1(B_1)$ with $0\not\in
  \mathop{\rm supp}u$.  Passing to polar coordinates, we have that
\begin{align}\label{eq:4}
    \int_{B_1} \bigg( |\nabla u(x)|^2-&
\frac{a(\frac{x}{|x|})}{|x|^2}|u(x)|^2\bigg)\,dx+
    \frac{N-2}{2}\int_{\partial B_1}|u(x)|^2\,dS\\
  \notag=& \int_{{\mathbb S}^{N-1}}\bigg(
  \int_0^{1}r^{N-1}|\partial_ru(r,\theta)|^2\,dr\bigg)\,dS(\theta)
  +\frac{N-2}2\int_{{\mathbb S}^{N-1}}|u(1,\theta)|^2\,dS(\theta)\\
  \notag& \quad+ \int_0^{1}\frac{r^{N-1}}{r^2}\bigg(\int_{{\mathbb
      S}^{N-1}}\left[|\nabla_{{\mathbb S}^{N-1}}u(r,\theta)|^2-a(\theta)
    |u(r,\theta)|^2\right]\,dS(\theta)\bigg)\,dr.
\end{align}
For all $\theta\in{\mathbb S}^{N-1}$, let $\varphi_{\theta}\in
C^\infty(0,1)$ be defined by $\varphi_{\theta}(r)=u(r,\theta)$, and
$\widetilde\varphi_{\theta}\in C^\infty(B_1)$ be the radially
symmetric function given by
$\widetilde\varphi_{\theta}(x)=\varphi_{\theta}(|x|)$.  We notice that
$0\not\in \mathop{\rm supp}\widetilde\varphi_{\theta}$.
The Hardy
inequality with boundary term proved in \cite{wz} yields
\begin{align}\label{eq:5}
  \int_{{\mathbb S}^{N-1}}&\bigg(
  \int_0^{1}r^{N-1}|\partial_ru(r,\theta)|^2\,dr
+\frac{N-2}2\,|u(1,\theta)|^2\bigg)\,dS(\theta)\\
\notag&=
  \int_{{\mathbb S}^{N-1}}\bigg(
  \int_0^{1}r^{N-1}|\varphi_{\theta}'(r)|^2\,dr
+\frac{N-2}2\,|\varphi_{\theta}(1)|^2\bigg)\,dS(\theta)\\
  \notag& =\frac1{\omega_{N-1}}  \int_{{\mathbb S}^{N-1}}\bigg(
\int_{B_1}|\nabla\widetilde \varphi_{\theta}(x)|^2\,dx
+\frac{N-2}2\int_{\partial B_1}|\widetilde \varphi_{\theta}(x)|^2\,dS
\bigg)\,dS(\theta)
\\&\notag\geq
\frac1{\omega_{N-1}}  \bigg(\frac{N-2}2\bigg)^{\!\!2}
  \int_{{\mathbb S}^{N-1}}\bigg(\int_{B_1}\frac{|\widetilde
    \varphi_{\theta}(x)|^2}{|x|^2}\,dx\bigg)\,dS(\theta)\\
\notag&=\bigg(\frac{N-2}2\bigg)^{\!\!2}
  \int_{{\mathbb S}^{N-1}}\bigg(
  \int_0^{1}\frac{r^{N-1}}{r^2}|u(r,\theta)|^2\,dr\bigg)\,dS(\theta)
=\bigg(\frac{N-2}2\bigg)^{\!\!2}\int_{B_1}\frac{|u(x)|^2}{|x|^2}\,dx.
\end{align}
where $\omega_{N-1}$ denotes the volume of the unit sphere ${\mathbb
  S}^{N-1}$, i.e. $\omega_{N-1}=\int_{{\mathbb S}^{N-1}}dS(\theta)$.
On the other hand, from the definition of $\mu_1(a)$ it follows
that, for every $r\in(0,1)$,
\begin{equation}\label{eq:3}
  \int_{{\mathbb
      S}^{N-1}}\!\!\left[|\nabla_{{\mathbb S}^{N-1}}u(r,\theta)|^2
-a(\theta)|u(r,\theta)|^2\right]dS(\theta)
  \geq \mu_1(a)\int_{{\mathbb S}^{N-1}}\!|u(r,\theta)|^2dS(\theta).
\end{equation}
From (\ref{eq:4}), (\ref{eq:5}), and (\ref{eq:3}), we deduce that
\begin{align*}
  \int_{B_1}\! \bigg( |\nabla u(x)|^2-
  \frac{a(\frac{x}{|x|})}{|x|^2}|u(x)|^2\bigg)\,dx+ \frac{N-2}{2}\int_{\partial
    B_1}\!|u(x)|^2\,dS\geq
  \left[\bigg(\frac{N-2}2\bigg)^{\!\!2}\!\!+\mu_1(a)\right]
  \int_{B_1}\!\frac{|u(x)|^2}{|x|^2}\,dx
\end{align*}
for all $u\in
  C^\infty(B_1)\cap H^1(B_1)$  with $0\not\in {\mathop{\rm
      supp}}u$,
which, by density, yields the stated inequality for all
$H^1(B_r)$-functions for $r=1$.
\end{pf}

\begin{Corollary} \label{c:hardyboundary-0}
For all $r>0$ and $u\in H^1(B_r)$, there holds
  \begin{align}\label{eq:10}
    \int_{B_r} |\nabla u(x)|^2\,dx+
    \frac{N-2}{2r}\int_{\partial B_r}|u(x)|^2\,dS
\geq \bigg(\frac{k-2}{2}\bigg)^{\!\!2}
      \int_{B_r}\frac{|u(x)|^2}{|x_J|^2}\,dx
  \end{align}
for any $J\in \mathcal A_k$ and
  \begin{align}\label{eq:10-Mb}
    \int_{B_r} |\nabla u(x)|^2\,dx+
    \frac{N-2}{2r}\int_{\partial B_r}|u(x)|^2\,dS
\geq \frac{(k-2)^2}{2}
      \int_{B_r}\frac{|u(x)|^2}{|x_{J_1}-x_{J_2}|^2}\,dx
  \end{align}
 for any $(J_1,J_2)\in \mathcal B_k$.
\end{Corollary}
\begin{pf}
Let $r>0$ and $u\in H^1(B_r)$.
Choosing $a$ as in the Example \ref{ex} with
$\alpha<\big(\frac{k-2}{2}\big)^2$, from Lemma
\ref{l:hardyboundary}, it follows that
  \begin{align*}
    \int_{B_r} \bigg( |\nabla u(x)|^2-
\frac{\alpha}{|x_J|^2}|u(x)|^2\bigg)\,dx+
    \frac{N-2}{2r}\int_{\partial B_r}|u(x)|^2\,dS\geq 0
  \end{align*}
hence
  \begin{align*}
\alpha\int_{B_r}  \frac{|u(x)|^2}{|x_J|^2}\,dx\leq
  \int_{B_r} |\nabla u(x)|^2\,dx+
    \frac{N-2}{2r}\int_{\partial B_r}|u(x)|^2\,dS.
  \end{align*}
  Letting $\alpha\to \big(\frac{k-2}{2}\big)^2$, \eqref{eq:10} follows.
To prove  \eqref{eq:10-Mb}, we may choose $a$ as in Example
\ref{ex:2} and proceed as in the proof of \eqref{eq:10}.
  \end{pf}

\begin{Corollary}\label{c:hardyboundary}
Let $a$ be as in (\ref{eq:ateta}) and assume that
(\ref{eq:lambdamin1}) holds. Then, for all $r>0$, $u\in H^1(B_r)$,
$J\in \mathcal A_k$ and $(J_1,J_2)\in \mathcal B_k$, there holds
\begin{multline}\label{coercivity}
\int_{B_r} |\nabla u(x)|^2 \, dx-\int_{B_r}
\frac{a(\frac{x}{|x|})}{|x|^2} |u(x)|^2 \,
dx+\Lambda(a)\frac{N-2}{2r}\int_{\partial B_r} |u(x)|^2 \, dS\\
 \geq
(1-\Lambda(a)) \int_{B_r} |\nabla u(x)|^2 \, dx \, ,
\end{multline}
\begin{gather} \label{eq:11}
\int_{B_r} |\nabla u(x)|^2 \, dx-\int_{B_r}
\frac{a(\frac{x}{|x|})}{|x|^2} |u(x)|^2 \,
dx+\frac{N-2}{2r}\int_{\partial B_r} |u(x)|^2 \, dS \\
\notag \geq
(1-\Lambda(a)) \left(\frac{k-2}{2}\right)^2 \int_{B_r}
\frac{|u(x)|^2}{|x_J|^2} \, dx,
\end{gather}
and
\begin{gather} \label{eq:11-MB}
\int_{B_r} |\nabla u(x)|^2 \, dx-\int_{B_r}
\frac{a(\frac{x}{|x|})}{|x|^2} |u(x)|^2 \,
dx+\frac{N-2}{2r}\int_{\partial B_r} |u(x)|^2 \, dS \\
\notag \geq (1-\Lambda(a)) \frac{(k-2)^2}{2} \int_{B_r}
\frac{|u(x)|^2}{|x_{J_1}-x_{J_2}|^2} \, dx
\end{gather}
with $\Lambda(a)$ as in (\ref{eq:1}).
\end{Corollary}

\begin{pf}
By scaling, it is enough to prove the inequalities for $r=1$. Let
$u\in C^\infty(B_1)\cap H^1(B_1)$ with $0\notin \mathop{\rm supp}
u$. Passing in polar
coordinates we obtain
\begin{align}\label{eq:4-bis}
    \int_{B_1} \bigg( |\nabla u(x)|^2-&
\frac{a(\frac{x}{|x|})}{|x|^2}|u(x)|^2\bigg)\,dx+\Lambda(a)
    \frac{N-2}{2}\int_{\partial B_1}|u(x)|^2\,dS\\
  \notag=& \int_{{\mathbb S}^{N-1}}\bigg(
  \int_0^{1}r^{N-1}|\partial_ru(r,\theta)|^2\,dr\bigg)\,dS(\theta)
  +\Lambda(a)\frac{N-2}2\int_{\SN}|u(1,\theta)|^2\,dS(\theta) \\
  \notag& \quad+ \int_0^{1}\frac{r^{N-1}}{r^2}\bigg(\int_{{\mathbb
      S}^{N-1}}\left[|\nabla_{{\mathbb S}^{N-1}}u(r,\theta)|^2-a(\theta)
    |u(r,\theta)|^2\right]\,dS(\theta)\bigg)\,dr.
\end{align}
By \eqref{supremum} and \eqref{eq:lambdamin1} we have
\begin{align*}
\int_{\SN}  &\big[\big|\nabla_{\mathbb S^{N-1}} u(r,\theta)\big|^2
    - a(\theta)|u(r,\theta)|^2\big]\,dS(\theta) \\
&\geq (1-\Lambda(a))\int_{\SN} \big|\nabla_{\SN} u(r,\theta)\big|^2
\, dS(\theta) -\Lambda(a)\left(\frac{N-2}{2}\right)^2 \int_{\SN}
|u(r,\theta)|^2 \, dS(\theta)
\end{align*}
which inserted into \eqref{eq:4-bis} gives
\begin{align*} 
    & \int_{B_1} \bigg( |\nabla u(x)|^2-
\frac{a(\frac{x}{|x|})}{|x|^2}|u(x)|^2\bigg)\,dx+\Lambda(a)
    \frac{N-2}{2}\int_{\partial B_1}|u(x)|^2\,dS \geq (1-\Lambda(a)) \int_{B_1} |\nabla u(x)|^2 \, dx\\
&  +\Lambda(a)
\left[\int_{\SN}\bigg(\int_0^{1}r^{N-1}|\partial_ru(r,\theta)|^2\,dr
+\frac{N-2}2 |u(1,\theta)|^2\bigg)\,dS(\theta)
-\left(\frac{N-2}2\right)^2\int_{B_1} \frac{|u(x)|^2}{|x|^2} \, dx
\right].
\end{align*}
Now, inequality \eqref{coercivity} follows immediately from
\eqref{eq:5}.

From \eqref{coercivity} and \eqref{eq:10} we obtain
\begin{align*}
   & \int_{B_1}  \bigg( |\nabla u(x)|^2-
\frac{a(\frac{x}{|x|})}{|x|^2}|u(x)|^2\bigg)\,dx+
    \frac{N-2}{2}\int_{\partial B_1}|u(x)|^2\,dS \\
   & \geq (1-\Lambda(a)) \left(\int_{B_1} |\nabla u(x)|^2 \, dx
+\frac{N-2}2 \int_{\partial B_1} |u(x)|^2 \, dS\right)\geq
(1-\Lambda(a)) \left(\frac{k-2}2\right)^2 \int_{B_1}
\frac{|u(x)|^2}{|x_J|^2} \, dx
\end{align*}
for all $J\in \mathcal A_k$ and for all $u\in C^\infty(B_1)\cap
H^1(B_1)$ with $0\notin \mathop{\rm supp}u$.

On the other hand by \eqref{coercivity} and \eqref{eq:10-Mb} we
obtain
\begin{align*}
    \int_{B_1} & \bigg( |\nabla u(x)|^2-
\frac{a(\frac{x}{|x|})}{|x|^2}|u(x)|^2\bigg)\,dx+
    \frac{N-2}{2}\int_{\partial B_1}|u(x)|^2\,dS \\
   & \geq (1-\Lambda(a)) \left(\int_{B_1} |\nabla u(x)|^2 \, dx
+\frac{N-2}2 \int_{\partial B_1} |u(x)|^2 \, dS\right) \\
& \geq (1-\Lambda(a)) \frac{(k-2)^2}2 \int_{B_1}
\frac{|u(x)|^2}{|x_{J_1}-x_{J_2}|^2} \, dx
\end{align*}
for all $(J_1,J_2)\in \mathcal B_k$ and for all $u\in
C^\infty(B_1)\cap H^1(B_1)$ with $0\notin \mathop{\rm supp}u$.

By density the stated inequalities follow for any $u\in
H^1(B_1)$.
\end{pf}

From (\ref{eq:10bi}) and (\ref{coercivity}), we can derive a
Hardy-Sobolev type inequality which takes into account the boundary
terms; to this aim, the following lemma is needed.

\begin{Lemma}\label{l:se}
Let  $\widetilde S_N>0$ be the best constant of the Sobolev embedding
$H^1(B_1)\subset L^{2^*}(B_1)$, i.e.
\begin{equation}\label{eq:9}
  \widetilde S_N:=
  \inf_{v\in H^1(B_1)\setminus\{0\}}\dfrac{\int_{B_1} \big(|\nabla u(x)|^2+
    |u(x)|^2\big)\,dx}{\Big(\int_{B_1}
    |u(x)|^{2^*}dx\Big)^{\!2/2^*}}.
\end{equation}
Then, for every $r>0$ and $u\in
  H^1(B_r)$, there holds
\begin{align}\label{eq:6}
  \int_{B_r} \bigg(|\nabla u(x)|^2+ \frac{|u(x)|^2}{|x|^2}\bigg)\,dx
\geq \widetilde S_N
  \bigg(\int_{B_r}|u(x)|^{2^*}dx\bigg)^{\!\!2/2^*}.
  \end{align}
  \end{Lemma}
\begin{pf}
  Inequality (\ref{eq:6}) follows simply by scaling from the
  definition of $\widetilde S_N$.
\end{pf}

The following boundary Hardy-Sobolev inequality holds true.
\begin{Corollary}\label{c:hardyboundarySOB}
Let $a$ be as in (\ref{eq:ateta}) and assume that
(\ref{eq:lambdamin1}) holds. Then, for all $r>0$ and $u\in H^1(B_r)$,
 there holds
\begin{multline}\label{coercivitySOB}
\int_{B_r} |\nabla u(x)|^2 \, dx-\int_{B_r}
\frac{a(\frac{x}{|x|})}{|x|^2} |u(x)|^2 \,
dx+\frac{1+\Lambda(a)}2\,\frac{N-2}{2r}\int_{\partial B_r} |u(x)|^2 \, dS\\
 \geq
\frac{\widetilde S_N}2\min\bigg\{
1-\Lambda(a),\mu_1(a)+\bigg(\frac{N-2}{2}\bigg)^{\!\!2}\bigg\}
 \bigg(\int_{B_r}|u(x)|^{2^*}dx\bigg)^{\!\!2/2^*},
\end{multline}
where $\widetilde S_N$ is defined in (\ref{eq:9}).
\end{Corollary}
\begin{pf}
Inequality (\ref{coercivitySOB}) follows simply by summing up (\ref{eq:10bi})
and (\ref{coercivity}) and using Lemma \ref{l:se}.
\end{pf}

\section{A Pohozaev-type identity}\label{sec:pohoz-type-ident}

\noindent

In order to approximate $L_a:=-\Delta_{\mathbb S^{N-1}}- a$ with
operators with bounded coefficients, for all $\lambda\in \R$, we define

\begin{equation}\label{eq:18}
a_\lambda(\theta):=
\left\{
\begin{array}{ll}
{\displaystyle \sum_{J\in{\mathcal A}_k} \frac{\alpha_J}{|\theta_J|^2+\lambda}+
\sum_{(J_1,J_2)\in{\mathcal B}_k}\frac{\alpha_{J_1\, J_2}}{|\theta_{J_1}-\theta_{J_2}|^2+\lambda}}
& \qquad \text{if } \lambda>0 \\ \\
a(\theta) & \qquad \text{if } \lambda\leq 0
\end{array}
\right.
\end{equation}
in such a way that $a_\lambda\in L^\infty({\mathbb S}^{N-1})$ for any $\lambda>0$.
We notice that $a_\lambda\in{\mathcal F}$ for any $\lambda\in \R$.

Since we are interested in the asymptotics of solutions at $0$, we
focus our attention on a ball $B_{r_0}$ which is sufficiently small to
ensure positivity of the quadratic forms associated to equation
(\ref{u}) and to some proper approximations of (\ref{u}) in $B_{r_0}$.
Let $u$ be a solution of (\ref{u}), with the perturbation potential
$h$ satisfying \eqref{H1-2} and the nonlinear term $f$ satisfying
\eqref{F}.  If condition (\ref{eq:lambdamin1}) holds, there exists
$r_0>0$ such that
\begin{multline}\label{eq:r_0}
B_{r_0}\subseteq\Omega\quad\text{and}\quad \Lambda(a)
+C_h
  r_0^\e\binom{N}{k}
\bigg(\frac{2}{k-2}\bigg)^{\!\!2}\bigg(1+\binom{N-k}{k}\bigg)\\
+C_f S^{-1}\left[(\omega_{N-1}/N)^{\frac2N}r_0^{2}+
 \|u\|_{L^{2^*}(B_{r_0})}^{2^*-2}\right]
<1,
\end{multline}
with $a$ as (\ref{eq:ateta}), $\Lambda(a)$ as in (\ref{supremum}) and
$\binom{N-k}{k}=0$ whenever $N<2k$.

\begin{Lemma}\label{l:approx}
  Let $\Omega\subset\R^N$, $N\geq 3$, be a bounded open set such that
  $0\in\Omega$, and let $a$ satisfy (\ref{eq:ateta}) and
  (\ref{eq:lambdamin1}).  Suppose that $h$ satisfies \eqref{H1-2},
$f$ satisfies \eqref{F}, $u$ is a $H^1(\Omega)$-weak solution
to (\ref{u}) in $\Omega$, and  $r_0>0$ is as in (\ref{eq:r_0}).
Then there exists $\bar\lambda>0$ such that, for every
$\lambda \in (0,\bar \lambda)$, the Dirichlet boundary value problem
\begin{equation} \label{uenne}
\begin{cases}
  -\Delta v(x)- \dfrac{a_\lambda\big(\frac{x}{|x|}\big)}{|x|^2}\,v(x)
  =h_\lambda(x)v(x)+f(x,v(x)), \quad
  &\text{in} \ B_{r_0},\\[10pt]
  v\big|_{\partial B_{r_0}}=u\big|_{\partial B_{r_0}},\quad &\text{on}
  \ \partial B_{r_0},
\end{cases}
\end{equation}
with
$$
h_\lambda(x)=
\begin{cases} \min\{1/\lambda,
\max\{-1/\lambda,h(x)\}\}, \quad&\text{if } \lambda>0, \\
h(x), \quad&\text{if } \lambda\leq 0,
\end{cases}
$$
admits a weak solution $u_\lambda\in H^1(B_{r_0})$ such that
$$
u_\lambda\to u \quad\text{in }H^1(B_{r_0})\quad\text{as } \lambda\to 0^+.
$$
\end{Lemma}
\begin{pf}
Let $\tilde v$ be the
unique $H^1(B_{r_0})$-weak solution to the problem
\begin{equation} \label{Lax}
\begin{cases}
  -\Delta \tilde v-\dfrac{a\big(\frac{x}{|x|}\big)}{|x|^2}\,
  \tilde v(x)=h(x) \tilde v,&\text{in }B_{r_0},\\[10pt]
  \tilde v=u\big|_{\partial B_{r_0}},&\text{on }\partial B_{r_0}.
\end{cases}
\end{equation}
The existence and uniqueness of such a $\tilde v$ can be proven by
introducing the continuous bilinear form ${\mathcal Q}:H^1_0(B_{r_0})\times
H^1_0(B_{r_0})\to\R$
\begin{align*}
  {\mathcal Q}(w_1,w_2):=\int_{B_{r_0}}\bigg[ \nabla w_1(x)\cdot\nabla
  w_2(x)-\bigg(
\frac{a(\frac{x}{|x|})}{|x|^2}+ h(x)\bigg) w_1(x)w_2(x)\bigg]\,dx,
\end{align*}
and the continuous functional $\Psi\in H^{-1}(B_{r_0})$
\begin{equation*}
{\phantom{\Big\langle}}_{H^{-1}(B_{r_0})}
\!\big\langle \Psi,w\big\rangle_{H^1_0(B_{r_0})}=-\int_{B_{r_0}} \!\! \nabla u(x)\cdot \nabla w(x)
dx+\int_{B_{r_0}} \frac{a(\frac{x}{|x|})}{|x|^2} u(x)w(x) dx
+\int_{B_{r_0}} h(x)u(x)w(x)  dx \, .
\end{equation*}
By \eqref{H1-2}, (\ref{eq:ssw}), (\ref{eq:MPHardy}), and
\eqref{eq:bound}, we have
\begin{multline} \label{eq:coerc}
\mathcal Q(w,w)=\int_{B_{r_0}}\bigg(|\nabla
  w(x)|^2-\frac{a(\frac{x}{|x|})}{|x|^2}w^2(x)-
  h(x)w^2(x)\bigg)\,dx\\
  \geq \int_{B_{r_0}}\bigg(|\nabla
  w(x)|^2-\frac{a(\frac{x}{|x|})}{|x|^2}w^2(x)-
  C_h \bigg( \sum_{J\in{\mathcal A}_k}|x_{J}|^{-2+\e}+
\sum_{(J_1,J_2)\in{\mathcal B}_k}|x_{J_1}-x_{J_2}|^{-2+\e}
  \bigg) \,w^2(x)\bigg)\,dx\\
 \geq
  \bigg[1-\Lambda(a)-C_hr_0^\e\binom{N}{k}
\left(\frac{2}{k-2}\right)^2 \left(1+\binom{N-k}{k}\right)\bigg]
\int_{B_{r_0}}\!\!|\nabla w(x)|^2\,dx
\end{multline}
for all $w\in H^1_0(B_{r_0})$.
By \eqref{eq:coerc}, \eqref{eq:lambdamin1} and \eqref{eq:r_0} it
follows that the bilinear form $\mathcal Q$ is coercive. The
Lax-Milgram lemma yields existence and uniqueness of a solution
$v\in H^1_0(B_{r_0})$ of the variational problem
$$
\mathcal Q(v,w)={\phantom{\Big\langle}}_{H^{-1}(B_{r_0})}
\!\big\langle \Psi,w\big\rangle_{H^1_0(B_{r_0})} \qquad \text{for any } w\in
H^1_0(B_{r_0}).
$$
Then the function $\tilde v:=v+u$ is the unique solution of
\eqref{Lax}.

Let us now define the map $\Phi:\R\times H^1_0(B_{r_0})\to H^{-1}(B_{r_0})$ as
\begin{gather*}
\notag  \Phi(\lambda,w)=-\Delta
w-\frac{a_\lambda\big(\frac{x}{|x|}\big)}{|x|^2}\, w -h_\lambda(x)
w-f(x,\tilde
v+w)+\left(\frac{a\big(\frac{x}{|x|}\big)}{|x|^2}+h(x)
-\frac{a_\lambda\big(\frac{x}{|x|}\big)}{|x|^2}-h_\lambda(x)\right)
\tilde v .
\end{gather*}
By \eqref{eq:ateta}, \eqref{eq:ssw}, \eqref{eq:MPHardy}, \eqref{H1-2}
and \eqref{F}, the function $\Phi$ is continuous and its first
variation with respect to the $w$ variable
$$
\Phi'_w:\R\times H^1_0(B_{r_0})\to \mathcal
L(H^1_0(B_{r_0}),H^{-1}(B_{r_0}))
$$
is also continuous.
We claim that
$$
\Phi(0,u-\tilde v)=0\text{  in }H^{-1}(B_{r_0})\quad\text{and}\quad
\Phi'_w(0,u-\tilde v)\in \mathcal
L\big(H^1_0(B_{r_0}),H^{-1}(B_{r_0})\big)\text{ is an isomorphism}.
$$
The first claim is an immediate consequence of the definition of $u$
and $\tilde v$. Let us prove the second one.  By \eqref{F},
\eqref{eq:bound}, and H\"older and Sobolev inequalities, for every
$w\in H^1_0(B_{r_0})$ we obtain
\begin{align*}
 &\hskip-20pt{\phantom{\Big\langle}}_{H^{-1}(B_{r_0})}
\!\Big\langle \Phi'_w(0,u-\tilde v)w,w\Big\rangle_{H^1_0(B_{r_0})} \\
& =
 \int_{B_{r_0}} |\nabla w(x)|^2\, dx-\int_{B_{r_0}}
\frac{a\big(\frac{x}{|x|}\big)}{|x|^2} w^2(x)\, dx-\int_{B_{r_0}}
h(x) w^2(x)\, dx
-\int_{B_{r_0}} f'_s(x,u(x)) w^2(x)\, dx \\
& \geq \int_{B_{r_0}} |\nabla w(x)|^2\, dx-\int_{B_{r_0}}
\frac{a\big(\frac{x}{|x|}\big)}{|x|^2} w^2(x)\, dx
-\int_{B_{r_0}} h(x) w^2(x)\, dx \\
& \qquad \quad -C_f
\int_{B_{r_0}} \big(1+|u(x)|^{2^*-2}\big) w^2(x) \, dx \\
& \geq (1-\Lambda(a))\int_{B_{r_0}} |\nabla w(x)|^2\,
dx \\
& \qquad \quad -C_h r_0^\e\binom{N}{k}
\bigg(\frac{2}{k-2}\bigg)^{\!\!2}\bigg(1+\binom{N-k}{k}\bigg)
\int_{B_{r_0}} |\nabla w(x)|^2 \, dx
\\
& \qquad \quad -C_f S^{-1}
\left[(\omega_{N-1}/N)^{\frac2N}r_0^{2}+
 \|u\|_{L^{2^*}(B_{r_0})}^{2^*-2}\right]
\int_{B_{r_0}} |\nabla w(x)|^2 \, dx \, .
\end{align*}
The above estimate, together with \eqref{eq:r_0}, shows that the
quadratic form $w\mapsto \langle \Phi'_w(0,u-\tilde v)w,w\rangle$ is
positive definite over $H^1_0(B_{r_0})$. Then the Lax-Milgram lemma
applied to the continuous and coercive bilinear form $(w_1,w_2)\mapsto
\!\!\phantom{.}_{H^{-1}(B_{r_0})}\big\langle \Phi'_w(0,u-\tilde
v)w_1,w_2\big\rangle_{H^1_0(B_{r_0})}$ ensures that the operator
$\Phi'_w(0,u-\tilde v)\in \mathcal L(H^1_0(B_{r_0}),H^{-1}(B_{r_0}))$
is an isomorphism and hence our second claim is proved.

We are now in position to apply the Implicit Function Theorem to
the map $\Phi$, thus showing the existence of $\bar\lambda>0$,
$\rho>0$, and of a continuous function
$$
g:(-\bar\lambda,\bar\lambda)\to B(u-\tilde v,\rho)
$$
with $B(u-\tilde v,\rho)=\{w\in H^1_0(B_{r_0}): \|w-u+\tilde
v\|_{H^1_0(B_{r_0})}<\rho \}$, such that
$\Phi(\lambda,g(\lambda))=0$ for all $\lambda\in
  (-\bar\lambda,\bar\lambda)$
and, if $(\lambda,w)\in
(-\bar\lambda,\bar\lambda)\times B(u-\tilde v,\rho)$ and
$\Phi(\lambda,w)=0$, then $w=g(\lambda)$. The function
$u_\lambda:=g(\lambda)+\tilde v$  solves
\eqref{uenne} for any $\lambda\in (0,\bar\lambda)$. Moreover, by
the continuity of $g$ over the interval
$(-\bar\lambda,\bar\lambda)$ and the fact that $g(0)=u-\tilde v$,
$u_\lambda-u=g(\lambda)-u+\tilde v\to 0$ in $H^1_0(B_{r_0})$ as
$\lambda\to 0^+$. This proves that $u_\lambda\to u$ in
$H^1(B_{r_0})$ as $\lambda\to 0^+$.
\end{pf}

\begin{remark}\label{rem:l1w11}
  We notice that, if $f\in L^1(\Omega)$ for some $\Omega\subset\R^N$
  bounded open set such that $0\in\Omega$, then, for every $r>0$ such
  that $B_{r}\subseteq\Omega$,
\begin{align*}
  \int_{B_{r}}|f(x)|\,dx=\int_0^{r}\bigg(\int_{\partial
    B_s}|f|\,dS\bigg)ds<+\infty,
\end{align*}
and hence the function $s\mapsto \int_{\partial B_s}|f|\,dS$ belongs
to $L^1(0,r)$ and is the weak derivative of the
$W^{1,1}(0,r)$-function $s\to \int_{B_{s}}|f(x)|\,dx$. In particular,
for every $u\in H^1(\Omega)$ and every $J\in \mathcal A_k$,
$(J_1,J_2)\in{\mathcal B}_k$,
 the
$L^1(0,r)$-function
$$
s\mapsto
\int_{\partial B_s}|\nabla u(x)|^2\,dS, \quad\text{respectively }s\mapsto
\int_{\partial B_s}\frac{u^2(x)}{|x_J|^{2}}\,dS,
\quad
s\mapsto
\int_{\partial B_s}\frac{u^2(x)}{|x_{J_1}-x_{J_2}|^{2}}\,dS,
$$
  is the weak derivative of the
$W^{1,1}(0,r)$-function
$$
s\to \int_{B_{s}}|\nabla u(x)|^2\,dx, \quad\text{respectively
}s\mapsto \int_{B_{s}}\frac{u^2(x)}{|x_J|^{2}}\,dx,
\quad
s\mapsto
\int_{B_s}\frac{u^2(x)}{|x_{J_1}-x_{J_2}|^{2}}\,dx.
$$
\end{remark}

\noindent Solutions to (\ref{u}) satisfy the following Pohozaev-type identity.

\begin{Theorem} \label{t:pohozaev} Let $\Omega\subset\R^N$, $N\geq 3$,
  be a bounded open set such that $0\in\Omega$. Let $a$ satisfy
  (\ref{eq:ateta}), (\ref{eq:lambdamin1}), and
$u$ be a $H^1(\Omega)$-weak solution to (\ref{u}) in $\Omega$
with $h$ satisfying
  \eqref{H1-2} and $f$ satisfying \eqref{F}. Then
\begin{multline}\label{eq:poho}
  -\frac{N-2}2\int_{B_r}\bigg[ |\nabla u(x)|^2
  -\frac{a(\frac{x}{|x|})}{|x|^2}u^2(x)\bigg]\,dx+\frac{r}{2}\int_{\partial
    B_r}\bigg[ |\nabla u(x)|^2
  -\frac{a(\frac{x}{|x|})}{|x|^2}u^2(x)\bigg]\,dS\\
  =r\int_{\partial B_r}\bigg|\frac{\partial u}{\partial
    \nu}\bigg|^2\,dS
  -\frac{1}2\int_{B_r}(\nabla h(x)\cdot x)u^2(x)\,dx
  -\frac{N}2\int_{B_r}h(x)\,u^2(x)\,dx
+
  \frac{r}2\int_{\partial B_r}h(x)u^2(x)\,dS  \\
+r\int_{\partial B_r} F(x,u(x))\, dS-\int_{B_r}
[\nabla_xF(x,u(x))\cdot x+NF(x,u(x))]\, dx
\end{multline}
and
\begin{multline}\label{eq:poho2}
  \int_{B_r}\bigg(|\nabla
  u(x)|^2-\frac{a(\frac{x}{|x|})}{|x|^2}u^2(x)\bigg)\,dx\\
  =\int_{\partial B_r}u\frac{\partial
    u}{\partial\nu}\,dS+\int_{B_r}h(x)u^2(x)\,dx+
    \int_{B_r} f(x,u(x))u(x)\, dx,
  \end{multline}
  for a.e. $r\in(0, r_0)$, where $r_0>0$ satisfies (\ref{eq:r_0}) and
$\nu=\nu(x)$ is the unit outer normal vector $\nu(x)=\frac{x}{|x|}$.
\end{Theorem}

\begin{pf}
  Let $a_\lambda$ as in (\ref{eq:18}), $r_0$ as in (\ref{eq:r_0}), and 
 $u_\lambda$,  $h_\lambda$  as
  in Lemma \ref{l:approx}. Since $a_\lambda$ and $h_\lambda$ are
  bounded for every $\lambda>0$ the following Pohozaev identity
\begin{multline}\label{eq:pohoapp}
  -\frac{N-2}2\int_{B_r}\bigg[ |\nabla u_\lambda(x)|^2
  -\frac{a_\lambda(\frac{x}{|x|})}
  {|x|^2}u_\lambda^2(x)\bigg]\,dx
  +\frac{r}{2}\int_{\partial B_r}\bigg[ |\nabla u_\lambda(x)|^2
  -\frac{a_\lambda(\frac{x}{|x|})}{|x|^2}u_\lambda^2(x)\bigg]\,dS\\
  =r\int_{\partial B_r}\bigg|\frac{\partial u_\lambda}{\partial
    \nu}\bigg|^2\,dS+ \int_{B_r}h_\lambda(x)u_\lambda(x)\,(x\cdot {\nabla
    u_\lambda(x)})\big)\,dx \\
+r\int_{\partial B_r} F(x,u_\lambda(x))\, dx-\int_{B_r}
[\nabla_xF(x,u_\lambda(x))\cdot x+NF(x,u_\lambda(x))]\, dx
\end{multline}
holds for all $r\in(0, r_0)$, see 
Proposition \ref{p:poho}.
Furthermore, testing (\ref{uenne}) with $u_{\lambda}$, integrating by parts,
and using the regularity of $u_{\lambda}$ outside the origin, we obtain that
\begin{multline}\label{eq:poho2app}
  \int_{B_r}\bigg( |\nabla u_{\lambda}(x)|^2
  -\frac{a_\lambda(\frac{x}{|x|})}
  {|x|^2}u_\lambda^2(x)\bigg)\,dx\\
=\int_{\partial B_r}u_{\lambda}
\frac{\partial u_{\lambda}}{\partial
    \nu}\,dS
+\int_{B_r}h_\lambda(x)u_{\lambda}^2(x)\,dx
+\int_{B_r} f(x,u_\lambda(x))u_\lambda(x)\, dx
\end{multline}
for all $r\in(0, r_0)$.

From the convergence of $u_{\lambda}$ to $u$ in $H^1(B_{r_0})$ as
$\lambda\to 0^+$ proved in Lemma \ref{l:approx}, inequalities
(\ref{eq:10}--\ref{eq:10-Mb}), and the Dominated Convergence
Theorem, it follows that
$$
\frac{a_\lambda(\frac{x}{|x|})}
  {|x|^2}u_{\lambda}^2-\frac{a(\frac{x}{|x|})}
  {|x|^2}u^2=
\frac{a_\lambda(\frac{x}{|x|})}
  {|x|^2}(u_{\lambda}+u)(u_{\lambda}-u)
+\frac{a_\lambda(\frac{x}{|x|})-a(\frac{x}{|x|})}
  {|x|^2}u^2
\to 0
$$
in $L^1(B_{r_0})$ as $\lambda\to 0^+$,
i.e.
\begin{multline}\label{eq:anlim}
  \lim_{\lambda\to 0^+}\int_{B_{r_0}}\bigg|
  \frac{a_\lambda(\frac{x}{|x|})}
  {|x|^2}u_{\lambda}^2(x)-\frac{a(\frac{x}{|x|})}
  {|x|^2}u^2(x)\bigg|dx\\
  = \lim_{\lambda\to0^+}\int_0^{r_0}\bigg[
\int_{\partial B_s}\bigg|
  \frac{a_\lambda(\frac{x}{|x|})}
  {|x|^2}u_{\lambda}^2(x)-\frac{a(\frac{x}{|x|})}
  {|x|^2}u^2(x)\bigg|
dS \bigg]ds=0.
\end{multline}
From (\ref{eq:anlim}) we deduce that
$$
\int_{B_r}\frac{a_\lambda(\frac{x}{|x|})}
  {|x|^2}u_{\lambda}^2(x)\,dx
\to\int_{B_r}\frac{a(\frac{x}{|x|})}
  {|x|^2}u^2(x)\,dx \quad\text{as } \lambda\to 0^+\quad\text{for all }r\in(0,r_0).
$$
and, along a sequence $\lambda_n\to 0^+$,
\begin{equation} \label{seq1}
\int_{\partial B_r}\frac{a_{\lambda_n}(\frac{x}{|x|})}
  {|x|^2}u_{\lambda_n}^2\,dS
\to\int_{\partial B_r}\frac{a(\frac{x}{|x|})}
  {|x|^2}u^2\,dS \quad\text{as } n\to +\infty\quad\text{for a.e. }r\in(0,r_0).
\end{equation}
On the other hand, from
\begin{align*}
  \lim_{\lambda\to 0^+}\int_{B_{r_0}}|\nabla (u_{\lambda}-u)(x)|^2\,dx =
  \lim_{\lambda\to 0^+}\int_0^{r_0}\bigg[\int_{\partial B_s}|\nabla (u_{\lambda}-u)|^2\,dS
  \bigg]ds=0,
\end{align*}
we deduce that, along a sequence converging monotonically to zero still
denoted by $\lambda_n$,
\begin{equation} \label{seq2}
\int_{\partial B_r}
  |\nabla u_{\lambda_n}|^2dS\to \int_{\partial B_r}
  |\nabla u|^2dS\quad\text{as } n\to +\infty\quad\text{for a.e. }r\in(0,r_0)
\end{equation}
and
\begin{equation} \label{seq2-1}
\int_{\partial B_r}
  \bigg|\frac{\partial u_{\lambda_n}}{\partial
    \nu}\bigg|^2dS\to \int_{\partial B_r}
  \bigg|\frac{\partial u}{\partial
    \nu}\bigg|^2dS\quad\text{as } n\to +\infty\quad\text{for a.e. }r\in(0,r_0).
\end{equation}
Let us fix $\lambda>0$ and $r\in (0,r_0)$. Since
$$
  \int_0^r\bigg[\int_{\partial B_s} h_\lambda(x)|u_{\lambda}(x)|^2\,dS\bigg]ds
= \int_{B_r}h_\lambda(x)|u_\lambda(x)|^2\,dx<+\infty,
$$
there exists a sequence $\{\delta_k\}_{k\in\N}\subset(0,r)$ such that
$\lim_{k\to+\infty}\delta_k=0$ and
\begin{equation}\label{eq:33cyl}
\delta_k \int_{\partial B_{\delta_k}}h_\lambda(x)|u_{\lambda}(x)|^2
\,dS\longrightarrow 0\quad \text{as }k\to+\infty.
\end{equation}
Recalling that $u_{\lambda}\in C^{1,\tau}_{\rm loc}(B_{r_0}\setminus\{0\})$
for some $\tau\in(0,1)$, integration by parts yields
\begin{align*}
\int_{B_r\setminus B_{\delta_k}}&h_\lambda(x)u_{\lambda}(x)\,(x\cdot {\nabla
    u_{\lambda}(x)})\,dx=-\frac12\int_{B_r\setminus B_{\delta_k}}
(\nabla h_\lambda(x)\cdot x)u_{\lambda}^2(x)\,dx\\
&-\frac{N}2\int_{B_r\setminus B_{\delta_k}}h_\lambda(x)
u_{\lambda}^2(x)\,dx+\frac{r}2\int_{\partial B_r}h_\lambda(x)
u_{\lambda}^2(x)\,dS-\frac{\delta_k}2\int_{\partial B_{\delta_k}}h_\lambda(x)
u_{\lambda}^2(x)\,dS.
\end{align*}
Letting $k\to+\infty$, by \eqref{eq:33cyl} and \eqref{H1-2}, we obtain
\begin{align*}
\int_{B_r}&h_\lambda(x)u_{\lambda}(x)\,(x\cdot {\nabla
    u_{\lambda}(x)})\,dx=-\frac12\int_{B_r}
(\nabla h_\lambda(x)\cdot x)u_{\lambda}^2(x)\,dx\\
&-\frac{N}2\int_{B_r}h_\lambda(x)
u_{\lambda}^2(x)\,dx+\frac{r}2\int_{\partial B_r}h_\lambda(x)
u_{\lambda}^2(x)\,dS.
\end{align*}
Arguing as above, using \eqref{H1-2} we can prove that
\begin{align*}
  &\lim_{\lambda\to 0^+}
    \int_{B_r}(\nabla h_\lambda(x)\cdot x)u_{\lambda}^2(x)\,dx
 =\int_{B_r}(\nabla h(x)\cdot x)u^2(x)\,dx
  \quad\text{for all }r\in(0,r_0),\\
  &\lim_{\lambda\to 0^+}
\int_{B_r} h_\lambda(x)u_{\lambda}^2(x)\,dx=
\int_{B_r} h(x)u^2(x)\,dx
  \quad\text{for all }r\in(0,r_0),
\end{align*}
and, along a sequence $\lambda_n\to 0^+$,
\begin{equation} \label{seq3}
\lim_{n\to +\infty}
\int_{\partial B_r}h_{\lambda_n}(x)u_{\lambda_n}^2(x)\,dS=
\int_{\partial B_r}h(x)u^2(x)\,dS\quad\text{for a.e. }r\in(0,r_0).
\end{equation}
It remains to study the convergence of the terms in \eqref{eq:pohoapp}
and \eqref{eq:poho2app} related to the nonlinearity $f$.  By
\eqref{F}, convergence of $u_{\lambda}$ to $u$ in $H^1(B_{r_0})$, and
the Dominated Convergence Theorem, we have that
\begin{equation*}
\lim_{\lambda\to 0^+}
\int_{B_r} [\nabla_x F(x,u_\lambda(x))\cdot x+N F(x,u_\lambda(x))] \, dx
=\int_{B_r} [\nabla_x F(x,u(x))\cdot x+N F(x,u(x))] \, dx \, ,
\end{equation*}
\begin{equation*}
\lim_{\lambda\to 0^+}
\int_{B_r}
f(x,u_\lambda(x))u_\lambda(x) \, dx=\int_{B_r}
f(x,u(x))u(x) \, dx,
\end{equation*}
for all $r\in (0,r_0)$, and along a sequence $\lambda_n\to0^+$,
\begin{equation} \label{seq4}
r\int_{\partial B_r} F(x,u_{\lambda_n}(x))\, dx\to
r\int_{\partial B_r} F(x,u(x))\, dx  \quad \text{as } n\to +\infty
\end{equation}
for a.e. $r\in(0,r_0)$.

Therefore, we can pass to the limit in \eqref{eq:pohoapp} and in
(\ref{eq:poho2app}) along a sequence $\lambda_n\to0^+$ such that
\eqref{seq1}, \eqref{seq2}, (\ref{seq2-1}), \eqref{seq3}, and
\eqref{seq4} hold true, thus obtaining \eqref{eq:poho} and
(\ref{eq:poho2}).
\end{pf}

\section{A Brezis-Kato type estimate}\label{sec:brezis-kato-type}

Throughout this section, we let 
$\Omega\subset \R^N$, $N\geq 3$, be a bounded open set such that
$0\in\Omega$, $a$ satisfy \eqref{eq:ateta}, \eqref{eq:lambdamin1}, $h$
satisfy \eqref{H1-2}, and $V\in L^1_{\rm
  loc}(\Omega)$  satisfy the form-bounded condition
$$
\sup_{v\in
  H^1(\Omega)\setminus\{0\}}\frac{\int_\Omega|V(x)|v^2(x)\,dx}
{\|v\|^2_{H^1(\Omega)}}<+\infty,
$$
see \cite{Mazja}. The above condition (which is in particular
satisfied by $L^{N/2}$ and $L^{N/2,\infty}$ functions, potentials of
the form (\ref{eq:V(x)}), etc.) in particular implies that for every
$u\in H^1(\Omega)$, $Vu\in H^{-1}(\Omega)$.  We assume that $u\in
H^1(\Omega)$ is a weak solution to
\begin{equation} \label{eq:ahV}
-\Delta u(x)-\frac{a(\frac{x}{|x|})}{|x|^2} u(x)=h(x)u(x)+V(x)u(x),
\quad \text{in } \Omega.
\end{equation}
  In the spirit of \cite[Theorem
2.3]{SM}, we prove the following Brezis-Kato type result.

\begin{Proposition} \label{SMETS}
If $V_+\in
L^{N/2}(\Omega)$, letting
$$
q_{\rm lim}:=\left\{
\begin{array}{ll}
\frac{2^*}{2} \min\left\{\frac{4}{\Lambda(a)}-2,2^*\right\}, &
\qquad \text{if }\quad \Lambda(a)>0, \\
\frac{(2^*)^2}2, & \qquad \text{if } \quad \Lambda(a)=0 ,
\end{array}
\right.
$$
then for every $1\leq q<q_{\rm lim}$ there exists $r_q>0$ depending
on $q,N,k,a,h$ such that $B_{r_q}\subset \Omega$ and
$u\in L^q(B_{r_q})$.
\end{Proposition}

\begin{pf}
For any $2<\tau<\frac{2}{2^*} q_{\rm lim}$ define
$C(\tau):=\frac{4}{\tau+2}$ and let $\ell_\tau>0$ be
large enough so that
\begin{equation} \label{elleq}
\bigg(\int\limits_{V_+(x) \geq\ell_\tau}\!\!\!\!\!\!
V_+^{\frac{N}{2}}\!(x) \, dx
\bigg)^{\!\!\frac{2}{N}}<\frac{S(C(\tau)-\Lambda(a))}2.
\end{equation}
Let $r>0$ be such that $B_r\subset \Omega$.  For any $w\in
H^1_0(B_{r})$, by H\"older and Sobolev inequalities and
\eqref{elleq}, we have
\begin{align} \label{estV}
\int_\Omega
V(x)|w(x)|^2\,dx &
\leq \int \limits_{B_r\cap
\{V_+(x) \leq\ell_\tau\}}\!\!\!\!\!\!\ V_+(x)|w(x)|^2 \, dx + \int
\limits_{B_r\cap \{V_+(x) \geq\ell_\tau\}}\!\!\!\!\!\!\
V_+(x)|w(x)|^2 \, dx
\\
\notag & \leq \ell_\tau \int_{B_r} |w(x)|^2 \, dx +
\bigg(\int\limits_{V_+(x) \geq\ell_\tau}\!\!\!\!\!\!
V_+^{\frac{N}{2}}(x) \, dx \bigg)^{\!\!\frac{2}{N}} \left(\int_{B_r}
|w(x)|^{2^*}\, dx\right)^{\!\!\frac{2}{2^*}}
\\
\notag & \leq \ell_\tau \int_{B_r} |w(x)|^2 \,
dx+\frac{C(\tau)-\Lambda(a)}2 \int_{B_r} |\nabla w(x)|^2 \, dx .
\end{align}
Let $\eta\in C^\infty_c (B_r)$ be such that $\eta\equiv 1$ in
$B_{r/2}$ and define $v(x):=\eta(x)u(x)\in H^1_0(B_r)$. Then $v$
is a $H^1(\Omega)$-weak solution of the equation
\begin{equation} \label{eq:g}
-\Delta v(x)-\frac{a(\frac{x}{|x|})}{|x|^2}
v(x)=h(x)v(x)+V(x)v(x)+g(x), \quad \text{in } \Omega,
\end{equation}
where $g(x)=-u(x)\Delta \eta(x)-2\nabla u(x)\cdot \nabla\eta(x)\in
L^2(B_r)$.
For any $n\in\N$, $n\geq 1$, let us define the function
$v^n:=\min\{|v|,n\}$. Testing \eqref{eq:g} with $(v^n)^{\tau-2} v\in
H^1_0(B_r)$ we obtain
\begin{multline} \label{vn}
\int_{B_r}  (v^n(x))^{\tau-2} |\nabla v(x)|^2 \, dx +(\tau-2)
\int_{B_r} (v^n(x))^{\tau-3}|v(x)| |\nabla v(x)|^2
\chi_{\{|v(x)|<n\}}(x) \, dx \\
 -\int_{B_r} \frac{a(\frac{x}{|x|})}{|x|^2} (v^n(x))^{\tau-2}
v^2(x)\,
dx \\
 \ =\int_{B_r} h(x)(v^n(x))^{\tau-2} v^2(x)\, dx
+\int_{B_r} V(x)(v^n(x))^{\tau-2} v^2(x)\, dx + \int_{B_r}
g(x)(v^n(x))^{\tau-2} v(x)\, dx \ .
\end{multline}
Since
$$
|\nabla((v^n)^{\frac{\tau}{2}-1}v)|^2= (v^n)^{\tau-2} |\nabla
v|^2+\frac{(\tau-2)(\tau+2)}4 (v^n)^{\tau-3} |v| |\nabla v|^2
\chi_{\{|v(x)|<n\}} \ ,
$$
then by \eqref{vn}, \eqref{eq:bound}, \eqref{H1-2}, and
\eqref{estV} with $w=(v^n)^{\frac{\tau}{2}-1}v$ we obtain
\begin{gather} \label{vn2}
C(\tau)  \int_{B_r} |\nabla((v^n(x))^{\frac{\tau}{2}-1}v(x))|^2 \, dx \\
\notag  \leq \int_{B_r} \frac{a(\frac{x}{|x|})}{|x|^2}
((v^n(x))^{\frac{\tau}{2}-1}v(x))^2 \, dx+ \int_{B_r}
h(x)((v^n(x))^{\frac{\tau}{2}-1}v(x))^2\,
dx  \\
\notag   \quad +\int_{B_r} V(x)((v^n(x))^{\frac{\tau}{2}-1}v(x))^2\,
dx + \int_{B_r} g(x)(v^n(x))^{\tau-2} v(x)\, dx \\
\notag \leq  \bigg[\Lambda(a)+C_h r^\e\binom{N}{k}
\bigg(\frac{2}{k-2}\bigg)^{\!\!2}\bigg(1+\binom{N-k}{k}\bigg)+\frac{C(\tau)-\Lambda(a)}2\bigg]
\int_{B_r} |\nabla((v^n(x))^{\frac{\tau}{2}-1}v(x))|^2 \, dx
\\
\notag \quad +\ell_\tau \int_{B_r} (v^n(x))^{\tau-2} (v(x))^2 \,dx
+\int_{B_r} |g(x)|(v^n(x))^{\tau-2} |v(x)|\, dx \ .
 \end{gather}
Let us consider the last term in the right hand side of
\eqref{vn2}. Since $g\in L^2(B_r)$, then by H\"older inequality
\begin{align*} 
\int_{B_r} & |g(x)|(v^n(x))^{\tau-2} |v(x)|\, dx \leq
\|g\|_{L^2(\Omega)} \left(\int_{B_r} (v^n(x))^{2\tau-4}|v(x)|^2 \,
dx\right)^{\!\!\frac{1}2}\\
\notag & =\|g\|_{L^2(\Omega)} \left(\int_{B_r}
(v^n(x))^{\frac{2(\tau-1)(\tau-2)}\tau} (v^n(x))^{\frac{2(\tau-2)}\tau}|v(x)|^2
\, dx\right)^{\!\!\frac{1}2} \\
\notag &\leq \|g\|_{L^2(\Omega)} \left(\int_{B_r}
|(v^n(x))^{\frac{\tau}2-1} v(x)|^{\frac{4(\tau-1)}\tau} \,
dx\right)^{\!\!\frac{1}2}
\end{align*}
and, since $\frac{4(\tau-1)}\tau<2^*$ for any $\tau<\frac{2}{2^*} q_{\rm
lim}$, by H\"older inequality, Sobolev embedding, and Young
inequality, we obtain
\begin{align} \label{est-g}
\int_{B_r} & |g(x)|(v^n(x))^{\tau-2} |v(x)|\, dx\\
\notag & \leq \|g\|_{L^2(\Omega)}
\bigg(\frac{\omega_{N-1}}N\bigg)^{\frac{1}2-\frac{2(\tau-1)}{2^*\tau}}
r^{\frac{N}2-\frac{2N(\tau-1)}{2^*\tau}} \left(\int_{B_r}
|(v^n(x))^{\frac{\tau}2-1} v(x)|^{2^*} \,
dx\right)^{\frac{2(\tau-1)}{2^*\tau}} \\
\notag & \leq \|g\|_{L^2(\Omega)}
\bigg(\frac{\omega_{N-1}}N\bigg)^{\frac{1}2-\frac{2(\tau-1)}{2^*\tau}}
r^{\frac{N}2-\frac{(N-2)(\tau-1)}{\tau}} S^{-\frac{\tau-1}\tau}
\left(\int_{B_r} |\nabla((v^n(x))^{\frac{\tau}2-1} v(x))|^{2} \,
dx\right)^{\frac{\tau-1}{\tau}}  \\
 \notag & \leq\frac{1}{\tau} \|g\|^\tau_{L^2(\Omega)} +\frac{\tau-1}\tau
\bigg(\frac{\omega_{N-1}}N\bigg)^{\frac{\tau}{2(\tau-1)}-\frac{2}{2^*}}
r^{\frac{N\tau}{2(\tau-1)}-N+2}  S^{-1} \int_{B_r}
|\nabla((v^n(x))^{\frac{\tau}2-1} v(x))|^{2} \, dx.
\end{align}
Inserting \eqref{est-g} into \eqref{vn2} we obtain
\begin{align*} 
\bigg[& \frac{C(\tau)-\Lambda(a)}2-C_h r^\e\binom{N}{k}
\bigg(\frac{2}{k-2}\bigg)^{\!\!2}\bigg(1+\binom{N-k}{k}\bigg) \\
\notag & -\frac{\tau-1}\tau
\bigg(\frac{\omega_{N-1}}N\bigg)^{\frac{\tau}{2(\tau-1)}-\frac{2}{2^*}}
r^{\frac{N\tau}{2(\tau-1)}-N+2}  S^{-1} \bigg]
 \int_{B_r}
|\nabla((v^n(x))^{\frac{\tau}2-1} v(x))|^{2} \, dx \\
\notag & \leq \frac{1}{\tau} \|g\|^\tau_{L^2(\Omega)} +\ell_\tau \int_{B_r}
(v^n(x))^{\tau-2} (v(x))^2 \,dx
\end{align*}
and by Sobolev embedding we also have
\begin{align} \label{lim-q}
S\bigg[& \frac{C(\tau)-\Lambda(a)}2-C_h r^\e\binom{N}{k}
\bigg(\frac{2}{k-2}\bigg)^{\!\!2}\bigg(1+\binom{N-k}{k}\bigg) \\
\notag & -\frac{\tau-1}\tau
\bigg(\frac{\omega_{N-1}}N\bigg)^{\frac{\tau}{2(\tau-1)}-\frac{2}{2^*}}
r^{\frac{N\tau}{2(\tau-1)}-N+2}  S^{-1} \bigg]
\bigg( \int_{B_r}
(v^n(x))^{\frac{2^*}2\tau-2^*} |v(x)|^{2^*} \, dx\bigg)^{\!\!\frac2{2^*}} \\
\notag & \leq \frac{1}{\tau} \|g\|^\tau_{L^2(\Omega)} +\ell_\tau \int_{B_r}
(v^n(x))^{\tau-2} (v(x))^2 \,dx \ .
\end{align}
 Since $\tau<\frac{2}{2^*}
q_{\rm lim}$ then $C(\tau)-\Lambda(a)$ is positive and
$\frac{N\tau}{2(\tau-1)}-N+2$ is also positive. Hence we may fix $r$
small enough in such a way that the left hand side of
\eqref{lim-q} becomes positive. Since $v\in L^\tau(B_{r})$, letting
$n\to +\infty$, the right hand side of \eqref{lim-q} remains
bounded and hence by Fatou Lemma we infer that $v\in
L^{\frac{2^*}2 \tau}(B_{r})$. Since $\eta\equiv 1$ in $B_{r/2}$ we
may conclude that $u\in L^{\frac{2^*}2 \tau}(B_{r/2})$. This
completes the proof of the lemma.
\end{pf}

\section{The Almgren type frequency
  function}\label{sec:monot-prop}

Let $u$ be a weak $H^1(\Omega)$-solution to equation (\ref{u}) in
a bounded domain $\Omega\subset \R^N$ containing the origin with
$a$ satisfying \eqref{eq:ateta} and \eqref{eq:lambdamin1}, $h$
satisfying \eqref{H1-2} and $f$ satisfying \eqref{F}.

By \eqref{F} and Sobolev embedding, we infer that the function
$$
V(x):=\left\{
\begin{array}{ll}
\frac{f(x,u(x))}{u(x)}, & \qquad \text{if } u(x)\neq 0, \\
0, & \qquad \text{if } u(x)=0,
\end{array}
\right.
$$
belongs to $L^{N/2}(\Omega)$ and hence we may apply Proposition
\ref{SMETS} to the function $u$. Therefore, throughout this
section, we may fix 
\begin{equation}\label{eq:46}
2^*<q<q_{\rm lim}
\end{equation}
 and $r_q$ as in
Proposition \ref{SMETS} in such a way that $u\in L^q(B_{r_q})$.

 By Remark \ref{rem:l1w11}, the function
\begin{equation} \label{H(r)}
H(r)=\frac{1}{r^{N-1}}\int_{\partial B_r}|u|^2 \, dS
\end{equation}
belongs to $L^1_{{\rm loc}}(0,R)$ for every $R>0$ such that $B_R\subseteq\Omega$.
 It is also easy to verify that
\begin{equation}\label{eq:heq}
  H(r)=\int_{{\mathbb S}^{N-1}}|u(r\theta)|^2\,dS(\theta)
\quad\text{for a.e. }r\in(0,R).
\end{equation}
Further regularity of $H$ is established  in the following lemma.

\begin{Lemma} \label{l:hprime}
  Let $\Omega\subset\R^N$, $N\geq 3$, be a bounded open set such that
  $0\in\Omega$ and let
 $u$ be a weak $H^1(\Omega)$-solution to equation (\ref{u}) in
$\Omega$  with
$a$ satisfying \eqref{eq:ateta} and \eqref{eq:lambdamin1}, $h$
satisfying \eqref{H1-2}, and $f$ satisfying \eqref{F}.
 If $H$ is the function defined in (\ref{H(r)}) and
  $R>0$ is such that $B_R\subseteq\Omega$, then $H\in W^{1,1}_{{\rm
      loc}}(0,R)$ and
\begin{equation}\label{H'}
  H'(r)=\frac{2}{r^{N-1}} \int_{\partial
    B_r} u\,\frac{\partial   u}{\partial \nu}
  \,dS
\end{equation}
in a distributional sense and for a.e. $r\in (0,R)$.
\end{Lemma}
 \begin{pf}
Since $u,\frac{\partial u}{\partial \nu}\in L^2(B_R)$,
by Remark \ref{rem:l1w11}, we
have that
$$
r\mapsto \frac{2}{r^{N-1}} \int_{\partial
    B_r} u\,\frac{\partial   u}{\partial \nu}
  \,dS \in L^1_{{\rm loc}}(0,R).
$$
If $0<s<r<R$, by Fubini's Theorem we obtain
\begin{align*}
  \int_s^r \frac{2}{t^{N-1}} \left(\int_{\partial B_t}
    u\,\frac{\partial u}{\partial \nu}\,dS\right)dt &=\int_s^r \left(
    \int_{{\mathbb S}^{N-1}}2u(t\theta)\frac{\partial u}{\partial
      \nu}(t\theta) \,dS(\theta)\right)dt\\
&= \int_{{\mathbb
      S}^{N-1}}\left( \int_s^r2u(t\theta)\frac{\partial u}{\partial
      \nu}(t\theta) \,dt\right)\,dS(\theta).
\end{align*}
From classical Brezis-Kato \cite{BK} estimates, standard bootstrap,
and elliptic regularity theory, it follows that $u\in C^{1,\tau}_{\rm
  loc}(\Omega\setminus\widetilde\Sigma)$ for some
$\tau\in(0,1)$. Hence, for every $\theta\in{\mathbb
  S}^{N-1}\setminus\Sigma$, and consequently for a.e.
$\theta\in{\mathbb S}^{N-1}$, $\frac{\partial u}{\partial
  \nu}(t\theta)=\frac{d}{dt}u(t\theta)$ for every $t\in(s,r)$.
Therefore, in view of (\ref{eq:heq}), we deduce that
\begin{align*}
  \int_s^r \frac{2}{t^{N-1}} \left(\int_{\partial B_t}
    u\,\frac{\partial u}{\partial \nu}\,dS\right)dt &= \int_{{\mathbb
      S}^{N-1}}\left( \int_s^r
    \frac{d}{dt}|u(t\theta)|^2\,dt\right)\,dS(\theta)\\
  &=\int_{{\mathbb S}^{N-1}}\big(|u(r\theta)|^2-
|u(s\theta)|^2\big)\,dS(\theta)=H(r)-H(s)
\end{align*}
thus proving that $H\in
  W^{1,1}_{{\rm loc}}(0,R)$ and that its weak derivative is given by (\ref{H'}).
\end{pf}

\noindent Now we show that, if $u\not\equiv 0$,
 $H(r)$ does not vanish for every $r\in(0,r_0)$.

 \begin{Lemma} \label{welld}
  Let $\Omega\subset\R^N$, $N\geq 3$, be a bounded open set such that
  $0\in\Omega$,  $a$ satisfy (\ref{eq:ateta}) and (\ref{eq:lambdamin1}),
and  $u\not\equiv 0$ be a weak
$H^1(\Omega)$-solution to (\ref{u}) in $\Omega$, with  $h$
verifying \eqref{H1-2} and $f$ as in  \eqref{F}. Then
$H(r)>0$ for any $r\in (0,r_0)$, where $H=H(r)$ is defined by (\ref{H(r)}) and
$r_0>0$ satisfies (\ref{eq:r_0}).
\end{Lemma}

\begin{pf} Suppose by contradiction that there exists $R\in(0,r_0)$ such that
$H(R)=0$. Then  $u= 0$ a.e. on $\partial B_{R}$ and thus $u\in H^1_0(B_R)$.
  Multiplying both sides of (\ref{u}) by
  $u$ and integrating by parts over $B_{R}$ we obtain
\begin{align*}
\int_{B_{R}} |\nabla u(x)|^2 dx -\int_{B_{R}}
\frac{a(\frac{x}{|x|})}{|x|^2} |u(x)|^2 \, dx = \int_{B_{R}}
h(x)\,|u(x)|^2\, dx+\int_{B_R} f(x,u(x))u(x)\, dx.
\end{align*}
Proceeding as in \eqref{eq:coerc} and using \eqref{F}, H\"older and
Sobolev inequalities, we obtain
\begin{align*} 
0&=\int_{B_{R}}\bigg(|\nabla
  u(x)|^2-\frac{a(\frac{x}{|x|})}{|x|^2}u^2(x)-
  h(x)u^2(x)-f(x,u(x))u(x)\bigg)\,dx\\
\notag & \geq
  \bigg[1-\Lambda(a)-C_h r_0^\e
\binom{N}{k}
 \left(\frac{2}{k-2}\right)^{\!\!2} \left(1+\binom{N-k}{k}\right)\bigg]
\int_{B_{R}}\!\!|\nabla u(x)|^2\,dx \\
& -C_f S^{-1}\left[(\omega_{N-1}/N)^{\frac2N}r_0^{2}+
 \|u\|_{L^{2^*}(B_{r_0})}^{2^*-2}\right]
\int_{B_{R}}\!\!|\nabla u(x)|^2\,dx ,
\end{align*}
which, together with \eqref{eq:r_0}, implies $u\equiv 0$ in
$B_{R}$.  Since $u\equiv 0$ in a neighborhood of the origin, we
may apply, away from the singular set
$\widetilde\Sigma$ (which has zero measure), classical
unique continuation principles for second order elliptic equations
with locally bounded coefficients (see e.g.  \cite{wolff}) to
conclude that $u=0$ a.e. in $\Omega$, a contradiction.
\end{pf}

\noindent
We also consider the function $D:(0,r_0)\to\R$ defined as
\begin{equation} \label{D(r)}
D(r)=\frac{1}{r^{N-2}} \int_{B_r}
  \bigg(|\nabla u(x)|^2-\frac{a\big(\frac{x}{|x|}\big)}{|x|^2}|u(x)|^2
    -h(x)\,|u(x)|^2-f(x,u(x))u(x) \bigg)dx,
\end{equation}
where $r_0$ is defined in \eqref{eq:r_0}. The regularity of the
function $D$ is established in the following lemma.

\begin{Lemma}\label{l:dprime}
  Let $\Omega\subset\R^N$, $N\geq 3$, be a bounded open set such that
  $0\in\Omega$. Let $a$ satisfy (\ref{eq:ateta}) and
  (\ref{eq:lambdamin1}), and $u$ be a
  weak $H^1(\Omega)$-solution to (\ref{u}), with  $h$ satisfying
 \eqref{H1-2} and $f$ satisfying \eqref{F}. Then the
  function $D$ defined in (\ref{D(r)}) belongs to  $W^{1,1}_{{\rm\
      loc}}(0, r_0)$ and
\begin{align}\label{D'F}
  D'(r)=&\frac{2}{r^{N-1}} \bigg[ r\int_{\partial B_r}
  \left|\frac{\partial u}{\partial \nu}\right|^2 dS -\frac{1}2
  \int_{B_r} (\nabla h(x)\cdot x)|u(x)|^2\,dx -\int_{B_r}
  h(x)|u(x)|^2\,dx
  \bigg] \\
  \notag & \qquad +r^{1-N}\int_{B_r}
  \big((N-2)f(x,u(x))u(x)-2NF(x,u(x))-
  2\nabla_xF(x,u(x))\cdot x\big)\, dx \\
  \notag & \qquad +r^{2-N}\int_{\partial B_r}
  \big(2F(x,u(x))-f(x,u(x))u(x)\big)\, dS
\end{align}
in a distributional sense and for a.e. $r\in (0, r_0)$.
\end{Lemma}

\begin{pf}
For any $r\in (0,r_0)$ let
\begin{align} \label{I(r)}
I(r)= \int_{B_r} \bigg(|\nabla u(x)|^2
-\frac{a\big(\frac{x}{|x|}\big)}{|x|^2}|u(x)|^2
-h(x)\,|u(x)|^2-f(x,u(x))u(x) \bigg) \, dx.
\end{align}
From Remark \ref{rem:l1w11},  we deduce that $I\in
W^{1,1}(0,r_0)$ and
\begin{equation} \label{I'(r)}
I'(r) = \int_{\partial B_r} \bigg(|\nabla u(x)|^2
-\frac{a(\frac{x}{|x|})}{|x|^2}|u(x)|^2
    - h(x)|u(x)|^2 -f(x,u(x))u(x)\bigg) \,dS
\end{equation}
 for a.e. $r\in (0,r_0)$ and in the distributional sense.
Therefore $D\in W^{1,1}_{{\rm loc}}(0,r_0)$ and, plugging
\eqref{eq:poho}, (\ref{I(r)}), and (\ref{I'(r)}) into
\begin{align*}
  D'(r)=r^{1-N}[-(N-2)I(r)+rI'(r)],
\end{align*}
we obtain (\ref{D'F}) for a.e. $r\in (0,r_0)$ and in the distributional sense.
\end{pf}

\noindent By virtue of Lemma \ref{welld},
if $u$ is a weak
$H^1(\Omega)$-solution to (\ref{u}), $u\not\equiv 0$,
the \emph{Almgren type frequency
  function}
\begin{equation} \label{N(r)}
{\mathcal N}(r)={\mathcal N}_{u,h,f}(r)=\frac{D(r)}{H(r)}
\end{equation}
is well defined in $(0,r_0)$.
 Collecting Lemmas \ref{l:hprime} and \ref{l:dprime}, we
 compute the derivative of ${\mathcal N}$.

 \begin{Lemma} \label{mono} Let $\Omega\subset\R^N$, $N\geq 3$, be a
   bounded open set such that $0\in\Omega$, $a$ satisfy
   (\ref{eq:ateta}) and (\ref{eq:lambdamin1}), and $u\not\equiv 0$ be
   a weak $H^1(\Omega)$-solution to (\ref{u}), with $h$ satisfying
   \eqref{H1-2} and $f$ satisfying \eqref{F}. Then the function
   ${\mathcal N}$ defined in (\ref{N(r)}) belongs to $W^{1,1}_{{\rm
       loc}}(0, r_0)$ and
\begin{align}\label{formulona}
{\mathcal N}'(r)=\nu_1(r)+\nu_2(r)
\end{align}
in a distributional sense and for a.e. $r\in (0,r_0)$,
where
\begin{align}\label{eq:nu1}
\nu_1(r)=&\frac{2r\Big[
    \left(\int_{\partial B_r} \left|\frac{\partial
          u}{\partial\nu}\right|^2 dS\right) \cdot
    \left(\int_{\partial B_r} |u|^2 dS\right)-\left(\int_{\partial
        B_r} u\frac{\partial  u}{\partial \nu}
      dS\right)^{\!2} \Big]} {\left(\int_{\partial B_r} |u|^2
    dS\right)^2}
\end{align}
and
\begin{align}\label{eq:nu2}
\nu_2(r)=
&
-\frac{\int_{B_r}
(2h(x)+\nabla h(x)\cdot x)|u(x)|^2\,dx}{\int_{\partial B_r}|u|^2 \, dS}
+\frac{r\int_{\partial B_r}\big(2F(x,u(x))-f(x,u(x))u(x)\big)\, dS}
{\int_{\partial B_r}|u|^2 \, dS} \\[10pt]
&\notag
+\frac{\int_{B_r} \big((N-2)f(x,u(x))u(x)-2NF(x,u(x))
-2\nabla_xF(x,u(x))\cdot x\big)\, dx}
{\int_{\partial B_r}|u|^2 \, dS} .
\end{align}
\end{Lemma}

\begin{pf} From Lemmas \ref{l:hprime}, \ref{welld}, and
  \ref{l:dprime}, it follows that ${\mathcal N}\in W^{1,1}_{{\rm
      loc}}(0,r_0)$. From (\ref{eq:poho2}), (\ref{D(r)}), and
  (\ref{H'}) we infer
\begin{equation} \label{ul}
D(r)=\frac{1}{2}r H'(r)
\end{equation}
for a.e. $r\in (0,r_0)$.
From (\ref{ul}) we have that
$$
{\mathcal N}'(r)=\frac{D'(r)H(r)-D(r)H'(r)}{(H(r))^2}
=\frac{D'(r)H(r)-\frac{1}{2} r (H'(r))^2}{(H(r))^2}
$$
and the proof of the lemma easily follows from (\ref{H'}) and
(\ref{D'F}).
\end{pf}

\noindent We now prove that ${\mathcal N}(r)$ admits a finite limit as
$r\to 0^+$. To this aim, the following estimate plays a crucial role.

\begin{Lemma} \label{l:stimasotto}
Under the same assumptions as in Lemma \ref{mono},
there
exist $\tilde r\in (0,\min\{r_0,r_q\})$ and a positive constant
$\overline{C}=\overline{C}(N,k,a,h,f,u)>0$
depending on $N$, $k$, $a$, $h$, $f$,  $u$
but independent of $r$ such that
\begin{align}\label{eq:47}
\int_{B_r}
  \bigg(&|\nabla u(x)|^2-\frac{a\big(\frac{x}{|x|}\big)}{|x|^2}|u(x)|^2
    -h(x)\,|u(x)|^2-f(x,u(x))u(x) \bigg)dx\\[10pt]
   \notag \geq &-\frac{N-2}{2r}
\int_{\partial B_r} |u(x)|^2 dS \\[10pt]
&\notag+\overline{C}
\bigg(\sum_{J\in{\mathcal A}_k}\int_{B_r} \frac{|u(x)|^2}{|x_{J}|^{2}}\,dx
+\sum_{(J_1,J_2)\in \mathcal B_k}\int_{B_r} \frac{|u(x)|^2}{
|x_{J_1}-x_{J_2}|^{2}}+
 \bigg(\int_{B_r}|u(x)|^{2^*}dx\bigg)^{\!\!\frac2{2^*}}\bigg)
\end{align}
and
\begin{equation}\label{Nbelow}
   {\mathcal N}(r)>-\frac{N-2}{2}
 \end{equation}
for every $r\in(0,\tilde r)$.
\end{Lemma}
\begin{pf}
 By (\ref{eq:11}), (\ref{eq:11-MB}), and
  (\ref{coercivitySOB}), we have that
\begin{align*}
\int_{B_r}
  \bigg(&|\nabla u(x)|^2-\frac{a\big(\frac{x}{|x|}\big)}{|x|^2}|u(x)|^2
    -h(x)\,|u(x)|^2-f(x,u(x))u(x) \bigg)dx\\[10pt]
   \notag \geq &-\frac{N-2}{2r}
\frac{\binom{N}{k}\big(1+\binom{N-k}{k}\big)+\frac{1+\Lambda(a)}{2}}
{\binom{N}{k}\big(1+\binom{N-k}{k}\big)+1}
\int_{\partial B_r} |u(x)|^2 dS \\[10pt]
&\notag+
\frac{\left(\frac{k-2}2\right)^2}{\binom{N}{k}\big(1+\binom{N-k}{k}\big)+1}
  \bigg[1-\Lambda(a)
-C_hr^\e
\frac{\binom{N}{k}\big(1+\binom{N-k}{k}\big)+1}
{(\frac{k-2}{2})^{2}}\bigg]\times\\
&\notag\hskip4cm\times
\bigg(\sum_{J\in{\mathcal A}_k}\int_{B_r} \frac{|u(x)|^2}{|x_{J}|^{2}}\,dx
+\sum_{(J_1,J_2)\in \mathcal B_k}\int_{B_r} \frac{|u(x)|^2}{
|x_{J_1}-x_{J_2}|^{2}}\bigg)
\\[10pt]
&\notag+\bigg[
\frac12\widetilde S_N\frac{\min\big\{
1-\Lambda(a),\mu_1(a)+\big(\frac{N-2}{2}\big)^{2}\big\}}
{\binom{N}{k}\big(1+\binom{N-k}{k}\big)+1}\\
&\notag\hskip2cm
-C_f
\left(\Big(\frac{\omega_{N-1}}N\Big)^{\!\frac2N}r^{2}+
 \|u\|_{L^{2^*}(B_r)}^{2^*-2}\right)\bigg]
 \bigg(\int_{B_r}|u(x)|^{2^*}dx\bigg)^{\!\!2/2^*}
\end{align*}
for every $r\in(0,r_0)$. Since $\Lambda(a)<1$, from the above estimate
it follows that we can choose $\tilde r\in (0,r_0)$ sufficiently small
such that estimate (\ref{eq:47}) holds for $r\in (0,\tilde r)$ for
some positive constant
$\overline{C}=\overline{C}(N,k,a,h,f,u)>0$.
   Estimate (\ref{eq:47}), together with
(\ref{H(r)}) and (\ref{D(r)}), yields
(\ref{Nbelow}).
\end{pf}

\begin{Lemma} \label{l:stima_nu2} Under the same assumptions as in
  Lemma \ref{mono}, let $\tilde r $ be as in Lemma \ref{l:stimasotto} and
  $\nu_2$ as in (\ref{eq:nu2}). Then there exist a positive constant
  $C_1>0$ depending on $N,q,C_f,C_h,\overline{C},\tilde
  r,\|u\|_{L^q(B_{\tilde r})}$ and a function $g\in L^1(0,\tilde r)$,
  $g\geq 0$ a.e. in $(0,\tilde r)$, such that
$$
|\nu_2(r)|\leq C_1\left[{\mathcal
      N}(r)+\frac{N}{2}\right]\Big(r^{-1+\e}+r^{-1+
\frac{2(q-2^*)}{q}}+g(r)\Big)
$$
for a.e. $r\in (0,\tilde r)$ and
$$
\int_0^r g(s)\,ds\leq \frac {\|u\|_{L^{2^*}(\Omega)}^{2^*(1-\alpha)}}{1-\alpha}
\, r^{\frac{N(q-2^*)}q(\alpha-\frac{2}{2^*})}
$$
for all $r\in (0,\tilde r)$ and for some $\alpha$ satisfying 
$\frac{2}{2^*}<\alpha<1$.
\end{Lemma}

\begin{pf}
From \eqref{H1-2} and  (\ref{eq:47}) we deduce that
\begin{align*}
\left|\int_{B_r} (2h(x)+\nabla h(x)\cdot x)|u(x)|^2\,dx\right|&\leq
2C_h r^\e
\bigg(\sum_{J\in{\mathcal A}_k}\int_{B_r} \frac{|u(x)|^2}{|x_{J}|^{2}}\,dx
+\sum_{(J_1,J_2)\in \mathcal B_k}\int_{B_r} \frac{|u(x)|^2}{
|x_{J_1}-x_{J_2}|^{2}}\bigg)\\
&\leq  2C_h\overline{C}^{-1}\,
r^{\e+N-2}\left[D(r)+{\textstyle{\frac{N-2}{2}}}H(r)\right],
\end{align*}
and, therefore, for any $r\in (0,\tilde r)$, we have that
\begin{align}\label{B00}
  \left|\frac{\int_{B_r} (2h(x)+\nabla h(x)\cdot x)|u(x)|^2\,dx}
    {\int_{\partial B_r}|u|^2 \, dS}\right|&\leq 2C_h\overline{C}^{-1}
  \,r^{-1+\e} \frac{D(r)+\frac{N-2}{2}H(r)}{H(r)} \\
&\notag=
  2C_h\overline{C}^{-1}\, r^{-1+\e}\left[{\mathcal N}(r)+
    \frac{N-2}{2}\right].
\end{align}
By \eqref{F}, H\"older's inequality, and (\ref{eq:47}), for some
constant ${\rm const\,}={\rm const\,}(N,C_f)>0$ depending on
$N,C_f$, and for all $r\in(0,\tilde r)$, there holds
\begin{align*}
  \bigg|&\int_{B_r} \big((N-2)f(x,u(x))u(x)-2NF(x,u(x))
  -2\nabla_xF(x,u(x))\cdot x\big)\, dx\bigg|
  \\
  &\leq {\rm const\,}\int_{B_r}\big(|u(x)|^2+|u(x)|^{2^*}\big)\,dx
  \\
  &\leq {\rm const\,}
  \left(\Big(\frac{\omega_{N-1}}N\Big)^{\!\frac2N}r^{2}+
     \|u\|_{L^{2^*}(B_r)}^{2^*-2}\right)
  \bigg(\int_{B_r}|u(x)|^{2^*}dx\bigg)^{\!\!2/2^*}
  \\
  &
\leq {\rm const\,}
 \left(\Big(\frac{\omega_{N-1}}N\Big)^{\!\frac2N}r^{2}+
    \Big(\frac{\omega_{N-1}}N\Big)^{\!\frac{2(q-2^*)}{Nq}}
    r^{\frac{2(q-2^*)}{q}} \|u\|_{L^{q}(B_{\tilde r})}^{2^*-2}\right)
\\
& \ \ \ \times \overline C^{-1}
r^{N-2}\left[D(r)+{\textstyle{\frac{N-2}{2}}}H(r)\right]
\end{align*}
and hence
\begin{multline}\label{eq:20}
  \left| \frac{\int_{B_r} \big((N-2)f(x,u(x))u(x)-2NF(x,u(x))
      -2\nabla_xF(x,u(x))\cdot x\big)\, dx}
    {\int_{\partial B_r}|u|^2 \, dS} \right|
    \\
\leq {\rm const\,}\overline C^{-1}
 \left(\Big(\frac{\omega_{N-1}}N\Big)^{\!\frac2N}
 \tilde r^{\frac{22^*}{q}}
 +\Big(\frac{\omega_{N-1}}N\Big)^{\frac{2(q-2^*)}{Nq}}
 \|u\|_{L^{q}(B_{\tilde r})}^{2^*-2}\right)
\\
\ \ \ \times r^{-1+\frac{2(q-2^*)}{q}}
\left[\mathcal N(r)+{\textstyle{\frac{N-2}{2}}}\right] .
\end{multline}
Let us fix $\frac{2}{2^*}<\alpha<1$. Then, by H\"older's inequality and
(\ref{eq:47}),
\begin{align}\label{eq:21}
  \bigg(\int_{B_r}&|u(x)|^{2^*}\,dx\bigg)^{\!\!\alpha} = 
\bigg(\int_{B_r}|u(x)|^{2^*}\,dx\bigg)^{\!\!\alpha-\frac{2}{2^*}}
  \bigg(\int_{B_r}|u(x)|^{2^*}\,dx\bigg)^{\!\!\frac2{2^*}}\\
  &\notag\leq
\Big(\frac{\omega_{N-1}}N\Big)^{\frac{q-2^*}q(\alpha-\frac{2}{2^*})} 
r^{\frac{N(q-2^*)}q(\alpha-\frac{2}{2^*})}
\|u\|_{L^q(B_{\tilde r})}^{2^*(\alpha-\frac{2}{2^*})} 
  \overline{C}^{-1}
  r^{N-2}
  \left[D(r)+{{\frac{N-2}{2}}}H(r)\right]\\
  &\notag=
  \overline{C}^{-1}\Big(\frac{\omega_{N-1}}N\Big)^{\!\!\frac{\beta}N}
  r^{-1+\beta}
\|u\|_{L^q(B_{\tilde r})}^{2^*(\alpha-\frac{2}{2^*})} 
\left[{\mathcal
      N}(r)+\frac{N-2}{2}\right]\bigg(\int_{\partial B_r}|u|^2 \, dS\bigg)
\end{align}
for all $r\in(0,\tilde r)$, where
$\beta=\frac{N(q-2^*)}q(\alpha-\frac{2}{2^*})>0$. From
\eqref{F}, (\ref{eq:21}), and (\ref{Nbelow}), there exists some ${\rm
  const\,}={\rm const\,}(N,q,C_f)>0$ depending on $N,q,C_f$ such
that, for all $r\in(0,\tilde r)$,
\begin{multline}\label{eq:27}
  \left|\frac{r\int_{\partial B_r}\big(2F(x,u(x))-f(x,u(x))u(x)\big)\,
      dx} {\int_{\partial B_r}|u|^2 \, dS} \right| \leq {\rm const\,}
  r\left(1+\frac{\int_{\partial B_r}|u|^{2^*} \, dS}{\int_{\partial
        B_r}|u|^2 \, dS}\right)\\
\leq {\rm const\,} r \left[{\mathcal N}(r)+\frac{N}{2}\right] 
  +\frac{\rm
    const\,}{\overline{C}
}\Big(\frac{\omega_{N-1}}N\Big)^{\!\!\frac{\beta}N}
 \|u\|_{L^q(B_{\tilde
        r})}^{{2^*}(\alpha-\frac{2}{2^*})} \left[{\mathcal
        N}(r)+\frac{N-2}{2}\right]\frac{r^{\beta}  \int_{\partial B_r}|u|^{2^*} \,
    dS}{\Big(\int_{B_r}|u(x)|^{2^*} \, dx\Big)^{\alpha}}.
\end{multline}
By a direct calculation, we have that
\begin{multline}\label{eq:22}
  \frac{r^{\beta}\int_{\partial B_r}|u|^{2^*} \!
    dS}{\Big(\int_{B_r}|u(x)|^{2^*} \! dx\Big)^{\!\alpha}}
    =\frac{1}{1-\alpha}\bigg[\frac{d}{dr}\bigg(r^{\beta}
  \bigg(\int_{B_r}|u(x)|^{2^*} \! dx\bigg)^{\!\!1-\alpha}\bigg) -
  \beta\,  r^{-1+\beta}
  \bigg(\int_{B_r}|u(x)|^{2^*} \! dx\bigg)^{\!\!1-\alpha}\bigg]
\end{multline}
in the distributional sense and for a.e. $r\in (0,\tilde r)$.
Since
$$
\lim_{r\to0^+}r^{\beta}
    \bigg(\int_{B_r}|u(x)|^{2^*} \! dx\bigg)^{\!\!1-\alpha}=0
$$
we deduce that the function
$$
r\mapsto
\frac{d}{dr}\bigg(r^{\beta}
    \bigg(\int_{B_r}|u(x)|^{2^*} \! dx\bigg)^{\!\!1-\alpha}\bigg)
$$
is integrable over $(0,\tilde r)$. Being
$$
r^{-1+\beta} \bigg(\int_{B_r}|u(x)|^{2^*} \!
dx\bigg)^{\!\!1-\alpha} =o(r^{-1+\beta})
$$
as $r\to 0^+$, we have that also the function
$$
r\mapsto
r^{-1+\beta} \bigg(\int_{B_r}|u(x)|^{2^*} \!
dx\bigg)^{\!\!1-\alpha}
$$
is integrable over $(0,\tilde r)$. Therefore, by (\ref{eq:22}), we deduce that
\begin{equation}\label{eq:28}
g(r):= \frac{r^{\beta}\int_{\partial B_r}|u|^{2^*} \!
  dS}{\Big(\int_{B_r}|u(x)|^{2^*} \! dx\Big)^{\alpha}}\in L^1(0,\tilde r)
\end{equation}
and
\begin{equation}\label{eq:29}
0\leq \int_0^r g(s)\,ds\leq \frac{\|u\|_{L^{2^*}(\Omega)}^{2^*(1-\alpha)}}{1-\alpha}
\, r^{\beta}
\end{equation}
for all $r\in (0,\tilde r)$. Collecting (\ref{B00}), (\ref{eq:20}),
(\ref{eq:27}), (\ref{eq:28}), and (\ref{eq:29}), we obtain the stated
estimate.
\end{pf}

\begin{Lemma} \label{l:stima_N_sopra} Under the same assumptions as in
  Lemma \ref{mono}, let $\tilde r $ be as in Lemma \ref{l:stimasotto}
  and ${\mathcal N}$ as in (\ref{N(r)}). Then there exist a positive
  constant $C_2>0$ depending on $N$,  $q$, $C_f$, $C_h$,
  $\overline{C}$, $\tilde r$, $\|u\|_{L^{q}(B_{\tilde r})}$, ${\mathcal
    N}(\tilde r)$, $\e$ such that
\begin{equation} \label{Nabove}
{\mathcal N}(r)\leq C_2
\end{equation}
for all $r\in (0,\tilde r)$.
\end{Lemma}
\begin{pf}
By Lemma \ref{mono}, Schwarz's inequality, and Lemma \ref{l:stima_nu2}, we
 obtain
\begin{equation}\label{eq:40}
\bigg({\mathcal N}+\frac N2\bigg)'(r)\geq
\nu_2(r)\geq
 -C_1\left[{\mathcal
      N}(r)+\frac{N}{2}\right]\Big(r^{-1+\e}
      +r^{-1+\frac{2(q-2^*)}{q}}+g(r)\Big)
\end{equation}
for a.e. $r\in (0,\tilde r)$.
After integration over $(r,\tilde r)$ it follows that
\begin{equation*}
{\mathcal N}(r)\leq -\frac N2+
\left({\mathcal
      N}(\tilde r)+\frac{N}{2}\right)
\exp\left(
C_1\bigg(
\frac{\tilde r^\e}{\e}+ \frac{q}{2(q-2^*)}\tilde
r^{\frac{2(q-2^*)}{q}}
+\int_0^{\tilde r} g(s)\,ds\bigg)\right)
\end{equation*}
for any $r\in (0,\tilde r)$, thus proving estimate (\ref{Nabove}).
\end{pf}

\begin{Lemma} \label{gamma}
Under the same assumptions as in Lemma \ref{mono}, the limit
$$
\gamma:=\lim_{r\rightarrow 0^+} {\mathcal N}(r)
$$
exists and is finite.
\end{Lemma}
\begin{pf}
  By Lemmas \ref{l:stima_nu2} and \ref{l:stima_N_sopra}, the function
  $\nu_2$ defined in (\ref{eq:nu2}) belongs to $L^1(0,\tilde r)$.
  Hence, by Lemma \ref{mono} and Schwarz's inequality, ${\mathcal N}'$
  is the sum of a nonnegative function and of a $L^1$-function on
  $(0,\tilde r)$.  Therefore
$$
{\mathcal N}(r)={\mathcal N}(\tilde r)-\int_r^{\tilde r} {\mathcal N}'(s)\, ds
$$
admits a limit as $r\rightarrow 0^+$ which is necessarily finite in view of
(\ref{Nbelow}) and (\ref{Nabove}).
\end{pf}

\noindent A first consequence of the above analysis on the Almgren's
frequency function is the following estimate of $H(r)$.
\begin{Lemma} \label{l:uppb} Under the same assumptions as in Lemma
  \ref{mono}, let $\gamma:=\lim_{r\rightarrow 0^+} {\mathcal N}(r)$ be as
  in Lemma \ref{gamma} and $\tilde r $ as in Lemma \ref{l:stimasotto}.
  Then there exists a constant $K_1>0$ such that
\begin{equation} \label{1stest}
H(r)\leq K_1 r^{2\gamma}  \quad \text{for all } r\in (0,\tilde r).
\end{equation}
On the other hand for any $\sigma>0$ there exists a constant
$K_2(\sigma)>0$ depending on $\sigma$ such that
\begin{equation} \label{2ndest} H(r)\geq K_2(\sigma)\,
  r^{2\gamma+\sigma} \quad \text{for all } r\in (0,\tilde r).
\end{equation}
\end{Lemma}

\begin{pf}
  By Lemma \ref{gamma}, ${\mathcal N}'\in L^1(0,\tilde r )$ and, by
  Lemma \ref{l:stima_N_sopra}, ${\mathcal N}$ is bounded, then from
  (\ref{eq:40}) and (\ref{eq:29}) it follows that
  \begin{equation} \label{qsopra}
{\mathcal N}(r)-\gamma=\int_0^r
    {\mathcal N}'(s) \, ds\geq -C_3 r^\delta
\end{equation}
for some constant $C_3>0$ and all $r\in(0, \tilde r)$, where
\begin{equation}\label{eq:45}
\delta=\min\bigg\{\e, 
\frac{N(q-2^*)}{q}\bigg(\alpha-\frac2{2^*}\bigg),
\frac{2(q-2^*)}q
\bigg\}
\end{equation}
with $\alpha$ as in Lemma \ref{l:stima_nu2}.

Therefore by (\ref{ul})
and (\ref{qsopra}) we deduce that, for $r\in(0,\tilde r)$,
$$
\frac{H'(r)}{H(r)}=\frac{2\,{\mathcal N}(r)}{r}\geq
\frac{2\gamma}{r}-2C_3 r^{-1+\delta},
$$
which, after integration over the interval $(r, \tilde r)$, yields
(\ref{1stest}).

Let us prove (\ref{2ndest}). Since $\gamma=\lim_{r\rightarrow 0^+}
{\mathcal N}(r)$, for any $\sigma>0$ there exists $r_\sigma>0$ such
that ${\mathcal N}(r)<\gamma+\sigma/2$ for any $r\in (0,r_\sigma)$ and
hence
$$
\frac{H'(r)}{H(r)}=\frac{2\,{\mathcal N}(r)}{r}<\frac{2\gamma+\sigma}{r}
\quad \text{for all } r\in (0,r_\sigma).
$$
Integrating over the interval $(r,r_\sigma)$ and by continuity of $H$
outside $0$, we obtain (\ref{2ndest}) for some constant $K_2(\sigma)$
depending on $\sigma$.
\end{pf}

\section{The blow-up argument}\label{sec:blow-up-argument}

Throughout this section we let $u$ be a weak $H^1(\Omega)$-solution to
equation (\ref{u}) in a bounded domain $\Omega\subset \R^N$ containing
the origin with $a$ satisfying \eqref{eq:ateta},
\eqref{eq:lambdamin1}, $h$ satisfying \eqref{H1-2}, and $f$ satisfying
\eqref{F}. Let $H$ and $D$ be the functions defined in (\ref{H(r)})
and (\ref{D(r)}) and $\tilde r$ be as in Lemma \ref{l:stimasotto}.

\begin{Lemma} \label{bound}
For $ \lambda\in(0,\tilde r)$, let
\begin{equation}\label{eq:defwlambda}
w^\lambda(x)=\frac{u(\lambda x)}{\sqrt{H(\lambda)}}.
\end{equation}
Then there exists $\bar r\in(0,\tilde r)$ depending on $N$, $k$,
$a$, $h$, $f$, $\e$, and $\|u\|_{L^{q}(B_{\tilde r})}$ such
that the set $\{w^\lambda\}_{\lambda\in (0,\bar r)}$ is bounded in
$H^1(B_1)$.
\end{Lemma}
\begin{pf}
From (\ref{eq:heq}) it follows that $\int_{\partial B_1}|w^{\lambda}|^2dS=1$.
 Moreover, by
scaling and (\ref{Nabove}),
\begin{multline}\label{eq:8}
  \int_{B_{1}} |\nabla w^\lambda(x)|^2 dx -\int_{B_{1}}
  \frac{a(\frac{x}{|x|})}{|x|^{2}}|w^\lambda(x)|^2\, dx
  -\lambda^2\int_{B_1} h(\lambda x)|w^\lambda(x)|^2\, dx\\
  -\frac{\lambda^2}{\sqrt{H(\lambda)}}\int_{B_1} f(\lambda
  x,\sqrt{H(\lambda)}w^\lambda(x)) w^\lambda(x)\, dx
  ={\mathcal N}(\lambda)\leq C_2
\end{multline}
for every $\lambda\in(0,\tilde r)$.
By \eqref{coercivity} applied to $w^\lambda$ we have that
\begin{multline}\label{eq:12}
  \int_{B_1}  |\nabla w^\lambda (x)|^2 \, dx-\int_{B_1}
  \frac{a(\frac{x}{|x|})}{|x|^{2}}|w^\lambda(x)|^2\, dx+\Lambda(a)\frac{N-2}2
  \int_{\partial B_1} |w^\lambda |^2 \, dS\\
\geq (1-\Lambda(a)) \int_{B_1} |\nabla w^\lambda (x)|^2 \, dx.
\end{multline}
Moreover by Corollary \ref{c:hardyboundary-0} we have
\begin{align}\label{eq:7}
\bigg| & \lambda^2\int_{B_1} h(\lambda x)|w^\lambda(x)|^2\, dx \bigg |
\leq C_h \lambda^\e
\bigg(\sum_{J\in{\mathcal
A}_k}
\int_{B_1} \frac{|w^\lambda(x)|^2}{|x_{J}|^2}\, dx
+\sum_{(J_1,J_2)\in \mathcal B_k}
\int_{B_1} \frac{|w^\lambda(x)|^2}
{|x_{J_1}-x_{J_2}|^{2}}\, dx\bigg) \\
&
\notag\leq
C_h\lambda^\e\binom{N}{k}
\bigg(1+\binom{N-k}{k}\bigg)
 \left[
\left(\frac{2}{k-2}\right)^2  \int_{B_1} |\nabla w^\lambda (x)|^2 \, dx
+\frac{2(N-2)}{(k-2)^2} \int_{\partial B_1} |w^\lambda|^2 \, dS \right].
\end{align}
From \eqref{F}, H\"older's inequality, Lemma \ref{l:se}, and Lemma
\ref{l:hardyboundary},
\begin{gather}\label{eq:26}
  \frac{\lambda^2}{\sqrt{H(\lambda)}}\bigg|\int_{B_1} f(\lambda
  x,\sqrt{H(\lambda)}w^\lambda(x)) w^\lambda(x)\, dx\bigg|
\\
  \notag\leq C_f\lambda^2\int_{B_1}|w^\lambda(x)|^2\,dx
  +C_f\lambda^2(H(\lambda))^{\frac {2^*}2-1}\int_{B_1}|w^\lambda(x)|^{2^*}\,dx
\\
  \notag\leq C_f
\bigg(
\big({\textstyle{\frac{\omega_{N-1}}N}}
  \big)^{\frac 2N}
\lambda^2+\lambda^2(H(\lambda))^{\frac
    {2^*}2-1}\bigg(\int_{B_1}|w^\lambda(x)|^{2^*}\,dx\bigg)^{\!\!2/N}\bigg)
\bigg(\int_{B_1}|w^\lambda(x)|^{2^*}\,dx\bigg)^{\!\!2/2^*}
  \\
  \notag\leq C_f
\widetilde S_N^{-1}
\bigg(\big({\textstyle{\frac{\omega_{N-1}}N}}
  \big)^{\frac 2N}
\lambda^2+\bigg(\int_{B_\lambda}|u(x)|^{2^*}\,dx\bigg)^{\!\!2/N}\bigg)
\bigg(\int_{B_1}\bigg( |\nabla
  w^\lambda(x)|^2+\frac{|w^\lambda(x)|^2}{|x|^2}\bigg)\,dx\bigg)
\\
\notag\leq
\frac{C_f}{\widetilde S_N}
\bigg(\big({\textstyle{\frac{\omega_{N-1}}N}}
  \big)^{\frac 2N}
\lambda^2+ \|u\|_{L^{q}(B_{\tilde r})}^{2^*-2}\big({\textstyle{\frac{\omega_{N-1}}N}}
  \big)^{\frac{2(q-2^*)}{qN}}\lambda^{\frac{2(q-2^*)}{q}}
\bigg)
\bigg(\int_{B_1}\bigg( |\nabla
  w^\lambda(x)|^2+\frac{|w^\lambda(x)|^2}{|x|^2}\bigg)\,dx\bigg)
 \\
  \notag\leq
\frac{C_f((N-2)^2+4)}{\widetilde S_N(N-2)^2}
\bigg(
\big({\textstyle{\frac{\omega_{N-1}}N}}
  \big)^{\frac 2N}
\lambda^2+
\|u\|_{L^{q}(B_{\tilde r})}^{2^*-2}\big({\textstyle{\frac{\omega_{N-1}}N}}
  \big)^{\frac{2(q-2^*)}{qN}}\lambda^{\frac{2(q-2^*)}{q}}
\bigg)
  \bigg(\int_{B_1}|\nabla w^\lambda(x)|^2\,dx\bigg) \\
  \notag+
\frac{2C_f}{\widetilde S_N(N-2)}
\bigg(\big({\textstyle{\frac{\omega_{N-1}}N}}
  \big)^{\frac 2N}
\lambda^2+ 
\|u\|_{L^{q}(B_{\tilde r})}^{2^*-2}\big({\textstyle{\frac{\omega_{N-1}}N}}
  \big)^{\frac{2(q-2^*)}{qN}}\lambda^{\frac{2(q-2^*)}{q}}
\bigg).
\end{gather}
From (\ref{eq:8}--\ref{eq:26}), we
deduce that
\begin{multline*}
  \bigg[1-\Lambda(a)-C_h\lambda^\e\binom{N}{k}
  \bigg(1+\binom{N-k}{k}\bigg)
  \left(\frac{2}{k-2}\right)^2\\
  \qquad- \frac{C_f((N-2)^2+4)}{\widetilde S_N(N-2)^2}
  \bigg(
\big({\textstyle{\frac{\omega_{N-1}}N}} \big)^{\frac 2N}
  \lambda^2+
\|u\|_{L^{q}(B_{\tilde r})}^{2^*-2}\big({\textstyle{\frac{\omega_{N-1}}N}}
  \big)^{\frac{2(q-2^*)}{qN}}\lambda^{\frac{2(q-2^*)}{q}}
\bigg)
  \bigg]\int_{B_1} |\nabla w^\lambda (x)|^2 \, dx\\
  \leq C_2+\Lambda(a)\frac{N-2}2+C_h
  \lambda^\e\binom{N}{k}\bigg(1+\binom{N-k}{k}\bigg)
  \,\frac{2(N-2)}{(k-2)^2}\\ + \frac{2C_f}{\widetilde S_N(N-2)}
  \bigg(
\big({\textstyle{\frac{\omega_{N-1}}N}} \big)^{\frac 2N}
  \lambda^2+
\|u\|_{L^{q}(B_{\tilde r})}^{2^*-2}\big({\textstyle{\frac{\omega_{N-1}}N}}
  \big)^{\frac{2(q-2^*)}{qN}}\lambda^{\frac{2(q-2^*)}{q}}
\bigg)
\end{multline*}
for every $\lambda\in(0,\tilde r)$, which
implies that $\{w^\lambda\}_{\lambda\in(0,\bar r)}$ is
bounded in $H^1(B_1)$ if $\bar r$ is chosen sufficiently small.
\end{pf}

In the next lemma we prove a {\em{doubling}} type result.

\begin{Lemma} \label{doubling}
There exists $C_4>0$ such that
\begin{equation} \label{dou1}
\frac{1}{C_4} H(\lambda)\leq H(R\lambda)\leq C_4 H(\lambda) \quad
\text{for any } \lambda\in (0,\tilde r/2)
\text{ and }
R\in \left[1,2\right],
\end{equation}
\begin{equation} \label{dou2}
\int_{B_R} |\nabla w^{\lambda}(x)|^2 dx\leq 2^{N-2} C_4 \int_{B_1}
|\nabla w^{R\lambda}(x)|^2 dx \quad
\text{for any } \lambda\in (0,\tilde r/2)
\text{ and }
R\in \left[1,2\right],
\end{equation}
and
\begin{equation} \label{dou3}
\int_{B_R} |w^{\lambda}(x)|^2 dx\leq 2^{N} C_4 \int_{B_1}
|w^{R\lambda}(x)|^2 dx \quad
\text{for any } \lambda\in (0,\tilde r/2)
\text{ and }
R\in \left[1,2\right],
\end{equation}
where $w^\lambda$ is defined in (\ref{eq:defwlambda}).
\end{Lemma}

\begin{pf}
By (\ref{Nbelow}), (\ref{Nabove}), and (\ref{ul}), it follows that
$$
-\frac{N-2}{r}\leq \frac{H'(r)}{H(r)}=\frac{2\,\mathcal
N(r)}{r}\leq \frac{2\,C_2}{r} \qquad {\rm for \ any } \
r\in (0,\tilde r).
$$
Let $R\in(1,2]$. For any $\lambda<\tilde r/R$, integration over
$(\lambda,R\lambda)$ and the fact  that $R\leq 2$ yield
$$
2^{2-N} H(\lambda)\leq H(R\lambda)\leq 4^{C_2}
H(\lambda) \qquad {\rm for \ any} \ \lambda\in (0,\tilde r/R).
$$
Since the above chain of inequalities trivially holds also for $R=1$,
the proof of (\ref{dou1}) is complete with $C_4=\max\{4^{C_2},2^{N-2}\}$.
 By scaling and (\ref{dou1}),
 we obtain that, for any  $\lambda\in (0,\tilde r/2)$ and $R\in [1,2]$,
\begin{align*}
\int_{B_R} |\nabla w^\lambda (x)|^2
dx&=\frac{\lambda^{2-N}}{H(\lambda)} \int_{B_{R\lambda}} |\nabla
u(x)|^2 dx \\
&= R^{N-2} \frac{H(R\lambda)}{H(\lambda)} \int_{B_1} |\nabla
w^{R\lambda} (x)|^2 dx \leq  R^{N-2} C_4 \int_{B_1} |\nabla
w^{R\lambda} (x)|^2 dx,
\end{align*}
thus providing (\ref{dou2}). In a similar way,
(\ref{dou3}) follows from (\ref{dou1}) by   scaling.
\end{pf}

\begin{Lemma} \label{Rlambda}
For every $\lambda\in(0,\tilde r)$, let $w^\lambda$ as in (\ref{eq:defwlambda}).
Then there exist $M>0$ and $\lambda_0>0$ such that for any
$\lambda\in (0,\lambda_0)$ there exists $R_\lambda \in
\left[1,2\right]$ such that
\begin{equation*}
\int_{\partial B_{R_\lambda}} |\nabla w^\lambda|^2 dS\leq M
\int_{B_{R_\lambda}} |\nabla w^\lambda(x)|^2 dx.
\end{equation*}
\end{Lemma}

\begin{pf}
We recall that, by Lemma \ref{bound}, the set
$\{w^\lambda\}_{\lambda\in (0,\bar r)}$ is bounded in
$H^1(B_1)$. Moreover by Lemma
\ref{doubling}, we have that the set
$\{w^\lambda\}_{\lambda\in (0,\bar r/2)}$ is
bounded in $H^1(B_2)$ and hence
\begin{equation} \label{limitata}
\limsup_{\lambda\rightarrow 0^+} \int_{B_2} |\nabla
w^\lambda(x)|^2 dx<+\infty.
\end{equation}
Let us denote, for every $\lambda\in (0,\bar r/2)$,
$$
f_{\lambda}(r)=\int_{B_r}|\nabla w^{\lambda}(x)|^2 \, dx.
$$
The function $f_{\lambda}$ is absolutely continuous in $[0,2]$
and its distributional derivative is given by
$$
f'_{\lambda}(r)=\int_{\partial B_r}|\nabla w^{\lambda}|^2 dS
\quad\text{for a.e. }r\in(0,2).
$$
Suppose by contradiction that for any $M>0$ there exists a
sequence $\lambda_n\rightarrow 0^+$ such that
\begin{equation} \label{maggiore}
\int_{\partial B_r} |\nabla w^{\lambda_n}|^2 dS> M \int_{B_r} |\nabla
w^{\lambda_n}(x)|^2 dx \quad \text {for all }  r\in \left[1,2\right],
\end{equation}
which may be rewritten as
\begin{equation} \label{f1}
f'_{\lambda_n}(r)>M f_{\lambda_n}(r) \quad \text{for  a.e. }  r\in
[1,2]\text{ and for  any }  n\in \N.
\end{equation}
Integrating \eqref{f1} over $[1,2]$ we obtain
\begin{equation*}
f_{\lambda_n}(2)>e^{M} f_{\lambda_n}(1) \qquad {\rm for \ any} \
n\in \N.
\end{equation*}
Letting $n\rightarrow +\infty$ we obtain
\begin{equation*}
\limsup_{n\rightarrow +\infty} f_{\lambda_n}(1)\leq e^{-M} \cdot
\limsup_{n\rightarrow +\infty} f_{\lambda_n}(2).
\end{equation*}
This implies
\begin{equation*}
\liminf_{\lambda \rightarrow 0^+} f_{\lambda}(1) \leq e^{-M} \cdot
\limsup_{\lambda \rightarrow 0^+} f_{\lambda}(2) \qquad {\rm for \
any } \ M>0.
\end{equation*}
Using \eqref{limitata} and letting $M\rightarrow +\infty$ we infer
$$
\liminf_{\lambda\rightarrow 0^+} \int_{B_{1}} |\nabla
w^{\lambda}(x)|^2 dx= \liminf_{\lambda\rightarrow 0^+}
f_{\lambda}(1)=0.
$$
Therefore, there exists a sequence $\widetilde\lambda_n\rightarrow 0$ such
that
\begin{equation} \label{lim=0}
\lim_{n\rightarrow +\infty} \int_{B_1} |\nabla w^{\widetilde\lambda_n}(x)|^2
dx=0
\end{equation}
and, up to a subsequence still denoted by $\widetilde\lambda_n$, we
may suppose that $w^{\widetilde\lambda_n}\rightharpoonup w$ in
$H^1(B_1)$ for some $w\in H^1(B_1)$. Notice that, for any $\lambda\in
(0,\tilde r)$,  $\int_{\partial B_1} |w^\lambda|^2 dS=1 $ and hence by
compactness of the trace map from $H^1(B_1)$ into $L^2(\partial B_1)$,
it follows that $\int_{\partial B_1} |w|^2 dS=1$. Moreover, by weak
lower semicontinuity and \eqref{lim=0}, we also have
$$
\int_{B_1} |\nabla w(x)|^2 dx\leq \lim_{n\rightarrow +\infty}
\int_{B_1} |\nabla w^{\widetilde\lambda_n}(x)|^2 dx=0
$$
from which it follows that $w\equiv {\rm const}$ in $B_1$. On the other
hand, for every $\lambda\in (0,\tilde r)$,
\begin{equation} \label{eqlam} -\Delta
  w^\lambda(x)-\frac{a(\frac{x}{|x|})}{|x|^2}w^\lambda(x)=
  \lambda^2h(\lambda x)\,w^\lambda(x)+
  \frac{\lambda^2}{\sqrt{H(\lambda)}}f(\lambda
    x,\sqrt{H(\lambda)}w^\lambda(x)) \quad \text{ in } B_{\tilde
      r/\lambda}.
\end{equation}
For every $\phi\in H^1_0(B_1)$, by \eqref{F} and H\"older's inequality,
\begin{multline}\label{eq:39}
  \frac{\lambda^2}{\sqrt{H(\lambda)}} \left|\int_{B_1}f(\lambda
    x,\sqrt{H(\lambda)}w^\lambda(x))\phi(x)\,dx\right|\\
  \leq C_f\lambda^2\int_{B_1}|w^\lambda(x)||\phi(x)|\,dx
  +C_f\lambda^2\int_{B_1}|u(\lambda
  x)|^{2^*-2}|w^\lambda(x)||\phi(x)|\,dx
  \\
  \leq C_f\lambda^2\|w^\lambda\|_{H^1(B_1)}\|\phi\|_{H^1(B_1)}
 +C_f
\big({\textstyle{\frac{\omega_{N-1}}N}}
  \big)^{\frac{2(q-2^*)}{qN}}\lambda^{\frac{2(q-2^*)}{q}}
  \|w^\lambda\|_{L^{2^*}(B_1)} \|\phi\|_{L^{2^*}(B_1)}
  \|u\|_{L^{q}(B_\lambda)}^{2^*-2}\\[7pt]
=o(1)\quad\text{as }\lambda\to 0^+
\end{multline}
and, by \eqref{H1-2} and Corollary \ref{c:hardyboundary-0},
\begin{multline}\label{eq:41}
  \lambda^2\left|\int_{B_1}h(\lambda
    x)\,w^\lambda(x)\phi(x)\,dx\right| \\\leq C_h\lambda^\e \binom{N}{k}
  \bigg(1+\binom{N-k}{k}\bigg) \left(\frac{2}{k-2}\right)^{\!\!2}
  \bigg(\int_{B_1}|\nabla w^\lambda(x)|^2\,dx+\frac{N-2}2\bigg)^{\!\!1/2}
  \bigg(\int_{B_1}|\nabla \phi(x)|^2\,dx\bigg)^{\!\!1/2}
  \\[7pt]
  =o(1)\quad\text{as }\lambda\to 0^+.
\end{multline}
From (\ref{eq:39}), (\ref{eq:41}), and weak convergence
$w^{\widetilde\lambda_n}\rightharpoonup w$ in $H^1(B_1)$, we can pass to
the limit in \eqref{eqlam} along the sequence $\widetilde\lambda_n$
and obtain  that $w$ is a $H^1(B_1)$-weak solution to the equation
$$
-\Delta w(x)-\frac{a(\frac{x}{|x|})}{|x|^2}w(x)=0 \quad \text{in }
\ B_1.
$$
Since $w$ is constant in $B_1$, this implies $w\equiv 0$ in $B_1$
which contradicts $\int_{\partial B_1} |w|^2 dS=1$.
\end{pf}

\begin{Lemma} \label{nablabound}
Let $w^\lambda$ and $R_\lambda$ be as in the statement of Lemma
\ref{Rlambda}. Then there exists  $\overline M>0$ such that
$$
\int_{\partial B_1} |\nabla w^{\lambda R_\lambda}|^2 dS\leq
\overline M \qquad {\rm for \ any } \ 0<\lambda<
\min\Big\{\lambda_0,\frac{\bar r}2\Big\}.
$$
\end{Lemma}

\begin{pf} We have
\begin{align*}
\int_{\partial B_1} |\nabla w^{\lambda R_{\lambda}}|^2 dS
=&\frac{(\lambda R_\lambda)^{2}}{H(\lambda R_\lambda)}
\int_{\partial B_1} |\nabla u(\lambda R_{\lambda} x)|^2 dS(x)
=\frac{\lambda^2 R_\lambda^{3-N}}{H(\lambda R_\lambda)}
\int_{\partial B_{R_\lambda}} |\nabla u(\lambda x)|^2 dS(x)\\
=& \frac{R_\lambda^{3-N}H(\lambda)}{H(\lambda
R_\lambda)}\frac{\lambda^2}{H(\lambda)} \int_{\partial
B_{R_\lambda}} |\nabla u(\lambda x)|^2 dS(x) =
\frac{R_\lambda^{3-N}H(\lambda)}{H(\lambda R_\lambda)}
\int_{\partial B_{R_\lambda}} |\nabla w^{\lambda}|^2 dS
\end{align*}
and, by (\ref{dou1}--\ref{dou2}), Lemma \ref{Rlambda}, Lemma
\ref{bound},
and the fact
that $1\leq R_\lambda\leq 2$, we finally obtain
\begin{align*}
\int_{\partial B_1} |\nabla w^{\lambda R_{\lambda}}|^2 dS \leq C_4
M \int_{B_{R_\lambda}} |\nabla w^{\lambda}(x)|^2 dx \leq 2^{N-2}
C_4^2 M \int_{B_1} |\nabla w^{\lambda R_{\lambda}}(x)|^2 dx\leq
\overline M<+\infty
\end{align*}
for any $0<\lambda<\min\big\{\lambda_0,\frac{\bar r}2\big\}$, thus
completing the proof.
\end{pf}

\begin{Lemma}\label{l:blow-up}
  Let $u$ be a weak $H^1(\Omega)$-solution to (\ref{u}), $u\not\equiv
  0$, in a bounded open set $\Omega\subset\R^N$, $N\geq 3$, with $a$
  satisfying (\ref{eq:ateta}) and (\ref{eq:lambdamin1}), $h$
  satisfying \eqref{H1-2}, and $f$ satisfying \eqref{F}. Let $\gamma$
  be as in Lemma \ref{gamma}.  Then
\begin{itemize}
\item[\rm (i)] there exists $k_0\in \N$ such that
  $\gamma=-\frac{N-2}2+\sqrt{\big(\frac{N-2}{2}\big)^{\!2}+\mu_{k_0}(a)}$;
\item[\rm (ii)] for every sequence $\lambda_n\to0^+$, there exist
a subsequence $\{\lambda_{n_k}\}_{k\in\N}$ and an eigenfunction
$\psi$ of the operator $-\Delta_{\SN}-a(\theta)$ associated to the
eigenvalue $\mu_{k_0}(a)$ such that $\|\psi\|_{L^{2}({\mathbb
S}^{N-1})}=1$ and
\[
\frac{u(\lambda_{n_k}x)}{\sqrt{H(\lambda_{n_k})}}\to
|x|^{\gamma}\psi\Big(\frac x{|x|}\Big)
\]
strongly in $H^1(B_1)$.
\end{itemize}
\end{Lemma}

\begin{pf}
  Let $\lambda_n \rightarrow 0^+$ and consider the sequence
  $w^{\lambda_n R_{\lambda_n}}$ as in (\ref{eq:defwlambda}) and
  $R_\lambda$ as in Lemma \ref{Rlambda}. By Lemmas \ref{bound} and
  \ref{doubling}, we have that the set $\{w^{\lambda
    R_\lambda}\}_{\lambda\in (0,\bar r/4)}$ is bounded in $H^1(B_2)$.
  Then there exists a subsequence $w^{\lambda_{n_k}
    R_{\lambda_{n_k}}}$ such that $w^{\lambda_{n_k}
    R_{\lambda_{n_k}}}\rightharpoonup w$ in $H^1(B_2)$ for some
  function $w\in H^1(B_2)$.  Due to compactness of the trace map from
  $H^1(B_1)$ into $L^2(\partial B_1)$, we obtain that $\int_{\partial
    B_1}|w|^2dS=1$. In particular $w\not\equiv 0$. Furthermore, weak
  convergence and (\ref{eq:39}--\ref{eq:41}) allow passing to the
  weak limit in the equation
\begin{multline}\label{eq:19}
  -\Delta
  w^{\lambda_{n_k}R_{\lambda_{n_k}}}(x)-\frac{a(\frac{x}{|x|})}{|x|^2}
  w^{\lambda_{n_k}R_{\lambda_{n_k}}}(x)=
  \lambda_{n_k}^2R_{\lambda_{n_k}}^2
  h(\lambda_{n_k}R_{\lambda_{n_k}} x)\,w^{\lambda_{n_k}R_{\lambda_{n_k}}}(x)\\+
  \frac{
    \lambda_{n_k}^2R_{\lambda_{n_k}}^2}{\sqrt{H({\lambda_{n_k}R_{\lambda_{n_k}}})}}
  f\Big({\lambda_{n_k}R_{\lambda_{n_k}}}
  x,\sqrt{H({\lambda_{n_k}R_{\lambda_{n_k}}})}\,w^{\lambda_{n_k}R_{\lambda_{n_k}}}(x)\Big)
\end{multline}
which holds in a weak sense in
$B_{\tilde r/(\lambda_{n_k}R_{\lambda_{n_k}})}\supset B_2$
thus yielding
\begin{equation} \label{eq:w}
-\Delta
  w(x)-\frac{a(\frac{x}{|x|})}{|x|^2}w(x)=0\quad\text{in }B_2.
\end{equation}
From Lemma \ref{nablabound} and density in $H^1(B_1)$ of
$C^{\infty}(\overline{B}_1)$-functions whose support is compactly
included in $\overline{B}_1\setminus\widetilde\Sigma$ with
$\widetilde\Sigma$ defined in (\ref{eq:sigmatilde}), it follows that,
for all $\phi\in H^1(B_1)$,
\begin{align}\label{eq:parti}
  \int_{B_1}&\bigg(\nabla
  w^{\lambda_{n_k}R_{\lambda_{n_k}}}(x)\cdot\nabla\phi(x)-
  \frac{a(\frac{x}{|x|})}{|x|^2}
  w^{\lambda_{n_k}R_{\lambda_{n_k}}}(x)\phi(x)
  \bigg)\,dx \\
  &\notag\quad =
  \lambda_{n_k}^2R_{\lambda_{n_k}}^2\int_{B_1}h(\lambda_{n_k}R_{\lambda_{n_k}}x)
  \,w^{\lambda_{n_k}R_{\lambda_{n_k}}}(x)\phi(x)\,dx+\int_{\partial
    B_1}\frac{\partial w^{\lambda_{n_k}R_{\lambda_{n_k}}}
  }{\partial\nu}\,\phi\,dS\\
  &\notag\qquad + \frac{ \lambda_{n_k}^2R_{\lambda_{n_k}}^2}
  {\sqrt{H({\lambda_{n_k}R_{\lambda_{n_k}}})}}
  \int_{B_1}f\Big({\lambda_{n_k}R_{\lambda_{n_k}}}
  x,\sqrt{H({\lambda_{n_k}R_{\lambda_{n_k}}})}
  \,w^{\lambda_{n_k}R_{\lambda_{n_k}}}(x)\Big) \phi(x)\,dx.
\end{align}
We notice that from (\ref{Nbelow}) it follows that
$\gamma\geq-\frac{N-2}2$. Then, by \eqref{F} and (\ref{1stest}),
\begin{align*}
  \frac{ \lambda^2}{\sqrt{H({\lambda })}} \left|\frac{f\big({\lambda }
    x,\sqrt{H({\lambda })}\,w^{\lambda }(x)\big)}{w^{\lambda }(x)}
 \right| & \leq
  C_f\frac{ \lambda^2}{\sqrt{H({\lambda })}} \Big(\sqrt{H({\lambda
    })}+
  (H({\lambda }))^{\frac{2^*-1}2}|w^{\lambda}(x)|^{2^*-2}\Big)\\
  &\notag\leq C_f \Big( \lambda^2+K_1^{\frac{2^*-2}2}
  \lambda^{2+\gamma(2^*-2)} |w^{\lambda }(x)|^{2^*-2}\Big)
  \\
  &\notag\leq {\rm const\,}\Big(1+|w^{\lambda }(x)|^{2^*-2}\Big)
\end{align*}
for all $\lambda\in(0,\tilde r)$.  Hence, if $s=\frac{q}{2^*-2}>N/2$
with $q$ as in (\ref{eq:46}), from (\ref{eq:defwlambda}) and 
Proposition \ref{SMETS}, we obtain that
\begin{align*}
  \left\| \frac{ \lambda^2}{\sqrt{H({\lambda })}} \frac{f\big({\lambda
      } x,\sqrt{H({\lambda })}\,w^{\lambda }(x)\big)}{w^{\lambda }(x)}
  \right\|_{L^s(B_2)} & \leq {\rm const\,}
  \bigg(1+\lambda^2(H({\lambda }))^{\frac{2^*-2}2} \bigg(\int_{B_2}|w^{\lambda
  }(x)|^{(2^*-2)s}dx\bigg)^{\!\!1/s}\bigg)\\
&= {\rm const\,}
  \bigg(1+\lambda^{2-\frac Ns}
 \bigg(\int_{B_{2\lambda}}|u(x)|^{q}dx\bigg)^{\!\!1/s}\bigg)=O(1)
\end{align*}
as $\lambda\to 0^+$.  Therefore from classical Brezis-Kato \cite{BK}
estimates (see also Theorem 8.6 part i)), classical bootstrap and
elliptic regularity theory, there holds
$$
w^{\lambda_{n_k}R_{\lambda_{n_k}}}\to w\quad\text{in }C^{1,\tau}_{\rm
  loc}(B_2\setminus \widetilde \Sigma),
$$
for any $\tau\in(0,1)$, which in particular yields
\begin{equation}\label{eq:convnormal}
  \frac{\partial w^{\lambda_{n_k}R_{\lambda_{n_k}}} }{\partial\nu} \to
  \frac{\partial w}{\partial\nu}\quad\text{in }C^{0,\tau}_{\rm
    loc}(\partial B_1\setminus \Sigma)\quad\text{and a.e. in }\partial B_1.
\end{equation}
From (\ref{eq:convnormal}) and Lemma \ref{nablabound}, it follows that
\begin{equation}\label{eq:convl2normal}
  \frac{\partial w^{\lambda_{n_k}R_{\lambda_{n_k}}} }{\partial\nu}
  \rightharpoonup
  \frac{\partial w}{\partial\nu}\quad\text{weakly in }L^2(\partial B_1).
\end{equation}
Passing to limit in (\ref{eq:parti}) and using (\ref{eq:convl2normal}) and
 (\ref{eq:39}--\ref{eq:41}), we
obtain that
\begin{align}\label{eq:partilim}
  \int_{B_1}&\bigg(\nabla w(x)\cdot\nabla\phi(x) -
  \frac{a(\frac{x}{|x|})}{|x|^2}w(x)\phi(x)\bigg)\,dx =\int_{\partial
    B_1}\frac{\partial w }{\partial\nu}\,\phi\,dS.
\end{align}
Subtracting (\ref{eq:partilim}) from (\ref{eq:parti}), choosing
$\phi=w^{\lambda_{n_k}R_{\lambda_{n_k}}}-w$, and arguing as in 
 (\ref{eq:39}--\ref{eq:41}) and
Corollary \ref{c:hardyboundary}, we obtain that
\begin{align}\label{eq:limnorm}
w^{\lambda_{n_k}R_{\lambda_{n_k}}}\to w
\quad\text{in }H^1(B_1).
\end{align}
For every $k\in\N$ and $r\in(0,1)$, let us define
\begin{multline*}
  D_k(r)\\
=\frac{1}{r^{N-2}} \int_{B_r}
  \bigg[\left|\nabla w^{\lambda_{n_k}R_{\lambda_{n_k}}}(x)\right|^2
-\frac{a(\frac{x}{|x|})}{|x|^2}|w^{\lambda_{n_k}R_{\lambda_{n_k}}}(x)|^2
-\lambda_{n_k}^2R_{\lambda_{n_k}}^2h(\lambda_{n_k}R_{\lambda_{n_k}}x)
|w^{\lambda_{n_k}R_{\lambda_{n_k}}}(x)|^2 \\
  -\frac{\lambda_{n_k}^2R_{\lambda_{n_k}}^2}{\sqrt{H(\lambda_{n_k} R_{\lambda_{n_k}})}}
 f\Big(\lambda_{n_k}R_{\lambda_{n_k}}
  x,\sqrt{H(\lambda_{n_k}R_{\lambda_{n_k}})}w^{\lambda_{n_k}R_{\lambda_{n_k}}}(x)\Big)
 w^{\lambda_{n_k}R_{\lambda_{n_k}}}(x)\bigg] \, dx
\end{multline*}
and
\begin{equation*}
H_k(r)=\frac{1}{r^{N-1}}\int_{\partial B_r}|w^{\lambda_{n_k}R_{\lambda_{n_k}}}|^2 \, dS.
\end{equation*}
We also define
\begin{equation} \label{Dw(r)}
  D_w(r)=\frac{1}{r^{N-2}} \int_{B_r}
  \left[\left|\nabla w(x)\right|^2
-\frac{a(\frac{x}{|x|})}{|x|^2}|w(x)|^2
  \right] \, dx \quad \text{for all } r\in (0,1)
\end{equation}
and
\begin{equation} \label{Hw(r)}
H_w(r)=\frac{1}{r^{N-1}}\int_{\partial B_r}|w|^2 \, dS
\quad  \text{for all } r\in (0,1).
\end{equation}
Using a change of variables, one sees that
\begin{equation}\label{NkNw}
{\mathcal
    N}_k(r):=\frac{D_k(r)}{H_k(r)}=\frac{D(\lambda_{n_k}R_{\lambda_{n_k}}r)}
{H(\lambda_{n_k}R_{\lambda_{n_k}}r)}
  ={\mathcal N}(\lambda_{n_k}R_{\lambda_{n_k}}r) \quad \text{for all } r\in (0,1).
\end{equation}
From (\ref{eq:7}), (\ref{eq:26}), and (\ref{eq:limnorm}), it follows that,
for any fixed $r\in (0,1)$,
\begin{equation} \label{convDk}
 D_k(r)\to D_w(r).
\end{equation}
On the other hand, by compactness of the trace embedding
$H^1(B_r)\hookrightarrow L^2(\partial B_r)$, we also have,
for any fixed $r\in (0,1)$,
\begin{equation} \label{convHk}
H_k(r)\to H_w(r).
\end{equation}
From Lemma \ref{l:hardyboundary} and classical unique continuation
principle for second order elliptic equations with locally bounded
coefficients (see e.g.  \cite{wolff}) applied away from the singular
set $\widetilde\Sigma$, it follows that $D_w(r)>-\frac{N-2}2H_w(r)$
for all $r\in(0,1)$.  Therefore, if for some $r\in(0,1)$, $H_w(r)=0$
then $D_w(r)>0$; passing to the limit in (\ref{NkNw}) and using
\eqref{convDk}-\eqref{convHk} this should give a contradiction with
Lemma \ref{gamma}. Hence $H_w(r)>0$ for all $r\in(0,1)$. Thus the
function
\[
{\mathcal N}_w(r):=\frac{D_w(r)}{H_w(r)}
\]
is well defined for $r\in (0,1)$.
This, together with (\ref{NkNw}), (\ref{convDk}), (\ref{convHk}), and
Lemma \ref{gamma}, shows that
\begin{equation} \label{Nw(r)} {\mathcal
    N}_w(r)=\lim_{k\to \infty} {\mathcal
    N}(\lambda_{n_k}R_{\lambda_{n_k}}r)=\gamma
\end{equation}
for all $r\in (0,1)$.  Therefore ${\mathcal N}_w$ is constant in
$(0,1)$ and hence ${\mathcal N}_w'(r)=0$ for any $r\in (0,1)$. Hence,
by (\ref{eq:w}) and Lemma \ref{mono} with $h\equiv0$, $f\equiv 0$, we obtain
\begin{equation*}
  \left(\int_{\partial B_r} \left|\frac{\partial
        w}{\partial\nu}\right|^2 dS\right) \cdot \left(\int_{\partial B_r} |w|^2
    dS\right)-\left(\int_{\partial B_r} w\frac{\partial w}{\partial \nu}\,
 dS\right)^{\!\!2}=0 \quad
  \text{for a.e. } r\in (0,1).
\end{equation*}
This shows that $w$ and $\frac{\partial w}{\partial \nu}$ have the
same direction as vectors in $L^2(\partial B_r)$ and hence there
exists $\eta=\eta(r)$ such that
\begin{equation}\label{eq:44}
\frac{\partial w}{\partial
\nu}(r,\theta)=\eta(r) w(r,\theta)
\quad\text{for a.e. $r\in(0,1)$, $\theta\in {\mathbb S}^{N-1}$}.
\end{equation}
Testing the above identity with $w(r,\theta)$, we have that
necessarily $\eta(r)=\frac{H_w'(r)}{2H_w(r)}$ implying that $\eta\in
L^1_{\rm loc}(0,1)$.  Moreover, since $w\in C^{1}_{\rm
  loc}(B_2\setminus \widetilde \Sigma)$, identity (\ref{eq:44}) also
holds, for all $\theta \in {\mathbb S}^{N-1}\setminus\Sigma$, in the sense
of absolutely continuous functions with respect to $r$ and,
after integration, we obtain
\begin{equation} \label{separate}
w(r,\theta)=e^{\int_1^r \eta(s)ds} w(1,\theta)
=\varphi(r) \psi(\theta) \quad  \text{for all }r\in(0,1), \ \theta\in
\SN\setminus\Sigma,
\end{equation}
where  $\varphi(r)=e^{\int_1^r \eta(s)ds}$ and $\psi(\theta)=w(1,\theta)$.
Since
$$
\mathcal -\Delta
w-\frac{a(\frac{x}{|x|})}{|x|^2}w=-\frac{\partial^2 w}{\partial
  r^2}-\frac{N-1}{r}\frac{\partial w}{\partial r} +\frac{1}{r^2}
L_aw,
$$
then (\ref{separate}) yields
$$
\left(-\varphi''(r)-\frac{N-1}{r} \varphi'(r) \right)\psi(\theta)
+\frac{\varphi(r)}{r^2} L_a\psi(\theta)=0.
$$
Taking $r$ fixed we deduce that $\psi$ is an eigenfunction of the
operator $L_a$. If $\mu_{k_0}(a)$ is the corresponding eigenvalue then
$\varphi(r)$ solves the equation
$$
-\varphi''(r)-\frac{N-1}{r}
\varphi(r)+\frac{\mu_{k_0}(a)}{r^2}\varphi(r)=0
$$
and hence $\varphi(r)$ is of the form
$$
\varphi(r)=c_1 r^{\sigma_{k_0}^+}+c_2 r^{\sigma_{k_0}^-}
$$
for some $c_1,c_2\in\R$, where
\begin{equation*}
  \sigma^+_{k_0}=-\frac{N-2}{2}+\sqrt{\bigg(\frac{N-2}
    {2}\bigg)^{\!\!2}+\mu_{k_0}(a)}\quad\text{and}\quad
  \sigma^-_{k_0}=-\frac{N-2}{2}-\sqrt{\bigg(\frac{N-2}{2}
    \bigg)^{\!\!2}+\mu_{k_0}(a)}.
\end{equation*}
Since the function
$\frac1{|x_J|}\big(|x|^{\sigma_{k_0}^-}\psi(\frac{x}{|x|})\big)\notin
L^2(B_1)$ for any $J\in \mathcal A_k$ and hence
$|x|^{\sigma_{k_0}^-}\psi(\frac{x}{|x|})\notin H^1(B_1)$, then
$c_2=0$ and  $\varphi(r)=c_1 r^{\sigma_{k_0}^+}$. Since
$\varphi(1)=1$, we obtain that $c_1=1$ and then
\begin{equation} \label{expw}
w(r,\theta)=r^{\sigma_{k_0}^+} \psi(\theta),  \quad
\text{for all }r\in (0,1)\text{ and }\theta\in \SN\setminus\Sigma.
\end{equation}
Consider now the sequence $w^{\lambda_{n_k}}$. Up to a further
subsequence still denoted by $w^{\lambda_{n_k}}$, we may suppose
that $w^{\lambda_{n_k}}\rightharpoonup \overline w$ for some
$\overline w\in H^1(B_1)$ and that $R_{\lambda_{n_k}}\rightarrow
\overline R$ for some $\overline R\in [1,2]$.

Strong convergence of $w^{\lambda_{n_k}R_{\lambda_{n_k}}}$  in
$H^1(B_1)$ implies that, up to a subsequence, both
$w^{\lambda_{n_k}R_{\lambda_{n_k}}}$ and
$|\nabla w^{\lambda_{n_k}R_{\lambda_{n_k}}}|$ are dominated by a
$L^2(B_1)$-function uniformly with respect to $k$. Moreover by
\eqref{dou1}, up to a subsequence we may assume that the limit
$$
l:=\lim_{k\rightarrow+\infty}
\frac{H(\lambda_{n_k}R_{\lambda_{n_k}})}{H(\lambda_{n_k})}
$$
exists and is finite.
Then, by the Dominated Convergence Theorem, we have
\begin{align*}
  & \lim_{k\rightarrow+\infty} \int_{B_1} w^{\lambda_{n_k}}(x)
  v(x) \, dx =\lim_{k\rightarrow +\infty}
  R_{\lambda_{n_k}}^N\int_{B_{1/R_{\lambda_{n_k}}}}
  w^{\lambda_{n_k}}(R_{\lambda_{n_k}}
  x)v(R_{\lambda_{n_k}}x)\, dx \\
  &=\lim_{k\rightarrow +\infty} R_{\lambda_{n_k}}^N
  \sqrt{\frac{H(\lambda_{n_k}R_{\lambda_{n_k}})}{H(\lambda_{n_k})}}
  \int_{B_1} \alchi_{B_{1/R_{\lambda_{n_k}}}}\!\!\!(x)\,
  w^{\lambda_{n_k}R_{\lambda_{n_k}}}(x)v(R_{\lambda_{n_k}}x)\,
  dx
  \\
  &=\overline R^N \sqrt l \int_{B_1} \alchi_{B_{1/\overline
      R}}(x)w(x) v(\overline R x) \, dx = \overline R^N
  \sqrt l \int_{B_{1/\overline R}} w(x) v(\overline R x) \,
  dx =\sqrt l \int_{B_1} w(x/\overline R)v(x) \, dx
\end{align*}
for any $v\in C^\infty(\R^N)$ with ${\rm supp} \,
v\subset B_1$. By a density argument, it follows that the
previous convergence also holds for all $v\in L^2(B_1)$.
This proves that $w^{\lambda_{n_k}}\rightharpoonup \sqrt l \,
w(\cdot/\overline R)$ in $L^2(B_1)$ (actually weakly in
$H^1(B_1)$) and in particular $\overline w=\sqrt l \,
w(\cdot/\overline R)$. Moreover
\begin{align*}
& \lim_{k\rightarrow+\infty} \int_{B_1} |\nabla
w^{\lambda_{n_k}}(x)|^2 \, dx =\lim_{k\rightarrow +\infty}
 R_{\lambda_{n_k}}^N\int_{B_{1/R_{\lambda_{n_k}}}}
|\nabla w^{\lambda_{n_k}}(R_{\lambda_{n_k}}
x)|^2 dx \\
&=\lim_{k\rightarrow +\infty} R_{\lambda_{n_k}}^{N-2}
\frac{H(\lambda_{n_k}R_{\lambda_{n_k}})}{H(\lambda_{n_k})}
\int_{B_1} \chi_{B_{1/R_{\lambda_{n_k}}}} |\nabla
w^{\lambda_{n_k}R_{\lambda_{n_k}}}(x)|^2 \, dx
\\
&=\overline R^{N-2} l \int_{B_1} \chi_{B_{1/\overline
R}}(x)|\nabla w(x)|^2 \, dx = \overline R^{N-2} l
\int_{B_{1/\overline R}} |\nabla w(x)|^2 \, dx = \int_{B_1}
|\sqrt l \,\nabla (w(x/\overline R))|^2 \, dx.
\end{align*}
This shows that $w^{\lambda_{n_k}}\to \overline w=\sqrt l \,
w(\cdot/\overline R)$ strongly in $H^1(B_1)$.
Furthermore, by \eqref{expw} and
the fact that $\int_{\partial B_1} |\overline w|^2
dS=\int_{\partial B_1} |w|^2 dS=1$, we deduce that $\overline w=w$.

It remains to prove part (i). By (\ref{expw}) and $\int_{{\mathbb
    S}^{N-1}}|\psi(\theta)|^2dS=1$
 we have that
\begin{multline*}
\int_{B_r} \bigg(|\nabla w(x)|^2-\frac{a(\frac{x}{|x|})}{|x|^2}|w(x)|^2
\bigg) \, dx=(\sigma_{k_0}^+)^2\int_0^r s^{N-1+2(\sigma_{k_0}^+-1)}\,ds\\+
\bigg(\int_0^r s^{N-1+2(\sigma_{k_0}^+-1)}\,ds\bigg)
\bigg(\int_{{\mathbb S}^{N-1}}\big(|\nabla_{\mathbb S^{N-1}}\psi(\theta)|^2-a(\theta)
|\psi(\theta)|^2\big)dS\bigg)\\
=\frac{(\sigma_{k_0}^+)^2+\mu_{k_0}(a)}{N+2(\sigma_{k_0}^+-1)}
r^{N+2(\sigma_{k_0}^+-1)}=\sigma_{k_0}^+r^{N+2(\sigma_{k_0}^+-1)}
\end{multline*}
and
\begin{multline*}
\int_{\partial B_r} |w(x)|^2dS=
r^{N-1}\int_{{\mathbb S}^{N-1}}|w(r\theta)|^2dS
=r^{N-1+2\sigma_{k_0}^+}\int_{{\mathbb
    S}^{N-1}}|\psi(\theta)|^2dS=r^{N-1+2\sigma_{k_0}^+}.
\end{multline*}
Therefore, by (\ref{Dw(r)}), (\ref{Hw(r)}), and (\ref{Nw(r)}), it follows
$$
\gamma={\mathcal N}_w(r)=\frac{D_w(r)}{H_w(r)}=\frac{r
 \int_{B_r}\big(|\nabla w(x)|^2-\frac{a(x/|x|)}{|x|^2}|w(x)|^2
\big) \, dx}{\int_{\partial B_r} |w|^2 dS } =\sigma_{k_0}^+.
$$
This completes the proof of the lemma.
\end{pf}

Let us now describe the behavior of $H(r)$ as $r\to 0^+$.
\begin{Lemma} \label{l:limite}
Under the same assumptions as in Lemma
  \ref{mono} and letting $\gamma:=\lim_{r\rightarrow 0^+} {\mathcal
    N}(r)\in \R$ as in Lemma \ref{gamma}, the limit
\[
\lim_{r\to0^+}r^{-2\gamma}H(r)
\]
exists and it is finite.
\end{Lemma}
\begin{pf}
In view of (\ref{1stest}) it is sufficient to prove that the limit
exists. By (\ref{H(r)}), (\ref{ul}), and Lemma~\ref{gamma} we have
$$
\frac{d}{dr} \frac{H(r)}{r^{2\gamma}} =-2\gamma r^{-2\gamma-1}
H(r)+r^{-2\gamma} H'(r) =2r^{-2\gamma-1} (D(r)-\gamma
H(r))=2r^{-2\gamma-1} H(r) \int_0^r {\mathcal N}'(s) ds.
$$
Let $\nu_1$ and $\nu_2$ be as in (\ref{eq:nu1}) and (\ref{eq:nu2}).
After integration over $(r,\tilde r)$,
\begin{equation}\label{inte}
  \frac{H(\tilde r)}{\tilde r^{2\gamma}}-
  \frac{H(r)}{r^{2\gamma}}=\int_r^{\tilde r} 2s^{-2\gamma-1}
  H(s) \left( \int_0^s \nu_1(t) dt \right) ds +\int_r^{\tilde r} 2s^{-2\gamma-1}
  H(s) \left( \int_0^s \nu_2(t) dt \right) ds.
\end{equation}
By Schwarz's inequality we have that $\nu_1(t)\geq 0$ and hence
$$
\lim_{r\to 0^+} \int_r^{\tilde r} 2s^{-2\gamma-1} H(s) \left( \int_0^s
  \nu_1(t) dt \right) ds
$$
exists.  On the other hand, by
 (\ref{1stest}), Lemma \ref{l:stima_nu2}, and (\ref{Nabove}), we
deduce that
\begin{multline*}
  \left| s^{-2\gamma-1} H(s) \left( \int_0^s \nu_2(t) dt \right)
  \right|\leq K_1C_1\Big(C_2+\frac N2\Big)s^{-1}\int_0^s\Big(
  t^{-1+\e}+t^{-1+\frac{2(q-2^*)}{q}}+g(t)\Big)\,dt\\
  \leq K_1C_1\Big(C_2+\frac N2\Big)s^{-1}\bigg(
  \frac{s^\e}{\e}
+\frac{q}{2(q-2^*)}
s^{\frac{2(q-2^*)}{q}}+\frac {\|u\|_{L^{2^*}(\Omega)}^{2^*(1-\alpha)}}{1-\alpha}
\, s^{\frac{N(q-2^*)}q(\alpha-\frac{2}{2^*})}
\bigg)
\end{multline*}
for all $s\in(0,\widetilde r)$, which proves that $s^{-2\gamma-1}
H(s) \left( \int_0^s \nu_2(t) dt \right)\in L^1(0,\widetilde r)$.
We may conclude that both terms in the right hand side of
(\ref{inte}) admit a limit as $r\to 0^+$ thus completing the proof
of the lemma.
\end{pf}

The next step of our asymptotic analysis relies on the proof that
$\lim_{r\to 0^+} r^{-2\gamma} H(r)$ is indeed strictly positive.
In the sequel we denote by $\psi_i$ a $L^2$-normalized eigenfunction
of the operator $L_a=-\Delta_{\SN}-a$ associated to the $i$-th eigenvalue
$\mu_i(a)$, i.e.
\begin{equation} \label{eq:2rad}
\begin{cases}
L_a\psi_i(\theta)
=\mu_i(a) \,\psi_i(\theta),&\text{in }{\mathbb S}^{N-1},\\[3pt]
\int_{{\mathbb S}^{N-1}}|\psi_i(\theta)|^2\,dS(\theta)=1.
\end{cases}
\end{equation}
Moreover, we  choose the $\psi_i$'s in such a way that
 the set $\{\psi_i\}_{i\in \N}$
forms an orthonormal basis of $L^2(\SN)$.

Let $u$ be a nontrivial weak $H^1(\Omega)$-solution to (\ref{u}).
From Lemma \ref{l:blow-up}, we deduce that, under assumptions
(\ref{eq:ateta}), (\ref{eq:lambdamin1}), and (\ref{H1-2}--\ref{F}), there exist
$j_0,m\in\N$, $j_0,m\geq 1$ such that $m$ is the multiplicity of the
eigenvalue $\mu_{j_0}(a)=\mu_{j_0+1}(a)=\cdots=\mu_{j_0+m-1}(a)$ and
\begin{equation}\label{eq:30}
  \gamma=\lim_{r\rightarrow 0^+} {\mathcal N}(r)=
-\frac{N-2}{2}+\sqrt{\bigg(\frac{N-2}
    {2}\bigg)^{\!\!2}+\mu_{i}(a)},
  \quad i=j_0,\dots,j_0+m-1.
\end{equation}
Let ${\mathcal E}_0$ be the eigenspace
of the operator $L_{a}$ associated to the eigenvalue
$\mu_{j_0}(a)$, so that the set $\{\psi_i\}_{i=j_0,\dots,j_0+m-1}$ is an
 orthonormal basis of ${\mathcal E}_0$.

\begin{Lemma}\label{l:stima}
  Let $\Omega\subset\R^N$, $N\geq 3$, be a bounded open set such that
  $0\in\Omega$,  $a$ such that (\ref{eq:ateta}) and (\ref{eq:lambdamin1})
hold, and
  $h,f$ as in  (\ref{H1-2}--\ref{F}). If $u$  is a  weak
  $H^1(\Omega)$-solution to (\ref{u}), then
$$
\sup_{\substack{ i=j_0,\dots,j_0+m-1,\
    J\in{\mathcal A}_k,\\
    (J_1,J_2)\in {\mathcal B}_k,\ \lambda\in(0,\bar r) }}\frac{
  {\displaystyle{\int_{B_\lambda}}}\bigg(\dfrac{|u(x)|
}{|x_J|^{2-\e}}+
\dfrac{|u(x)|}{|x_{J_1}-x_{J_2}|^{2-\e}}+
|f(x,u(x))|\bigg)
\big|\psi_i(\frac{x}{|x|})\big|
\,dx}
{{\raisebox{-3pt}{$\lambda^{N-2+\delta+\gamma}$}}} <+\infty,
$$
where $\bar r$ is as in Lemma \ref{bound} and $\delta>0$ is defined in
(\ref{eq:45}).
\end{Lemma}
\begin{pf}
From Lemma \ref{bound} and Corollary \ref{c:hardyboundary-0},
it follows that, for some
positive constant $C_5$ independent of $\lambda$, $J$, $(J_1,J_2)$, and $i$,
$$
\int_{B_1 }\Big(|x_J|^{-2+\e}+|x_{J_1}-x_{J_2}|^{-2+\e}\Big)
|w^\lambda(x)|
\big|\psi_i\big({\textstyle{\frac{x}{|x|}}}\big)\big|\,dx\leq C_5
$$
for all $i=j_0,\dots,j_0+m-1$, $J\in{\mathcal A}_k$, $(J_1,J_2)\in
{\mathcal B}_k$, and $\lambda\in(0,\bar r)$, where $w^\lambda$ is
defined in (\ref{eq:defwlambda}). Moreover, arguing as in
(\ref{eq:39}), by (\ref{eq:45}), we have as $\lambda\to 0^+$
\begin{align*}
  \frac{\lambda^2}{\sqrt{H(\lambda)}}
\int_{B_1}\left|f(\lambda
    x,u(\lambda x))\right|
\big|\psi_i\big({\textstyle{\frac{x}{|x|}}}\big)\big|\,dx \leq C_6
\lambda^{\frac{2(q-2^*)}q}=O(\lambda^\delta)
\end{align*}
where $C_6>0$ is a positive constant. The conclusion follows from
(\ref{1stest}) and a change of variable.~\end{pf}

\begin{Lemma} \label{l:phii} Let $\Omega\subset\R^N$, $N\geq
  3$, be a bounded open set such that $0\in\Omega$, $a$ satisfy
  (\ref{eq:ateta}) and (\ref{eq:lambdamin1}), and $u\not\equiv 0$ be a
  weak $H^1(\Omega)$-solution to (\ref{u}), with $h,f$ satisfying
  {\rm (\ref{H1-2}--\ref{F})}. Let
$\gamma:=\lim_{r\rightarrow 0^+} {\mathcal
    N}(r)$ be as in Lemma \ref{gamma} and
$j_0,m\in\N$ as in
(\ref{eq:30}), i.e.
 $m$ is the
multiplicity of the eigenvalue
$\mu_{j_0}(a)=\mu_{j_0+1}(a)=\cdots=\mu_{j_0+m-1}(a)$ and
(\ref{eq:30}) holds for all $i=j_0,\dots,j_0+m-1$.
Then the function $\varphi_i$ defined as
\begin{equation} \label{78}
  \varphi_i(\lambda):=\int_{{\mathbb S}^{N-1}}u(\lambda\,\theta)
  \psi_i(\theta)\,dS(\theta),
\quad\text{with }\psi_i\text{ as in (\ref{eq:2rad})},
\end{equation}
satisfies, as $\lambda\to0^+$ ,
\begin{multline}\label{eq:35}
  \varphi_i(\lambda)=\lambda^{\gamma} \bigg( R^{-\gamma}\varphi_i(R)+
  \frac{2-N-\gamma}{2-N-2\gamma} \int_\lambda^Rs^{-N+1-\gamma}\Upsilon_i(s)ds\\
  - \frac{\gamma R^{-N+2-2\gamma}}{2-N-2\gamma} \int_0^R
  s^{\gamma-1}\Upsilon_i(s)\,ds \bigg)+O(\lambda^{\gamma+\delta})
\end{multline}
for every $i\in \{j_0,\dots,j_0+m-1\}$  and
$R>0$ such that $\overline{B_R}\subset \Omega$,
where $\delta$ is defined in
(\ref{eq:45}) and
\begin{equation}\label{eq:fi}
\Upsilon_i(\lambda):=
\int_{B_\lambda}\Big(h(x)u(x)
+f(x,u(x))\Big)\psi_i\bigg(\frac{x}{|x|}\bigg)\,dx.
\end{equation}
\end{Lemma}
\begin{pf}
Let $R>0$ be such that $\overline{B_R}\subset \Omega$.
For any $\lambda\in (0,R)$, we 
expand $\theta\mapsto u(\lambda\theta)\in L^2({\mathbb S}^{N-1})$
 in  Fourier series with respect to the orthonormal
basis $\{\psi_i\}$ of $L^2({\mathbb S}^{N-1})$ defined in \eqref{eq:2rad}, i.e.

\begin{equation} \label{77}
u(\lambda\,\theta)=\sum_{i=1}^\infty\varphi_i(\lambda)\psi_i(\theta)
\quad \text{in }L^2({\mathbb S}^{N-1}),
\end{equation}
with $\varphi_i$ is defined in (\ref{78}).  On the other hand,
$\int_{B_R}(|x_J|^{-2+\e}+|x_{J_1}-x_{J_2}|^{-2+\e} )u^2(x)\,dx<+\infty$
for all $J\in{\mathcal A}_k$ and $(J_1,J_2)\in {\mathcal B}_k$ by
Corollary \ref{c:hardyboundary-0}, hence $\int_{{\mathbb S}^{N-1}}
(|\theta_J|^{-2+\e}+|\theta_{J_1}-\theta_{J_2}|^{-2+\e})
u^2(\lambda\theta)dS(\theta)$ is finite for all $J\in{\mathcal A}_k$,
$(J_1,J_2)\in {\mathcal B}_k$, and a.e. $\lambda\in (0,R)$, which,
together with Lemma \ref{l:hardysphere}, implies that $\theta\mapsto
h(\lambda\theta)u(\lambda\theta)\in H^{-1}({\mathbb S}^{N-1})$ for
a.e. $\lambda\in (0,R)$. Moreover, by \eqref{F}, it is also
easy to verify that $\theta\mapsto
f(\lambda\theta,u(\lambda\theta))\in H^{-1}({\mathbb S}^{N-1})$ for
a.e. $\lambda\in (0,R)$.
Therefore, we may write
\begin{equation}\label{77bis}
  h(\lambda\,\theta)u(\lambda\,\theta)+f(\lambda\,\theta,u(\lambda\,\theta))=
  \sum_{i=1}^\infty\zeta_i(\lambda)\psi_i(\theta) \quad \text{in
  }H^{-1}({\mathbb S}^{N-1})\text{ for a.e. }\lambda\in
  (0,R)\end{equation} where
\begin{align}\label{78bis}
  \zeta_i(\lambda)&= {}_{H^{-1}({\mathbb S}^{N-1})}\big\langle
  h(\lambda\cdot)u(\lambda\cdot)+ f(\lambda\cdot,u(\lambda\cdot)) ,
  \psi_i
  \big\rangle_{H^1({\mathbb S}^{N-1})}\\
  &\notag = \int_{{\mathbb S}^{N-1}} \Big(
  h(\lambda\,\theta)u(\lambda\,\theta)
  +f(\lambda\,\theta,u(\lambda\,\theta))\Big) \psi_i(\theta)
  \,dS(\theta).
\end{align}
We notice that, in view of Remark \ref{rem:l1w11}, $\zeta_i\in
L^{1}_{{\rm loc}}(0,R)$ and
\begin{equation}\label{eq:31}
  \lambda^{N-1}\zeta_i(\lambda)=
\Upsilon_i'(\lambda)\quad\text{a.e. in }(0,R),
\end{equation}
where $\Upsilon_i$ is defined in (\ref{eq:fi}).
Since $u$ solves \eqref{u}, by \eqref{eq:2rad} we obtain that,
 for any $i\in \N$, $\varphi_i$ solves
\begin{equation*}
  -\varphi_i''(\lambda)-\frac{N-1}{\lambda}\varphi_i^\prime(\lambda)+
  \frac{\mu_i(a)}{\lambda^2}\varphi_i(\lambda)=
  \zeta_i(\lambda)\quad\text{in the sense of distributions in }(0,R),
\end{equation*}
which can be also written as
\begin{equation*}
  -\left(
    \lambda^{N-1+2\sigma_i}\big(\lambda^{-\sigma_i}\varphi_i(\lambda)\big)'\right)'
  =\lambda^{N-1+\sigma_i}
  \zeta_i(\lambda)\quad\text{in the sense of distributions in }(0,R),
\end{equation*}
where
\begin{equation*}
  \sigma_i=-\frac{N-2}{2}+\sqrt{\bigg(\frac{N-2}
    {2}\bigg)^{\!\!2}+\mu_i(a)}.
\end{equation*}
Integrating by parts the right hand side and taking into account
(\ref{eq:31}), we obtain that there exists $c_i \in \R$ such that
\begin{equation*}
\big(
    \lambda^{-\sigma_i}\varphi_i(\lambda)\big)'
  =-\lambda^{-N+1-\sigma_i}\Upsilon_i(\lambda)-
\sigma_i\lambda^{-N+1-2\sigma_i}\bigg(
c_i+\int_\lambda^Rs^{\sigma_i-1}\Upsilon_i(s)ds\bigg)
\end{equation*}
in the sense of distributions in $(0,R)$, thus implying that
$\varphi_i\in W^{2,1}_{\rm loc}(0,R)$. A further integration yields
\begin{align}\label{eq:42}
  \varphi_i(\lambda)&=\lambda^{\sigma_i} \bigg(
  R^{-\sigma_i}\varphi_i(R)+
  \int_\lambda^Rs^{-N+1-\sigma_i}\Upsilon_i(s)ds\bigg)\\[5pt]
  \notag&\qquad\qquad+\sigma_i\lambda^{\sigma_i}\int_\lambda^R
  s^{-N+1-2\sigma_i}\bigg(
  c_i+\int_s^Rt^{\sigma_i-1}\Upsilon_i(t)dt\bigg)ds\\[5pt]
  \notag&=\lambda^{\sigma_i} \bigg( R^{-\sigma_i}\varphi_i(R)+
  \int_\lambda^Rs^{-N+1-\sigma_i}\Upsilon_i(s)ds
  +\frac{\sigma_ic_iR^{-N+2-2\sigma_i}}{2-N-2\sigma_i} \bigg)-
  \frac{\sigma_ic_i\lambda^{-N+2-\sigma_i}}{2-N-2\sigma_i}\\[5pt]
  \notag&\qquad\qquad
  +\frac{\sigma_i\lambda^{\sigma_i}}{2-N-2\sigma_i}\int_\lambda^R
  t^{-N+1-\sigma_i}\Upsilon_i(t)\,dt-
  \frac{\sigma_i\lambda^{-N+2-\sigma_i}}{2-N-2\sigma_i}\int_\lambda^R
  t^{\sigma_i-1}\Upsilon_i(t)\,dt\\[5pt]
  \notag&=\lambda^{\sigma_i} \bigg( R^{-\sigma_i}\varphi_i(R)+
  \frac{2-N-\sigma_i}{2-N-2\sigma_i}
  \int_\lambda^Rs^{-N+1-\sigma_i}\Upsilon_i(s)ds
  +\frac{\sigma_ic_iR^{-N+2-2\sigma_i}}{2-N-2\sigma_i} \bigg)\\[5pt]
  \notag&\qquad\qquad
+  \frac{\sigma_i\lambda^{-N+2-\sigma_i}}{N-2+2\sigma_i}
\bigg(c_i+\int_\lambda^R
  t^{\sigma_i-1}\Upsilon_i(t)\,dt\bigg).
\end{align}
Let $j_0,m\in\N$ be as in
(\ref{eq:30}), i.e.
 $m$ is the
multiplicity of the eigenvalue
$$
\mu_{j_0}(a)=\mu_{j_0+1}(a)=\cdots=\mu_{j_0+m-1}(a)
$$
and
\begin{equation}\label{eq:15}
  \gamma=\lim_{r\rightarrow 0^+} {\mathcal N}(r)=\sigma_{i},
  \quad i=j_0,\dots,j_0+m-1,
\end{equation}
see  Lemma \ref{l:blow-up}.
The Parseval identity yields
\begin{equation}\label{eq:17}
H(\lambda)=\int_{{\mathbb
    S}^{N-1}}|u(\lambda\,\theta)|^2\,dS(\theta)=
\sum_{i=1}^{\infty}|\varphi_i(\lambda)|^2,\quad\text{for all }0<\lambda\leq R.
\end{equation}
From Lemma \ref{l:stima}, it follows that
\begin{equation}\label{eq:32}
\Upsilon_i(\lambda)=O(\lambda^{N-2+\delta+\sigma_i})\quad
\text{for every }i\in\{j_0,\dots,j_0+m-1\}\quad
\text{as }\lambda\to0^+.
\end{equation}
From (\ref{eq:32}), it follows that
\begin{equation}\label{eq:34}
s\mapsto s^{-N+1-\sigma_i}\Upsilon_i(s)\in L^1(0,R)\quad
\text{for every }i\in\{j_0,\dots,j_0+m-1\}
\end{equation}
which yields
\begin{multline}\label{eq:33}
\lambda^{\sigma_i} \bigg( R^{-\sigma_i}\varphi_i(R)+
  \frac{2-N-\sigma_i}{2-N-2\sigma_i}
  \int_\lambda^Rs^{-N+1-\sigma_i}\Upsilon_i(s)ds
  +\frac{\sigma_ic_iR^{-N+2-2\sigma_i}}{2-N-2\sigma_i} \bigg)\\
=O(\lambda^{\sigma_i})=o(\lambda^{-N+2-\sigma_i})
\end{multline}
for all $i\in\{j_0,\dots,j_0+m-1\}$
as $\lambda\to0^+$.
From (\ref{eq:32}), it also follows that
$$
t\mapsto t^{\sigma_i-1}\Upsilon_i(t)\in L^1(0,R)\quad
\text{for every }i\in\{j_0,\dots,j_0+m-1\}.
$$
 From
 $\frac{u}{|x|}\in L^2(B_{R})$, we deduce that
$$
\int_0^Rr^{N-3}\varphi_i^2(r)\,dr<+\infty.
$$
 Then, since
$\int_0^Rr^{N-3}(r^{-N+2-\sigma_i})^2\,dr=+\infty$, from
(\ref{eq:42}) and (\ref{eq:33}) it follows that
$$
c_i+\int_0^R
  t^{\sigma_i-1}\Upsilon_i(t)\,dt=0
$$
and hence
\begin{align*}
  \varphi_i(\lambda)&=\lambda^{\sigma_i} \bigg( R^{-\sigma_i}\varphi_i(R)+
  \frac{2-N-\sigma_i}{2-N-2\sigma_i}
  \int_\lambda^Rs^{-N+1-\sigma_i}\Upsilon_i(s)ds
-
\frac{\sigma_iR^{-N+2-2\sigma_i}}{2-N-2\sigma_i}
\int_0^R
  t^{\sigma_i-1}\Upsilon_i(t)\,dt \bigg)\\[5pt]
  \notag&\qquad\qquad
-  \frac{\sigma_i\lambda^{-N+2-\sigma_i}}{N-2+2\sigma_i}
\int_0^\lambda
  t^{\sigma_i-1}\Upsilon_i(t)\,dt.
\end{align*}
On the other hand, from (\ref{eq:32}) it follows that
$$
\lambda^{-N+2-\sigma_i}\int_0^\lambda
  t^{\sigma_i-1}\Upsilon_i(t)\,dt=O(\lambda^{\sigma_i+\delta})
\quad\text{as }\lambda\to 0^+
$$
for all $i\in\{j_0,\dots,j_0+m-1\}$, thus completing the proof.
\end{pf}

\begin{Lemma} \label{l:limitepositivo} Let $\Omega\subset\R^N$, $N\geq
  3$, be a bounded open set such that $0\in\Omega$, $a$ satisfy
  (\ref{eq:ateta}) and (\ref{eq:lambdamin1}), and $u\not\equiv 0$ be a
  weak $H^1(\Omega)$-solution to (\ref{u}), with $h,f$ satisfying
  {\rm (\ref{H1-2}--\ref{F})}.  Then
\[
\lim_{r\to0^+}r^{-2\gamma}H(r)>0.
\]
\end{Lemma}

\begin{pf}
Let $R>0$ be such that $\overline{B_R}\subset \Omega$
and
$j_0,m\in\N$ as in
(\ref{eq:30}).
We argue by  contradiction and assume that
$\lim_{\lambda\to0^+}\lambda^{-2\gamma}H(\lambda)=0$.
Then, letting $\varphi_i$ as in (\ref{78}), (\ref{eq:17}) implies that
\begin{equation}\label{eq:13}
\lim_{\lambda\to0^+}\lambda^{-\gamma}\varphi_{i}(\lambda)=0
\quad\text{for all }i\in\{j_0,\dots,j_0+m-1\}.
\end{equation}
From Lemma \ref{l:phii}, (\ref{eq:34}), and (\ref{eq:13}),
we deduce that
\begin{align*}
 R^{-\gamma}\varphi_i(R) - \frac{\gamma R^{-N+2-2\gamma}}{2-N-2\gamma} \int_0^R
  s^{\gamma-1}\Upsilon_i(s)\,ds=-
  \frac{2-N-\gamma}{2-N-2\gamma} \int_0^Rs^{-N+1-\gamma}\Upsilon_i(s)ds
 \end{align*}
for all $i\in\{j_0,\dots,j_0+m-1\}$. Hence (\ref{eq:35})
can be rewritten as
\begin{align}\label{eq:36}
  \varphi_i(\lambda)=- \frac{2-N-\gamma}{2-N-2\gamma}
  \,\lambda^{\gamma}\int_0^\lambda s^{-N+1-\gamma}\Upsilon_i(s)ds
  +O(\lambda^{\gamma+\delta})\quad\text{as }\lambda\to 0^+
\end{align}
for all $i\in\{j_0,\dots,j_0+m-1\}$. From (\ref{eq:36}), (\ref{eq:15}), and
(\ref{eq:32}), we infer the estimate
\[
\varphi_{i}(\lambda)=O(\lambda^{\gamma+\delta})\quad\text{as
}\lambda\to0^+,\quad\text{for every }i\in\{j_0,\dots,j_0+m-1\},
\]
namely, setting $u^\lambda(\theta)=u(\lambda \theta)$,
\[
(u^\lambda,\psi_i)_{L^2({\mathbb S}^{N-1})}=
O(\lambda^{\gamma+\delta})\quad\text{as
}\lambda\to0^+,\quad\text{for every }i\in\{j_0,\dots,j_0+m-1\},
\]
and hence
\[
(u^\lambda,\psi)_{L^2({\mathbb S}^{N-1})}=
O(\lambda^{\gamma+\delta})\quad\text{as
}\lambda\to0^+,
\]
for every $\psi\in {\mathcal E}_0$,
being ${\mathcal E}_0$  the eigenspace
of the operator $L_{a}$ associated to the eigenvalue
$\mu_{j_0}(a)$.
  Let
$w^\lambda(\theta)=(H(\lambda))^{-1/2}u(\lambda \theta)$.  From (\ref{2ndest}),
there exists $C(\delta)>0$ such that $\sqrt{H(\lambda)}\geq
C(\delta)\lambda^{\gamma+\frac\delta 2}$ for $\lambda$ small, and therefore
\begin{equation}\label{eq:15-bis}
(w^\lambda,\psi)_{L^2({\mathbb S}^{N-1})}=
O(\lambda^{\delta/2})=o(1),\quad\text{as
}\lambda\to0^+
\end{equation}
for every $\psi\in{\mathcal E}_0$.  From Lemma \ref{l:blow-up}, for
every sequence $\lambda_n\to0^+$, there exist a subsequence
$\{\lambda_{n_j}\}_{j\in\N}$ and an eigenfunction $\widetilde
\psi\in{\mathcal E}_0$ such that
\begin{equation}\label{eq:16}
\int_{{\mathbb S}^{N-1}}|\widetilde\psi(\theta)|^2dS=1\quad\text{and} \quad
w^{\lambda_{n_j}}\to \widetilde \psi\quad\text{in } L^2({\mathbb
  S}^{N-1}).
\end{equation}
From (\ref{eq:15-bis}) and (\ref{eq:16}), we infer that
\[
0=\lim_{j\to+\infty}(w^{\lambda_{n_j}},\widetilde\psi)_{L^2({\mathbb S}^{N-1})}
=\|\widetilde\psi\|_{L^2({\mathbb S}^{N-1})}^2=1,
\]
thus reaching a contradiction.
\end{pf}

\noindent Combining Lemma \ref{l:blow-up} with Lemma \ref{l:limitepositivo},
we can now prove Theorem  \ref{Main-h}.

\medskip\noindent
\begin{pfn}{Theorem \ref{Main-h}}
  Identity \eqref{lim-N} follows immediately from Lemma
  \ref{l:blow-up}.  As in the statement of the theorem, let $m$ be the
  multiplicity of the eigenvalue $\mu_{k_0}(a)$ found in Lemma
  \ref{l:blow-up},  $j_0\in \N\setminus\{0\}$,  such that
$j_0\leq k_0\leq j_0+m-1$,
  $\mu_{j_0}(a)=\mu_{j_0+1}(a)=\dots=\mu_{j_0+m-1}(a)$, and 
  $\gamma=\lim_{r\to 0^+} \mathcal N(r)$.

  In order to prove \eqref{convergence}, let $\{\lambda_n\}_{n\in
    \N}\subset (0,\infty)$ be a sequence such that $\lambda_n\to 0^+$
  as $n\to +\infty$.  Then by Lemmas \ref{l:blow-up}, \ref{l:limite},
  and \ref{l:limitepositivo}, there exist a subsequence  $\lambda_{n_j}$ and
  $\beta_{j_0},\dots,\beta_{j_0+m-1}\in \R$
  such that
\begin{equation}\label{eq:23}
  \lambda_{n_j}^{-\gamma}u(\lambda_{n_j}x)\to
  |x|^\gamma\sum_{i=j_0}^{j_0+m-1} \beta_i\psi_{i}\bigg(\frac{x}{|x|}\bigg)
  \quad \text{in }
  H^{1}(B_1) \quad \text{as }j\to+\infty
\end{equation}
and $(\beta_{j_0},\beta_{j_0+1},\dots,\beta_{j_0+m-1})\neq(0,0,\dots,0)$,
which implies
\begin{equation}\label{eq:37}
\lambda_{n_j}^{-\gamma}u(\lambda_{n_j}\theta)\to
\sum_{i=j_0}^{j_0+m-1} \beta_i\psi_{i}(\theta)\quad \text{in }
L^2({\mathbb S}^{N-1}) \quad \text{as }j\to+\infty.
\end{equation}
We now prove
that the $\beta_i$'s depend neither on the sequence
$\{\lambda_n\}_{n\in\N}$ nor on its subsequence
$\{\lambda_{n_j}\}_{j\in\N}$.
Let us fix $R>0$ such that
$\overline{B_{R}}\subset\Omega$.  Defining $\varphi_i$
as in (\ref{78}), from
(\ref{eq:37}) it follows that, for any $i=j_0,\dots, j_0+m-1$,
\begin{equation}\label{eq:25}
\lambda_{n_j}^{-\gamma}\varphi_i(\lambda_{n_j}) =\int_{{\mathbb
    S}^{N-1}}\frac{u(\lambda_{n_j}\theta)}{\lambda_{n_j}^{\gamma}}
\psi_i(\theta)\,dS(\theta)
\to\sum_{\ell=j_0}^{j_0+m-1} \beta_\ell\int_{{\mathbb
    S}^{N-1}}\psi_{\ell}(\theta)\psi_i(\theta)\,dS(\theta)=\beta_i
\end{equation}
as $j\to+\infty$.  On the other hand, from Lemma \ref{l:phii},
it follows that,  for any $i=j_0,\dots, j_0+m-1$,
\begin{align*}
  \lambda^{-\gamma}\varphi_i(\lambda)\ \to\
  &
 R^{-\gamma}\varphi_i(R)+ \frac{2-N-\gamma}{2-N-2\gamma}
  \int_0^Rs^{-N+1-\gamma}\Upsilon_i(s)\,ds- \frac{\gamma
    R^{-N+2-2\gamma}}{2-N-2\gamma} \int_0^R s^{\gamma-1}\Upsilon_i(s) \,ds
\end{align*}
as $\lambda\to0^+$, with $\Upsilon_i$ as in (\ref{eq:fi}),
and therefore from (\ref{eq:25}) we deduce that
\begin{align*}
  \beta_i= R^{-\gamma}\varphi_i(R)+ \frac{2-N-\gamma}{2-N-2\gamma}
  \int_0^Rs^{-N+1-\gamma}\Upsilon_i(s)\,ds- \frac{\gamma
    R^{-N+2-2\gamma}}{2-N-2\gamma} \int_0^R s^{\gamma-1}\Upsilon_i(s) \,ds,
\end{align*}
 for any $i=j_0,\dots, j_0+m-1$.
Integration by parts and (\ref{eq:31}) allow rewriting the above
 formula as in (\ref{eq:38}).
In particular the $\beta_i$'s depend neither on the sequence
$\{\lambda_n\}_{n\in\N}$ nor on its subsequence
$\{\lambda_{n_k}\}_{k\in\N}$, thus implying that the convergence in
(\ref{eq:23}) actually holds as $\lambda\to 0^+$
and proving the theorem.~\end{pfn}

\section{Asymptotic behavior of
  eigenfunctions}\label{sec:asympt-behav-eigenf}

We describe the asymptotic behavior of eigenfunctions of the
operator $L_a=-\Delta_{\SN}-a$ near the singular set of the
function $a$. Actually, for simplicity we study the asymptotic
behavior of eigenfunctions near the south pole as an application
of Theorem \ref{Main-h} after a stereographic projection of $\SN$
onto $\R^{N-1}$ with respect to the ``north pole''.

Throughout this section we assume that $3\leq k\leq N-1$ and that $a$
 satisfies (\ref{eq:ateta}) and (\ref{eq:lambdamin1}).
Note that if $k=N$ then $a$ is constant  and hence the
eigenfunctions of $L_a$ are smooth.

By Lemma \ref{l:spe} the spectrum of $L_a$ consists of a diverging
sequence of eigenvalues $\mu_1(a)< \mu_2(a)\leq \dots \leq\mu_n(a)\leq
\dots$ each of them having finite multiplicity.

Let $\mu_i(a)$ be an eigenvalue of $L_a$ and let $\psi\in H^1(\SN)$
be a corresponding eigenfunction, i.e.
\begin{equation} \label{eig1}
-\Delta_{\SN} \psi(\theta)-a(\theta)\psi(\theta)=\mu_i(a)
\psi(\theta) \qquad {\rm in \ } \SN.
\end{equation}
Let $\Pi:\SN\setminus \{e_N\}\to \R^{N-1}$ be the standard
stereographic projection with respect to the ``north pole''. Here
by $e_N$, we denote the vector $(0,0,\dots,0,1)\in \R^N$. 

Let
$\phi:\R^{N-1}\to \R$ be the function given by
\begin{equation} \label{phi-s}
\phi(y)=\frac{4}{\big(|y|^2+1\big)^2} \qquad {\rm for \ all \ }
y\in \R^{N-1}.
\end{equation}
If 
$\theta\in \SN\setminus\{e_N\}$ and $x,z\in T_\theta \SN$
(by $T_\theta \SN$ we mean the tangent space to $\SN$ at $\theta$), then
$$
(x,z)_{T_\theta \SN}=\phi(\Pi(\theta))\
(d\Pi(\theta)[x],d\Pi(\theta)[z])_{\R^{N-1}}
$$
where the vector space $T_{\Pi(\theta)}\R^{N-1}$ is identified with $\R^{N-1}$. 
In the following lemma the equation satisfied by the projection of $\psi$ 
is deduced.

\begin{Lemma}\label{eq-conforme} 
Let $3\leq k\leq N-1$,  $a$ satisfy (\ref{eq:ateta}) and (\ref{eq:lambdamin1}),
 and let
  $\Pi$ and $\phi$ be respectively the stereographic projection with
  respect to the north pole and the function defined in (\ref{phi-s}).
  Let $\mu_i(a)$ be an eigenvalue of the operator $L_a$ and let
  $\psi\in H^1(\SN)$ be a corresponding eigenfunction. Then the
  function 
\begin{equation} \label{eq:psi-tilde}
\widetilde \psi(y):=(\phi(y))^{\frac{N-3}4} \psi(\Pi^{-1}(y))
\end{equation}
belongs to ${\mathcal D}^{1,2}(\R^{N-1})$ and weakly solves 
\begin{equation} \label{eq-proiettata}
-\Delta \widetilde \psi(y)-\frac{b(\frac{y}{|y|})}{|y|^2} \, \widetilde \psi(y)=\widetilde h(y)
\widetilde\psi(y) 
\end{equation}
where $b$ and $\widetilde h$ are defined by 
\begin{equation} \label{ref:b}
b(\omega)= \sum_{J\in \mathcal A_k, N\notin J}
\frac{\alpha_J}{|\omega_J|^2}+\sum_{(J_1,J_2)\in \mathcal B_k, N\notin
J_1\cup J_2} \frac{\alpha_{J_1 J_2}}{|\omega_{J_1}-\omega_{J_2}|^2},
\quad\text{for any }\omega\in  {\mathbb
  S}^{N-2}\setminus \Sigma_1,
\end{equation}
where 
\begin{align*}
\Sigma_1:=&\{ (\omega_1,\dots,\omega_{N-1})\in {\mathbb
  S}^{N-2}: \omega_J=0 \text{ for some }J\in{\mathcal A}_k^{N-1}\} \\
&\qquad\cup \{(\omega_1,\dots,\omega_{N-1})\in {\mathbb
  S}^{N-2}:
\omega_{J_1}=\omega_{J_2} \text{ for some } (J_1,J_2)\in \mathcal
B_k^{N-1}\},
\end{align*}
\begin{gather*}
{\mathcal A}_k^{N-1}\!\!=\big\{J\subseteq\{1,2,\dots,N-1\}:\#J=k\big\},\\
{\mathcal B}_k^{N-1}\!\!=\{(J_1,J_2)\in {\mathcal A}_k^{N-1}\times
{\mathcal A}_k^{N-1}\!\!: J_1\cap J_2=\emptyset, J_1<J_2 \},
\end{gather*}
and 
\begin{multline} \label{ref:htilde}
\widetilde h(y)=\phi(y)\left(\frac{(N-3)(N-1)}{4}+\mu_i(a)\right)
+\sum_{J\in \mathcal A_k, N\in J}
\frac{4\alpha_J}{4|y_{J'}|^2+(|y|^2-1)^2} \\
+\sum_{(J_1,J_2)\in \mathcal B_k, N\in J_1\setminus J_2}
\frac{4\alpha_{J_1
J_2}}{4|y_{J'_1}-y_{J'_2}|^2+(|y|^2-1-2y_{m_k})^2} \\
+\sum_{(J_1,J_2)\in \mathcal B_k, N\in J_2\setminus J_1}
\frac{4\alpha_{J_1
J_2}}{4|y_{J'_1}-y_{J'_2}|^2+(|y|^2-1-2y_{n_k})^2} \, .
\end{multline}
for a.e. $y\in \R^{N-1}$, where for any $(J_1,J_2)\in \mathcal B_k$,
$n_k=\max J_1$, $m_k=\max J_2$, and for any
 $J=\{n_1,\dots,n_k\}\in\mathcal A_k$, $n_1<n_2<\dots<n_k$, we denote
$J'=J\setminus\{n_k\}\in \mathcal A_{k-1}$.
\end{Lemma}

\begin{pf} 
The  conformal laplacian  on $\SN$ is
given by
$$
-\Delta_{\SN}+\frac{(N-3)(N-1)}{4},
$$
while, since $\R^{N-1}$ has zero scalar curvature, the conformal
laplacian  in $\R^{N-1}$ coincides with the usual Laplace operator.
Then for any function $\eta\in C^2(\SN\setminus \{e_N\})$ we have
\begin{equation} \label{conformal}
-\Delta_{\SN} \eta(\theta)
+\frac{(N-3)(N-1)}{4}\eta(\theta)=-\phi^{-\frac{N+1}{4}}
\Delta(\phi^{\frac{N-3}{4}}\cdot(\eta\circ \Pi^{-1}))\bigg|_
{\Pi(\theta)} 
\end{equation}
for every $\theta\in \SN\setminus \{e_N\}$.  For the definition of the
conformal Laplacian and for a proof of \eqref{conformal} see \cite[\S
3]{chang} or \cite[(1.2.27)]{Baum-Juhl}.

We claim that the function 
$\widetilde\psi$ defined in (\ref{eq:psi-tilde}) belongs to ${\mathcal D}^{1,2}(\R^{N-1})$  and weakly solves 
\begin{equation} \label{eq-proj}
-\Delta \widetilde \psi(y)-\phi(y)a(\Pi^{-1}(y)) \widetilde
\psi(y)=\left(\frac{(N-3)(N-1)}{4}+\mu_i(a)\right)\phi(y)\widetilde\psi(y)
\quad \text{in } \R^{N-1},
\end{equation}
i.e.
\begin{multline*}
  \int_{\R^{N-1}} \nabla\widetilde\psi(y)\cdot \nabla v(y)\, dy
  -\int_{\R^{N-1}} \phi(y)a(\Pi^{-1}(y))\widetilde \psi(y)v(y)\, dy\\
  =\bigg(\frac{(N-3)(N-1)}{4}+\mu_i(a)\bigg)\int_{\R^{N-1}}
  \phi(y)\widetilde\psi(y)v(y)\,
  dy \quad \text{for all } v\in \mathcal D^{1,2}(\R^{N-1}).
\end{multline*}
Indeed by \eqref{conformal}, integration by parts, and the change of
variables $y:=\Pi(\theta)\in \R^{N-1}$, for any $v_1,v_2\in
C^{\infty}_c(\R^{N-1})$ we have
\begin{equation} \label{S-R-nabla}
\int_{\R^{N-1}} \nabla v_1(y)\cdot\nabla v_2(y)\, dy=
\int_{\SN} \!\! \Big(
\nabla_{\SN} w_1(\theta)\cdot
\nabla_{\SN} w_2(\theta)+{\textstyle{\frac{(N-3)(N-1)}{4}}}
w_1(\theta)w_2(\theta)\Big)\, dS(\theta)
\end{equation}
with $w_j(\theta)=\phi(\Pi(\theta))^{-\frac{N-3}4} v_j(\Pi(\theta))$, $j=1,2$.   
Moreover
\begin{equation} \label{S-R-hardy-0}
\int_{\R^{N-1}} v_1(y)v_2(y)\phi(y)\, dy
=
\int_{\SN} w_1(\theta)w_2(\theta)\, dS(\theta) 
\end{equation}
and
\begin{equation} \label{S-R-hardy}
\int_{\R^{N-1}} 
a(\Pi^{-1}(y))v_1(y)v_2(y)\phi(y)\, dy
=
\int_{\SN} a(\theta) w_1(\theta)w_2(\theta)\, dS(\theta) 
\end{equation}
with $w_1,w_2$ as above. 
By  density,  (\ref{S-R-nabla}--\ref{S-R-hardy})
 actually hold for any 
$v_1,v_2\in \mathcal D^{1,2}(\R^{N-1})$ and hence the claim  follows.

We now write the function $a(\Pi^{-1}(y))$ in a more explicit way.
We recall that the function
$\Pi^{-1}:\R^{N-1}\to \SN$ is given by
\begin{equation*}
\Pi^{-1}(y)=\frac{2}{|y|^2+1}\, y+\frac{|y|^2-1}{|y|^2+1} \, e_N
\end{equation*}
where we identified $\R^{N-1}$ with the subspace of $\R^N$ of all
$x=(x_1,\dots,x_N)\in \R^N$ such that $x_N=0$. Therefore for any
$\theta\in \SN\setminus\{e_N\}$, if
$y=(y_1,\dots,y_{N-1})=\Pi(\theta)$, we have
\begin{equation*}
\theta_J=
\begin{cases}
\frac{2}{|y|^2+1} \, y_J=\left(\frac{2}{|y|^2+1} \, y_{J'},
\frac{2}{|y|^2+1}\, y_{n_k}\right),
&\text{if } N\notin J,\\[6pt]
\left(\frac{2}{|y|^2+1} \, y_{J'}, \frac{|y|^2-1}{|y|^2+1}\right),
&\text{if } N\in J \, .
\end{cases}
\end{equation*}
Hence, for every $J\in\mathcal A_k$,
\begin{equation*}
|\theta_J|^2=
\begin{cases}
\frac{4}{(|y|^2+1)^2} \, |y_J|^2, &\text{if } N\notin J,\\[7pt]
\frac{4|y_{J'}|^2+(|y|^2-1)^2}{(|y|^2+1)^2}, &\text{if } N\in J,
\end{cases}
\end{equation*}
and, for every $(J_1,J_2)\in \mathcal B_k$, with $J_1=\{n_1,\dots,n_k\}$
and $J_2=\{m_1,\dots,m_k\}$, 
\begin{equation} \label{cases}
|\theta_{J_1}-\theta_{J_2}|^2=
\begin{cases}
\frac{4}{(|y|^2+1)^2} \, |y_{J_1}-y_{J_2}|^2, &\text{if } N\notin J_1\cup
J_2,\\[8pt]
\frac{4|y_{J_1'}-y_{J_2'}|^2+(|y|^2-1-2y_{m_k})^2}{(|y|^2+1)^2},
&\text{if }
N\in J_1\setminus J_2, \\[8pt]
\frac{4|y_{J_1'}-y_{J_2'}|^2+(|y|^2-1-2y_{n_k})^2}{(|y|^2+1)^2},
&\text{if } N\in J_2\setminus J_1 \, .
\end{cases}
\end{equation}
By \eqref{cases} we obtain
\begin{multline} \label{phi-a} 
\phi(y)a(\Pi^{-1}(y))=\sum_{J\in \mathcal A_k, N\notin J}
\frac{\alpha_J}{|y_J|^2}+\sum_{(J_1,J_2)\in \mathcal B_k, N\notin J_1\cup J_2}
\frac{\alpha_{J_1 J_2}}{|y_{J_1}-y_{J_2}|^2}\\
+\sum_{J\in \mathcal A_k, N\in J}
\frac{4\alpha_J}{4|y_{J'}|^2+(|y|^2-1)^2}+\sum_{(J_1,J_2)\in
\mathcal B_k, N\in J_1\setminus J_2} \frac{4\alpha_{J_1
J_2}}{4|y_{J'_1}-y_{J'_2}|^2+(|y|^2-1-2y_{m_k})^2} \\
+\sum_{(J_1,J_2)\in \mathcal B_k, N\in J_2\setminus J_1}
\frac{4\alpha_{J_1
J_2}}{4|y_{J'_1}-y_{J'_2}|^2+(|y|^2-1-2y_{n_k})^2} \, .
\end{multline}
The conclusion follows from (\ref{cases}) and (\ref{phi-a}). 
\end{pf}

\noindent According with \eqref{eq:bound} we introduce the number 
\begin{equation}
\label{Lambda(b)}
\Lambda(b):=\sup_{v\in \mathcal D^{1,2}(\R^{N-1})\setminus\{0\}}\dfrac{{\displaystyle
{\int_{\R^{N-1}}{{|y|^{-2}}{b(y/|y|)}\,v^2(y)\,dy}}}}
{{\displaystyle{\int_{\R^{N-1}}{|\nabla v(y)|^2\,dy}}}} \, .
\end{equation}
\begin{Lemma} \label{Lambda(b)<1}
Let $3\leq k\leq N-1$,  $a$  satisfy (\ref{eq:ateta}) and (\ref{eq:lambdamin1}),
 and let $b$ be 
the corresponding function defined in \eqref{ref:b}. 
Then $\Lambda(b)<1$ with $\Lambda(b)$ as 
in (\ref{Lambda(b)}).  
\end{Lemma}

\begin{pf} Let $v\in C^\infty_c(\R^{N-1})$ and let
$w(\theta)=\phi(\Pi(\theta))^{-\frac{N-3}4} v(\Pi(\theta))$. 
Then by (\ref{S-R-nabla}-\ref{S-R-hardy}) and \eqref{supremum} we have
\begin{align*} 
  \int_{\R^{N-1}} & |\nabla v(y)|^2 dy-\frac{(N-3)(N-1)}{4}
  \int_{\R^{N-1}}
  \phi(y) |v(y)|^2 dy-\int_{\R^{N-1}} \phi(y)a(\Pi^{-1}(y))|v(y)|^2 dy \\
  &=\int_{\SN} |\nabla_{\SN} w(\theta)|^2 dS(\theta)-
  \int_{\SN} a(\theta)|w(\theta)|^2 dS(\theta) \\
  & \geq (1-\Lambda(a))\int_{\SN} |\nabla_{\SN} w(\theta)|^2
  dS(\theta)
  -\Lambda(a)\left(\frac{N-2}2\right)^2 \int_{\SN} |w(\theta)|^2dS(\theta) \\
  & =(1-\Lambda(a)) \left(\int_{\R^{N-1}} |\nabla v(y)|^2 dy
    -\frac{(N-3)(N-1)}{4} \int_{\R^{N-1}} \phi(y)|v(y)|^2 dy \right) \\
  & \quad -\Lambda(a)\left(\frac{N-2}2\right)^2 \int_{\R^{N-1}}
  \phi(y)|v(y)|^2 dy
\end{align*}
and, in turn,
\begin{align} \label{pos-b}
\Lambda(a)&\int_{\R^{N-1}}  |\nabla v(y)|^2 dy
+\Lambda(a)\left[\left(\frac{N-2}2\right)^2-\frac{(N-3)(N-1)}{4}\right]
\int_{\R^{N-1}} \phi(y)|v(y)|^2 dy \\
\notag
& \geq
\int_{\R^{N-1}} \phi(y)a(\Pi^{-1}(y))|v(y)|^2 dy 
=\int_{\R^{N-1}} \frac{b(y/|y|)}{|y|^2}|v(y)|^2 dy +
\int_{\R^{N-1}} R(y)|v(y)|^2 dy
\end{align}
where $R(y)=\phi(y)a(\Pi^{-1}(y))-\frac{b(y/|y|)}{|y|^2}$ is bounded
in a sufficiently small neighborhood of $0$ by (\ref{phi-a}).
On the other hand if we define, for any $\delta>0$,
$B_\delta\subset \R^{N-1}$ to be  the
  open ball of radius $\delta$ centered at the origin and 
\begin{equation} \label{sup-delta} 
C(\delta):=\sup_{
\substack{
v\in C^\infty(\R^{N-1})\setminus\{0\}\\
 {\rm supp}\, v\subset B_\delta}}
\dfrac{\int_{\R^{N-1}} 
\Big[\Big(\Lambda(a)\left(\frac{N-2}2\right)^2-\frac{(N-3)(N-1)}{4}\Big)
\phi(y)-R(y)
\Big]|v(y)|^2 dy}{\int_{\R^{N-1}} |\nabla v(y)|^2 dy}  ,
\end{equation}
then $C(\delta)\to 0$ as $\delta\to 0^+$.
Therefore, by \eqref{pos-b}, (\ref{sup-delta}), and (\ref{eq:lambdamin1}),
 we deduce that there exists $\delta_1>0$ such that for any  
$\delta\in (0,\delta_1)$ 
\begin{align*}
& \sup_{v\in C^\infty_c(\R^{N-1})\setminus\{0\}} 
\frac{\int_{\R^{N-1}} \frac{b(\frac{y}{|y|})}{|y|^2}\, |v(y)|^2 dy}
{\int_{\R^{N-1}} |\nabla v(y)|^2 dy}=
\sup_{\substack{v\in C^\infty(\R^{N-1})\setminus\{0\}\\
 {\rm supp}\, v\subset B_\delta}}
\frac{\int_{\R^{N-1}} \frac{b(\frac{y}{|y|})}{|y|^2}\, |v(y)|^2 dy}
{\int_{\R^{N-1}} |\nabla v(y)|^2 dy}
\leq \Lambda(a)+C(\delta)
<1.
\end{align*}
The conclusion follows by density.
\end{pf}

By Lemmas \ref{Lambda(b)<1} and \ref{l:spe}, we
deduce that the spectrum of the operator $L_b:=-\Delta_{\mathbb
  S^{N-2}}-b$ on $\mathbb S^{N-2}$ consists of real eigenvalues with
finite multiplicity $\mu_1(b)< \mu_2(b)\leq \dots \leq \mu_k(b)\leq
\dots$.

Let $\widetilde h$ be the function defined in \eqref{ref:htilde} and, 
according with \eqref{eq:14},
for any nontrivial $\mathcal D^{1,2}(\R^{N-1})$-solution 
$v$ of the equation
\begin{equation*}
-\Delta v(y)-\frac{b(\frac{y}{|y|})}{|y|^2}\, v(y)=\widetilde h(y) v(y), 
\end{equation*}
we define the corresponding Almgren's frequency function by
\begin{equation}
\mathcal{N}_{v,\widetilde h,0}(r) 
=\frac{r\int_{B_r} \big(|\nabla v(y)|^2
-\frac{b(y/|y|)}{|y|^2}v^2(y) -
    \widetilde h(y)v^2(y) \big) \, dy}{\int_{\partial B_r}|v(y)|^2 \, dS} \, .
\end{equation}
 
We are ready to prove the following asymptotic description of 
eigenfunctions.

\begin{Proposition} \label{autofunzioni} Let $3\leq k\leq N-1$, let
  $a$  satisfy (\ref{eq:ateta}), (\ref{eq:lambdamin1}),
  and let $b$ and $\widetilde h$ be 
  respectively defined in (\ref{ref:b}) and (\ref{ref:htilde}).  Let
  $\mu_i(a)$ be an eigenvalue of the operator $L_a$ and let $\psi\in
  H^1(\SN)\setminus\{0\}$ be an associated eigenfunction.  Let
  $\widetilde \psi\in \mathcal D^{1,2}(\R^{N-1})$ be the corresponding
  function defined in (\ref{eq:psi-tilde}).  Then there exists
  $\widetilde k_0\in \N$, $\widetilde k_0\geq 1$, such that
\begin{equation}\label{lim-N-b}
  \lim_{r\to 0^+} \mathcal
  N_{\widetilde \psi,\widetilde h,0}(r)=-\frac{N-3}2+
\sqrt{\left(\frac{N-3}2\right)^{\!\!2}
    +\mu_{\widetilde k_0}(b)}.
\end{equation}
Furthermore, if $\widetilde \gamma$ denotes the limit in
(\ref{lim-N-b}), $\widetilde m\geq 1$ is the multiplicity of the
eigenvalue $\mu_{\widetilde k_0}(b)$ and $\{\eta_j:\widetilde \ell_0\leq
j\leq \widetilde \ell_0+\widetilde m-1\}$ ($\widetilde \ell_0\leq \widetilde
k_0\leq \widetilde \ell_0+\widetilde m-1$) is an $L^2(\mathbb
S^{N-2})$-orthonormal basis for the eigenspace associated to
$\mu_{\widetilde k_0}(b)$, then
\begin{equation*}
  \lambda^{-\widetilde\gamma} \psi(\Pi^{-1}(\lambda\Pi(\theta)))\to
  4^{-\frac{N-3}{4}}|\Pi(\theta)|^
  {\widetilde\gamma}\sum_{j=\widetilde \ell_0}^{\widetilde \ell_0+\widetilde m-1}
  \widetilde\beta_j \eta_j\bigg(\frac{\Pi(\theta)}{|\Pi(\theta)|}\bigg)
 \quad \text{in } H^1_{\rm{loc}}(\SN\setminus\{e_N\}) \quad
  \text{as } \lambda\to 0^+,
\end{equation*}
where
\begin{align*}
  \widetilde\beta_j= \int_{{\mathbb S}^{N-2}}\!\bigg[
  R^{-\widetilde\gamma}\widetilde\psi(R\omega)+
  \int_{0}^R\frac{\widetilde
    h(s\omega)\widetilde\psi(s\omega)}{2\widetilde\gamma+N-3}
  \bigg(s^{1-\widetilde\gamma}-\frac{s^{\widetilde\gamma+N-2}}
{R^{2\widetilde\gamma+N-3}}\bigg)ds
  \bigg]\eta_{j}(\omega)\,dS(\omega),
 \end{align*}
 for all $R\in(0,\overline{R})$ for some $\overline{R}>0$, and 
$(\widetilde\beta_{\tilde \ell_0},\widetilde
 \beta_{\widetilde \ell_0+1},\dots, \widetilde\beta_{\widetilde
   \ell_0+\widetilde m-1})\neq(0,0,\dots,0)$.
\end{Proposition}

\begin{pf} Since $\psi$ is a solution of \eqref{eig1}, then, by Lemma
  \ref{eq-conforme}, $\widetilde \psi$ solves
  \eqref{eq-proiettata}. By Lemma \ref{Lambda(b)<1}, $\Lambda(b)<1$
  i.e.  the function $b$ satisfies the positivity condition required
  in Theorem \ref{Main-h}.  Moreover by \eqref{ref:htilde}, the
  function $\widetilde h\in C^{1}(\overline B_\delta)$ for some
  $\delta>0$ small enough.  Hence we may
  apply Theorem~\ref{Main-h} to the function $\widetilde \psi$ to conclude.
\end{pf}

\section{Pointwise  estimates}\label{sec:pointw-estim-under}

Let $\hat \sigma$ as in (\ref{def-sigma-hat}) and 
$\hat \psi_1\in H^1(\SN)$, $\|\hat
\psi_1\|_{L^2(\SN)}=1$,  be the first positive eigenfuntion of the
eigenvalue problem $L_a \psi=\mu_1(\hat a)\psi$ in $\SN$.
The following lemma holds true.
\begin{Lemma}\label{l:all_pos}
If $\hat a$ satisfies (\ref{eq:atetahat}) and  (\ref{eq:53}) ,  then
$$
\mu_1(\hat a)\leq 0,\quad \hat\sigma\leq 0,\quad\text{and} \quad
\inf_{\SN} \hat \psi_1>0.
$$
\end{Lemma}
\begin{pf}
  The fact that $\mu_1(\hat a)\leq 0$ follows easily by taking a
  constant function in the Rayleigh quotient minimized by
  $\mu_1(\hat a)$ (see (\ref{firsteig})).  Moreover, there exists $\delta>0$ 
such that, letting
$$
\Sigma_\delta:=\bigg(\bigcup_{\substack{J\in{\mathcal A}_k\\
\alpha_J>0}}\!\! 
\{ (\theta_1,\dots,\theta_{N})\in {\mathbb
  S}^{N-1}: |\theta_J|<\delta\}\bigg) \cup 
\bigg(\bigcup_{\substack{(J_1,J_2)\in \mathcal
B_k\\
\alpha_{J_1J_2}>0}}\!\!\!\!\{
(\theta_1,\dots,\theta_N)\in \SN:
|\theta_{J_1}-\theta_{J_2}|<\delta\}\bigg)
$$
there holds $\hat a(\theta)+\mu_1(\hat a)>0$ in $\Sigma_\delta$. 
By
classical elliptic regularity theory and maximum principles applied to
the equation satisfied by $\hat\psi_1$ in $\SN\setminus \Sigma_{\delta/2}$, we 
have that $\min_{\SN\setminus\Sigma_{\delta/2}}\hat\psi_1>0$ and 
 $\min_{\partial\Sigma_{\delta}}\hat\psi_1>0$.
Moreover, testing 
$$
\begin{cases}
-\Delta_{\SN}(\hat\psi_1-\min_{\partial\Sigma_{\delta}}\hat\psi_1)
=(\mu_1(\hat a)+\hat a(\theta))\hat\psi_1\geq 0,\quad 
\text{in }\Sigma_{\delta},\\
\hat\psi_1-\min_{\partial\Sigma_{\delta}}\hat\psi_1\geq 0\quad 
\text{on }\partial\Sigma_{\delta},
\end{cases}
$$
with $-(\hat\psi_1-\min_{\partial\Sigma_{\delta}}\hat\psi_1)^-$ we obtain that 
$\hat\psi_1\geq \min_{\partial\Sigma_{\delta}}\hat\psi_1$ in $\Sigma_{\delta}$.
\end{pf}

\noindent Let us introduce the weight function
\begin{equation} \label{PHI}
\rho(x)=|x|^{\hat\sigma} \hat\psi_1\bigg(\frac{x}{|x|}\bigg) \qquad {\rm
for \ all \ } x\in \R^N\setminus \widetilde \Sigma.
\end{equation}
From Lemma \ref{l:all_pos}, under assumptions
(\ref{eq:atetahat}) and  
(\ref{eq:53}), there holds
\begin{equation}\label{eq:d}
  d=d(\mathop{\rm diam}\Omega,N,\hat a):=
\sup_{\Omega\setminus \widetilde \Sigma}\rho^{2-2^*}\in (0,+\infty).
\end{equation}

We notice that $\rho\in H^1_{\rm loc}(\R^N)$ and introduce the
weighted space ${\mathcal D}^{1,2}_{\rho}(\R^N)$ as the completion of
$C^\infty_{\rm c}(\R^N)$ with respect to the norm
\begin{equation} \label{normD12}
  \|v\|_{{\mathcal D}^{1,2}_{\rho}(\R^N)}:
  =\bigg(\int_{\R^N} \rho^2(x)\big|\nabla v(x)\big|^2dx\bigg)^{\!\!1/2}
\end{equation}
and, similarly,  $\mathcal D^{1,2}_\rho
(\Omega)$ as the completion of $C^\infty_c(\Omega)$ with respect
to  \eqref{normD12}.

\begin{Lemma}\label{l:d12density}
 $C^\infty_c(\R^N\setminus \widetilde \Sigma)$ is dense in $\mathcal
D^{1,2}(\R^N)$.
\end{Lemma}
\begin{pf}
  By density of $C^\infty_c(\R^N)$ in $\Di$, it is enough to prove,
  for every $J\in{\mathcal A}_k$ and $(J_1,J_2)\in{\mathcal B}_k$, the
  density of $C^\infty_c(\R^N\setminus \{x_J=0\})$ and of
  $C^\infty_c(\R^N\setminus \{x_{J_1}=x_{J_2}\})$ in
  $C^\infty_c(\R^N)$ with respect to the norm $\|\cdot\|_{\Di}$. Let
  $\phi\in C^\infty(0,\infty)$ such that $\phi(t)=0$ for all $t\in
  (0,1)$ and $\phi(t)=1$ for all $t\in (2,\infty)$. If $J\in{\mathcal
    A}_k$ and $u\in C^\infty_c(\R^N)$, let $u_n(x)=\phi(n|x_J|)u(x)\in 
C^\infty_c(\R^N\setminus \{x_J=0\})$. Since
$$
\nabla u_n(x)-\nabla u(x)=
\nabla u(x)\big(\phi(n|x_J|)-1\big)+nu(x)\phi'(n|x_J|)\frac{x_J}{|x_J|},
$$
$$
\lim_{n\to+\infty}\int_{\R^N}|\nabla
u(x)|^2\big(\phi(n|x_J|)-1\big)^2\,dx =0
$$ 
by the Dominated Convergence Theorem, and 
$$
n^2\int_{\R^N}u^2(x)(\phi'(n|x_J|))^2dx
=n^{2-k}\int_{\R^N}u^2\Big(y_1,\dots,\frac{y_J}{n},\dots,y_N\Big)
(\phi'(|y_J|))^2dy
=O(n^{2-k})
$$
as $n\to+\infty$, we conclude that $u_n\to u$ in $\Di$, thus proving the density 
of  $C^\infty_c(\R^N\setminus \{x_J=0\})$ in
  $C^\infty_c(\R^N)$ and hence in $\Di$. The density 
of   $C^\infty_c(\R^N\setminus \{x_{J_1}=x_{J_2}\})$ can be proven in a similar way.
\end{pf}

\begin{Lemma}\label{l:d12rho}
If $\hat a$ satisfy (\ref{eq:atetahat}) and (\ref{eq:53}), then 
\begin{itemize}
\item[(i)] $C^\infty_c(\R^N\setminus \widetilde \Sigma)$ is dense in $\mathcal
D^{1,2}_\rho(\R^N)$;
\item[(ii)] $v\in \mathcal
D^{1,2}_\rho(\R^N)$ if and only if $\rho v\in \Di$;
\item[(iii)] for all $v\in \mathcal
D^{1,2}_\rho(\R^N)$
\begin{equation}\label{eq:48}
  \int_{\R^N} \rho^2(x) |\nabla v(x)|^2 dx=
  \int_{\R^N} \bigg(  |\nabla (\rho v)(x)|^2 dx  - 
  \frac{\hat a(\frac{x}{|x|})}{|x|^2}
  (\rho v)^2(x)\bigg) \, dx 
\end{equation}
\end{itemize}
\end{Lemma}
\begin{pf}
  We first prove that (\ref{eq:48}) holds for all $v\in
  C^\infty_c(\R^N\setminus \widetilde \Sigma)$. Indeed, 
by direct computation  $\rho$ solves
\begin{equation} \label{eq-rho}
-\Delta \rho(x)-\frac{\hat a(\frac{x}{|x|})}{|x|^2} \,\rho(x)=0 \qquad
\text{in } \R^N\setminus \widetilde \Sigma.
\end{equation}
Let $v\in C^\infty_c(\R^N\setminus \widetilde\Sigma)$ and put 
$u=\rho v$ so that $u\in
C^\infty_c(\R^N)\subset \mathcal D^{1,2}(\R^N)$. Then, testing 
\eqref{eq-rho} with $\rho v^2$ we obtain
\begin{equation} \label{rho-1} \int_{\R^N} \nabla
  \rho(x)\nabla(\rho(x)v^2(x)) \, dx-\int_{\R^N}
  \frac{\hat a(\frac{x}{|x|})}{|x|^2} \, \rho^2(x)v^2(x) \, dx=0.
\end{equation}
Moreover
\begin{equation} \label{rho-2} \nabla \rho\nabla(\rho v^2)=v^2 |\nabla
  \rho|^2+2\rho v \nabla \rho\nabla v
\end{equation}
and
\begin{equation} \label{rho-3} |\nabla u|^2=v^2|\nabla
  \rho|^2+2v\rho\nabla\rho\nabla v+\rho^2|\nabla v|^2.
\end{equation}
By \eqref{rho-1}-\eqref{rho-3} we then have
\begin{align}\label{eq:49}
  Q_{\hat a}(\rho v)&=\int_{\R^N} |\nabla u(x)|^2 dx -\int_{\R^N}
  \frac{\hat a(\frac{x}{|x|})}{|x|^2}
  u^2(x) \, dx \\
  &\notag=\int_{\R^N} \rho^2(x) |\nabla v(x)|^2 dx+\int_{\R^N} \nabla
  \rho(x)\nabla(\rho(x)v^2(x))\, dx-\int_{\R^N}
  \frac{\hat a(\frac{x}{|x|})}{|x|^2}
  \rho^2(x)v^2(x)\, dx\\
  &\notag= \int_{\R^N} \rho^2(x) |\nabla v(x)|^2 dx, \quad \text{for
    all }v\in C^\infty_c(\R^N\setminus \widetilde \Sigma).
\end{align}
To prove (i), by density of $C^\infty_c(\R^N)$ in ${\mathcal
  D}^{1,2}_\rho(\R^N)$, it is enough to prove the density of
$C^\infty_c(\R^N\setminus \widetilde \Sigma)$ in $C^\infty_c(\R^N)$
with respect to the norm $\|\cdot\|_{{\mathcal
    D}^{1,2}_\rho(\R^N)}$. Let $v\in C^\infty_c(\R^N)$. It is easy to
verify that $u=\rho v\in \Di$, hence, by Lemma \ref{l:d12density},
there exists a sequence $\{u_n\}_n\subset C^\infty_c(\R^N\setminus
\widetilde \Sigma)$ such that $u_n\to u$ in $\Di$.  Letting
$v_n=\frac{u_n}{\rho}$, we have that $v_n\in C^\infty_c(\R^N\setminus
\widetilde \Sigma)$ and, by (\ref{eq:49}),
$$
\int_{\R^N} \rho^2(x) |\nabla v_n(x)-\nabla v_m(x)|^2 dx=
Q_{\hat a}(u_n-u_m).
$$
Therefore, since $u_n$ is a Cauchy sequence in $\Di$ and, by
(\ref{eq:53}) and Lemma \ref{l:pos}, $(Q_{\hat a}(u))^{1/2}$ is an
equivalent norm in $\Di$, we conclude that $v_n$ is a
Cauchy sequence in $ {\mathcal D}^{1,2}_\rho(\R^N)$ and hence
converges to some $\tilde v\in {\mathcal D}^{1,2}_\rho(\R^N)$. Since
$v_n\to v$ a.e. in $\R^N$, we deduce that $\tilde v=v$ and then
$v_n\to v$ in ${\mathcal D}^{1,2}_\rho(\R^N)$. The proof of (i) is
thereby complete.

To prove (ii-iii), let $v\in {\mathcal D}^{1,2}_\rho(\R^N)$. By (i), there
exists a sequence $\{v_n\}_n\subset C^\infty_c(\R^N\setminus
\widetilde \Sigma)$ such that $v_n\to v$ in ${\mathcal
  D}^{1,2}_\rho(\R^N)$.  Letting $u_n=v_n\rho\in C^\infty_c(\R^N\setminus
\widetilde \Sigma)$, by (\ref{eq:49}) we have that 
\begin{equation}\label{eq:50}
\int_{\R^N} \rho^2(x) |\nabla v_n(x)|^2 dx
=
\int_{\R^N} |\nabla u_n(x)|^2 dx -\int_{\R^N}
  \frac{\hat a(\frac{x}{|x|})}{|x|^2}
  u_n^2(x) \, dx
\end{equation}
and $\|v_n-v_m\|^2_{{\mathcal D}^{1,2}_\rho(\R^N)}= Q_{\hat a}(u_n-u_m)$.
Therefore, since $v_n$ is a Cauchy sequence in ${\mathcal
  D}^{1,2}_\rho(\R^N)$ and $(Q_{\hat a}(u))^{1/2}$ is an equivalent norm in
$\Di$, we infer that $u_n$ is a Cauchy sequence in $\Di$ and hence
converges to some $u$ in $\Di$. Since $u_n=\rho v_n\to \rho v$ a.e. in
$\R^N$, we deduce that $\rho v=u\in\Di$. Moreover, we can pass to the
limit in (\ref{eq:50}), thus obtaining (\ref{eq:48}) and proving
(iii).  In a similar way, one can prove that if $u\in\Di$ then $\frac
u\rho\in {\mathcal D}^{1,2}_\rho(\R^N)$, thus completing the proof of
(ii).
\end{pf}

\noindent Thanks to Lemma \ref{l:pos}, (\ref{eq:53}), and the
standard Sobolev inequality, the number
\begin{equation*}
S(\hat a)=\inf_{u\in\mathcal D^{1,2}(\R^N)\setminus \{0\}}
\frac{Q_{\hat a}(u)}{\left(\int_{\R^N} |u(x)|^{2^*} \, dx\right)^{2/2^*}}
 \end{equation*}
 is strictly positive and provides the best constant in the following
 weighted Sobolev inequality.

 \begin{Lemma} \label{l:W-Sob} Let $N\geq k\geq 3$ and let $\hat a$ satisfy
  (\ref{eq:atetahat}) and (\ref{eq:53}). Then
\begin{equation} \label{eq:W-Sob}
\int_{\R^N} \rho^2(x) |\nabla v(x)|^2 \, dx\geq S(\hat a) \left(
\int_{\R^N} \rho^{2^*}(x) |v(x)|^{2^*} \, dx\right)^{2/2^*}
\end{equation}
for all $v\in {\mathcal D}^{1,2}_{\rho}(\R^N)$.
\end{Lemma}

\begin{pf}
It follows from Lemma \ref{l:d12rho} and the change $u=\rho v$. 
\end{pf}

We also define the weighted Sobolev space $H^1_\rho(\Omega)$ as
the completion of $V_\rho(\Omega)$ with respect to the norm
\begin{equation*}
\|v\|_{H^1_\rho(\Omega)}:=\left(\int_\Omega \rho^2(x)|\nabla
v(x)|^2\, dx+\int_\Omega \rho^2(x)v^2(x)\, dx\right)^{\!\!1/2}
\end{equation*}
where $V_\rho(\Omega)$ denotes the space of all functions
$v\in C^\infty(\Omega)\cap H^1(\Omega)$ such that
$$
\overline{ \{x\in \Omega:v(x)\neq 0\}}^{\, \Omega}\subset
\Omega\setminus \widetilde \Sigma.
$$
For  any $q\geq 1$, we also denote as
$L^q(\rho^{2^*},\Omega)$ the weighted $L^q$-space endowed with
the norm
$$
\|u\|_{L^q(\rho^{2^*},\Omega)}:=
\bigg(\int_{\Omega}\rho^{2^*}(x)|u(x)|^q\,dx\bigg)^{\!\!1/q}.
$$

\begin{Lemma} \label{preliminare} Let $N\geq k\geq 3$,
  $\Omega\subset\R^N$ be a bounded open set such that $0\in\Omega$,
  $\hat a$ satisfy (\ref{eq:atetahat}) and (\ref{eq:53}),
 and $h$ satisfy \eqref{H1-2}. Let $V\in L^1_{\rm
    loc}(\Omega\setminus \widetilde \Sigma)$ such that
\begin{equation}\label{eq:52}
\sup_{v\in
  H^1_\rho(\Omega)\setminus\{0\}}\frac{\int_\Omega\rho^{2^*}(x)|V(x)|v^2(x)\,dx}
{\|v\|^2_{H^1_\rho(\Omega)}}<+\infty,
\end{equation}
and $v\in
  H^1_\rho(\Omega)\cap L^q(\rho^{2^*},\Omega)$, $q>2$, be a weak 
solution to
\begin{equation}\label{eq:prel}
  -{\rm div}(\rho^2(x)\nabla v(x))=\big(\rho^2(x)h(x)+ \rho^{2^*}(x)V(x)\big)
  v(x).
\end{equation}
If 
\begin{equation}\label{sum-V} 
  V_+\in L^s(\rho^{2^*},\Omega) \qquad
  {\rm for \ some} \ s>N/2,
\end{equation}
then for any $\Omega'\Subset \Omega$ such that $0\in \Omega'$, $v\in
L^{\frac{2^*q}{2}}(\rho^{2^*},\Omega')$ and
\begin{multline}\label{eq:66}
  \|v\|_{L^{\frac{2^*q}2}(\rho^{2^*},\Omega')}\leq
  S(\hat a)^{-\frac1q}  \|v\|_{L^{q}(\rho^{2^*},\Omega)}
\bigg(\frac{20}{C(q)} \frac{d}{(\mathop{\rm dist}
    (\Omega',\partial\Omega))^2}+ \frac{4(q-2)d}{(\mathop{\rm
      dist}(\Omega',\partial\Omega))^2}
  +\frac{4\ell_q}{C(q)}\bigg)^{\!\!\frac1q},
\end{multline}
where $C(q):=\min\big\{\frac14, \frac4{q+4}\big\}$,
$d$ is as in (\ref{eq:d}), and 
\begin{multline}
\ell_q= \max\Bigg\{ \bigg(\frac{\max\{16,q+4\}}{S(\hat a)}\,
\|V_+\|_{L^s(\rho^{2^*},\Omega)} ^{2s/N}\bigg)^{\!\!\frac{N}{2s-N}},\\
\frac{d
C_h^{2/\e}\big({\textstyle{\frac2{k-2}}}\big)^{\frac{2(2-\e)}{\e}}\binom{N}{k}^{2/\e}
  \big(1+\binom{N-k}{k}\big)^{2/\e}}
{(1-\Lambda(\hat a))^{\frac{2-\e}{\e}}}\,(\max\{16,q+4\})^{\frac{2-\e}{\e}}
\Bigg\}.
\end{multline}
\end{Lemma}

\begin{pf} Let $w\in \mathcal D^{1,2}_\rho(\Omega)$. 
Then by Lemma \ref{l:W-Sob} we have
\begin{align} \label{eq:4bk}
  &\qquad\int_{\Omega} \rho^{2^*}(x)V_+(x)|w(x)|^2\,dx\\
  &\notag \leq\ell_q\!\!\!\!\!\!
  \int\limits_{V_+(x) \leq\ell_q}\!\!\!\!\!\!\rho^{2^*}(x)
  |w(x)|^2\,dx +\!\!\!\!\!\!  \int\limits_{V_+(x)\geq
    \ell_q}\!\!\!\!\!\!
  \rho^{2^*-2}(x)V_+(x)\rho^{2}(x)|w(x)|^2\,dx\\
  \notag &\leq \ell_q\int_{\Omega}\rho^{2^*}(x) |w(x)|^2\,dx
  +\bigg(\int_{\Omega}\rho^{2^*}(x)|w(x)|^{2^*}dx\bigg)^{\!\!\frac2{2^*}}
  \bigg( \int\limits_{V_+(x)\geq \ell_q}
  \rho^{2^*}(x)V_+^{\frac N2}(x)\,dx\bigg)^{\!\!\frac2N}\\
  \notag &\leq
\ell_q\int_{\Omega}\rho^{2^*}(x) |w(x)|^2\,dx+
\frac1{S(\hat a)} \bigg(\int_{\Omega}
  {\textstyle{\rho^{2}(x)\big|\nabla w(x)\big|^2 \,dx}}\bigg)
  \bigg( \int\limits_{V_+(x) \geq \ell_q}
  \rho^{2^*}(x)V_+^{\frac N2}(x)\,dx\bigg)^{\!\!\frac2N}.
\end{align}
Next, H\"older inequality and the definition of $\ell_q$ yield
\begin{align} \label{HOLDER}
\int\limits_{V_+(x) \geq \ell_q} &
  \rho^{2^*}(x)V_+^{\frac N2}(x)\,dx\leq
  \bigg(\int_{\Omega}\rho^{2^*}(x) V_+^s(x) \,dx
  \bigg)^{\!\!\frac{N}{2s}}\bigg( \int\limits_{V_+(x)\geq
    \ell_q} \!\!\!\!\!\rho^{2^*}(x)\,dx
  \bigg)^{\!\!\frac{2s-N}{2s}}\\
  \notag&\leq
  \bigg(\int_{\Omega}\rho^{2^*}(x)V_+^s(x) \,dx\bigg)^{\!\!\frac{N}{2s}}\bigg(
  \int\limits_{V_+(x) \geq
    \ell_q}\bigg(\frac{V_+(x)}{\ell_q} \bigg)^{\!\!s}
  \rho^{2^*}(x)\,dx\bigg)^{\!\!\frac{2s-N}{2s}}\\
  \notag&\leq
  \|V_+\|_{L^s(\rho^{2^*},\Omega)}^{s}\,\ell_q^{-s+\frac
    N2}\leq\left(\min\bigg\{\frac{S(\hat a)}{16}, \frac{
    S(\hat a)}{q+4}\bigg\}\right)^{\!\!\frac N2},
\end{align}
Inserting \eqref{HOLDER} into \eqref{eq:4bk} we obtain
for any $w\in \mathcal D^{1,2}_\rho(\Omega)$
\begin{align} \label{stimaLS1}
\int_{\Omega} \rho^{2^*}(x) V_+(x)& |w(x)|^2 \, dx\\
&\notag\leq
\ell_q\int_{\Omega}\rho^{2^*}(x) |w(x)|^2\,dx
+\frac12\min\left\{\frac{1}{8},\frac{2}{q+4}\right\}
\int_\Omega \rho^{2}(x) |\nabla w(x)|^2 dx.
\end{align}
On the other hand, letting $\delta_q=\Big(C_hd \binom{N}{k}
\big(1+\binom{N-k}{k}\big)\Big)^{1/(2-\e)}\ell_q^{-1/(2-\e)}$, from
(\ref{H1-2}), (\ref{eq:ssw}), (\ref{eq:MPHardy}), 
(\ref{eq:1}),
(\ref{eq:48}), and (\ref{eq:d}), for every $w\in \mathcal D^{1,2}_\rho(\Omega)$
 we can estimate
\begin{align}\label{eq:51}
  &\int_{\Omega} \rho^{2}(x)|h(x)| |w(x)|^2 \, dx\\
  \notag &\leq C_h \bigg[ \delta_q^\e \bigg( \sum_{J\in{\mathcal A}_k}
  \int_{|x_J|\leq \delta_q} \frac{\rho^2(x)w^2(x)}{|x_{J}|^2}\, dx
  +\sum_{(J_1,J_2)\in \mathcal B_k} \int_{|x_{J_1}-x_{J_2}|\leq
    \delta_q} \frac{\rho^2(x)w^2(x)}{|x_{J_1}-x_{J_2}|^{2}}\, dx
  \bigg)\\
  \notag &\qquad+ \delta_q^{-2+\e}d \bigg( \sum_{J\in{\mathcal A}_k}
  \int_{|x_J|\geq \delta_q} \rho^{2^*}(x)w^2(x)\, dx
  +\sum_{(J_1,J_2)\in \mathcal B_k} \int_{|x_{J_1}-x_{J_2}|\geq
    \delta_q} \rho^{2^*}(x)w^2(x)\, dx
  \bigg)\bigg]\\
  \notag &\leq C_h {\textstyle{
\binom{N}{k} \Big(1+\binom{N-k}{k}\Big)}} \bigg( \delta_q^\e
  \left({\textstyle{\frac{2}{k-2}}}
\right)^{\!\!2}(1-\Lambda(\hat a))^{-1}\int_{\Omega}\rho^2(x)|\nabla
  w(x)|^2dx\\
\notag &\qquad+\delta_q^{-2+\e} d
\int_\Omega \rho^{2^*}(x)w^2(x)\,dx\bigg)\\
\notag &\notag\leq
\ell_q\int_{\Omega}\rho^{2^*}(x) |w(x)|^2\,dx
+\frac12\min\left\{\frac{1}{8},\frac{2}{q+4}\right\}
\int_\Omega \rho^{2}(x) |\nabla w(x)|^2 dx.
\end{align}
Summing up (\ref{stimaLS1}) and (\ref{eq:51}), we obtain
\begin{multline}\label{stimaLS}
\int_{\Omega} \big(\rho^{2^*}(x) V_+(x)+|h(x)|\rho^2(x)\big)
 |w(x)|^2 \, dx\\
\leq
2\ell_q\int_{\Omega}\rho^{2^*}(x) |w(x)|^2\,dx
+\min\left\{\frac{1}{8},\frac{2}{q+4}\right\}
\int_\Omega \rho^{2}(x) |\nabla w(x)|^2 dx
\end{multline}
for all $w\in \mathcal D^{1,2}_\rho(\Omega)$.
As in \cite{FFT,FMT2} we define $v^n:=\min\{n,|v|\}\in H^1_\rho(\Omega)$ and we introduce
a cut-off function $\eta\in C^\infty_c(\Omega)$ satisfying
$$
\eta\equiv 1 \ {\rm in } \ \Omega' \qquad {\rm and } \qquad
|\nabla \eta|\leq \frac{2}{{\rm dist}(\Omega',\partial \Omega)}.
$$
Testing \eqref{eq:prel} with $\eta^2 (v^n)^{q-2}v\in \mathcal D^{1,2}_\rho (\Omega)$
we obtain

\begin{align*}
  &(q-2)\int_{\Omega} \rho^{2}(x)\eta^2(x)(v^n(x))^{q-2}
  \alchi_{\{y\in\Omega:|v(y)|<n\}}(x)|\nabla|v|(x)|^2\,dx\\
  &\quad\qquad+ \int_{\Omega}\rho^{2}(x)\eta^2(x) (v^n(x))^{q-2}|\nabla
  v(x)|^2\,dx \\
  &=\!\int_{\Omega}\!(\rho^{2^*}(x)V(x)+\rho^2(x)h(x))\eta^2(x)
  |v(x)|^2(v^n(x))^{q-2}\,dx \\
  &\qquad\quad-2\!\int_{\Omega}\!\rho^{2}(x)
  \eta(x)(v^n(x))^{q-2}
  v(x)\nabla v(x)\cdot \nabla \eta(x)\,dx\\
  &\leq \int_{\Omega}\! (\rho^{2^*}(x) V_+(x)+\rho^2(x)|h(x)|) \eta^2(x)
  |v(x)|^2(v^n(x))^{q-2}\,dx\\
  &\quad \qquad
+2\int_{\Omega}\!\rho^{2}(x) |\nabla\eta(x)|^2(v^n(x))^{q-2}
  |v(x)|^2\,dx
 +\frac12 \int_{\Omega}\!\rho^{2}(x)
  \eta^2(x)(v^n(x))^{q-2}
  |\nabla v(x)|^2\,dx
\end{align*}
and hence
\begin{align} \label{eq:19bk}
  & (q-2)\int_{\Omega}\rho^{2}(x)\eta^2(x)(v^n(x))^{q-2}
  |\nabla v^n(x)|^2\,dx+ \frac 12 \int_{\Omega}\rho^{2}(x)\eta^2(x)
  (v^n(x))^{q-2}|\nabla v(x)|^2\,dx \\
 \notag & \leq \int_{\Omega}\! \big(\rho^{2^*}(x) V_+(x)
+\rho^2(x)|h(x)|\big) \eta^2(x)
  |v(x)|^2(v^n(x))^{q-2}\,dx \\
&\notag\qquad+2\int_{\Omega}\! \rho^{2}(x)
|\nabla\eta(x)|^2(v^n(x))^{q-2}
  |v(x)|^2\,dx.
 \end{align}
By direct computation we also have
\begin{align*} 
& \left|\nabla\left((v^n)^{\frac{q-2}2}\eta v\right)\right|^2
\leq
\frac{(q+4)(q-2)}4 (v^n)^{q-2} \eta^2 |\nabla v^n|^2 \\
\notag
& \qquad
+2\eta^2 (v^n)^{q-2}
|\nabla v|^2 +2|\nabla \eta|^2 (v^n)^{q-2} |v|^2+
\frac{q-2}2 (v^n)^q |\nabla \eta|^2
\end{align*}
and hence by \eqref{eq:19bk} we obtain
\begin{multline}\label{disug}
  C(q)\int_\Omega \rho^{2}(x)
  \left|\nabla\left((v^n)^{\frac{q-2}2}\eta v\right)\right|^2 \, dx\\
  \leq \int_{\Omega}\!\big( \rho^{2^*}(x)
  V_+(x)+\rho^2(x)|h(x)|\big)\eta^2(x) |v(x)|^2(v^n(x))^{q-2}\,dx
  \\
  +2(C(q)+1)\int_{\Omega}\! \rho^{2}(x) (v^n(x))^{q-2} |v(x)|^2|\nabla
  \eta(x)|^2\,dx+ C(q)\frac{q-2}2 \int_{\Omega}\! \rho^{2}(x)
  (v^n(x))^{q}|\nabla \eta(x)|^2\,dx .
\end{multline}
Applying  estimate \eqref{stimaLS} to the function
$w=\eta (v^n)^{\frac{q-2}2} v$, by \eqref{disug} we have
\begin{align*} 
& \frac{C(q)}2
\int_\Omega
\rho^{2}(x) \left|\nabla\left((v^n)^{\frac{q-2}2}\eta v\right)\right|^2
\, dx
\leq
2\ell_q
\int_{\Omega}\! \rho^{2^*}(x) \eta^2(x) (v^n(x))^{q-2} |v(x)|^2 \, dx
\\
\notag &
+2(C(q)+1)\int_{\Omega}\!\rho^{2}(x)
(v^n(x))^{q-2} |v(x)|^2|\nabla \eta(x)|^2\,dx+
C(q)\frac{q-2}2 \int_{\Omega}\!\rho^{2}(x)
  (v^n(x))^{q}|\nabla \eta(x)|^2\,dx.
\end{align*}
and this with Lemma \ref{l:W-Sob} and (\ref{eq:d}) implies
\begin{align*}
&
\left(\int_\Omega \rho^{2^*}(x) |v^n(x)|^{2^*\frac{q-2}2} |v(x)|^{2^*}
\eta^{2^*}(x) \, dx\right)^{\frac{2}{2^*}}
\leq \frac{4\ell_q}{C(q)S(\hat a)}
\int_{\Omega}\!\rho^{2^*}(x) \eta^2(x) (v^n(x))^{q-2} |v(x)|^2 \, dx \\
\notag & \quad
+\frac{4(C(q)+1)d}{C(q)S(\hat a)}\int_{\Omega}\! \rho^{2^*}(x)
(v^n(x))^{q-2} |v(x)|^2|\nabla \eta(x)|^2\,dx
+\frac{(q-2)d}{S(\hat a)} \int_{\Omega}\! \rho^{2^*}(x)
  (v^n(x))^{q}|\nabla \eta(x)|^2\,dx.
\end{align*}
The proof of the lemma then follows letting $n\to +\infty$.
\end{pf}

\begin{Theorem} \label{Main} Let $N\geq k\geq 3$, $\Omega\subset\R^N$
  be a bounded open set such that $0\in\Omega$, $\hat a$ satisfy
(\ref{eq:atetahat}) and (\ref{eq:53}), $h$ as in \eqref{H1-2}, and
  $V\in L^1_{\rm loc}(\Omega\setminus \widetilde \Sigma)$ verify
  (\ref{eq:52}). 
\begin{itemize}
 \item[i)] If  $V_+\in
L^s(\rho^{2^*},\Omega)$
for some $s>N/2$,
then for any $\Omega'\Subset \Omega$ there exists a positive constant
$$
C_\infty=C_\infty(N,k,\hat a, h, \|V_+\|_{L^s(\rho^{2^*},\Omega)},
{\rm dist}(\Omega',\partial \Omega), {\rm diam}(\Omega))
$$
depending only on $N,k,\hat a, h,\|V_+\|_{L^s(\rho^{2^*},\Omega)},
{\rm dist}(\Omega',\partial \Omega)$ and  ${\rm diam}(\Omega)$,
such that for every solution $u\in H^1(\Omega)$ of 
\begin{equation} \label{uV} -\Delta
  u(x)-\frac{\hat a(\frac{x}{|x|})}{|x|^2} u(x)= \big(h(x)+\rho^{2^*-2}(x)
  V(x)\big) u(x) \quad \text{in } \Omega,
\end{equation}
 there holds
 $\rho^{-1}u\in L^\infty(\Omega')$ and
$$
\|\rho^{-1}u\|_{L^\infty(\Omega')}\leq C_\infty
\|u\|_{L^{2^*}(\Omega)}.
$$
\item[ii)] If $V_+\in
L^{N/2}(\rho^{2^*},\Omega)$, then for any $\Omega'\Subset \Omega$
and for any $s\geq 1$ there exists a positive constant
$$
C_s=C_s(N,k,\hat a, h, V, s, {\rm
dist}(\Omega',\partial \Omega), {\rm diam}(\Omega))
$$
depending only on
$N,k,\hat a,h,V, s, {\rm
dist}(\Omega',\partial \Omega)$,  ${\rm diam}(\Omega)$, such
that every solution $u\in H^1(\Omega)$ to \eqref{uV} satisfies
$\rho^{-1}u\in L^s(\rho^{2^*},\Omega')$ and
$$
\|\rho^{-1}u\|_{L^s(\rho^{2^*},\Omega')}\leq C_s
\|u\|_{L^{2^*}(\Omega)}.
$$
\end{itemize}
\end{Theorem}

\begin{pf}
Let $u\in H^1(\Omega)$ be a weak solution of
\eqref{uV},  $\Omega'\Subset \Omega$, and  $R>0$  such that
$$
\Omega'\Subset \Omega'+B(0,2R)\Subset \Omega.
$$
We claim that the function $v(x):=\rho^{-1}(x) u(x)$
belongs to $H^1_\rho(\Omega'+B(0,2R))$. Indeed, arguing as in Lemma
\ref{l:d12density}, we can prove that $V_{\rho}(\Omega)$ is dense in
$H^1(\Omega)$, hence there exists a sequence $\{u_n\}_{n\in\N}\subset
V_{\rho}(\Omega)$ such that $u_n\rightarrow u$ in $H^1(\Omega)$.  If
$\eta\in C^\infty_c(\Omega)$ is a cut-off function such that
$\eta\equiv 1$ in $\Omega'+B(0,2R)$, from (\ref{eq:48}) it follows that 
\begin{align*}
\int_{\Omega}  & |\nabla (\eta(x)(u_n(x)-u_m(x)))|^2 dx
-\int_{\Omega} \frac{\hat a(\frac{x}{|x|})}{|x|^2}
\eta^2(x)(u_n(x)-u_m(x))^2 \, dx\\
& = \int_{\Omega} \rho^2(x)
|\nabla(\eta(x)\rho^{-1}(x)(u_n(x)-u_m(x)))|^2 dx.
\end{align*}
This shows that $\{\rho^{-1}u_n\}$ is a Cauchy sequence in
$H^1_\rho(\Omega'+B(0,2R))$ which then  converges to
$v(x)=\rho^{-1}(x)u(x)$. In particular $v\in
H^1_\rho(\Omega'+B(0,2R))$.

By direct computation one also sees that $v$ is a weak solution of
\eqref{eq:prel}. By  Lemma \ref{preliminare}, 
proceeding exactly as in the proofs of \cite[Theorem
9.3]{FFT} and \cite[Theorem 1.2]{FMT2}, we arrive to the conclusion.
\end{pf}

\begin{remark}\label{r:erratum}
  The statement of Theorem 9.4 in our previous paper \cite{FFT} should
  be corrected as in the statement of Theorem \ref{t:BK}.  The missing
  point in Theorem 9.4. as it was stated in \cite{FFT} relies in the
  fact that the constant $\widetilde C_\infty$ such that
  $\||x|^{-\sigma}u\|_{L^\infty(\Omega')}\leq \widetilde
  C_\infty\|u\|_{L^{2^*}(\Omega)}$ depends on $u$, more precisely on
  the distribution function of $f(x,u)/u$.

In a similar way, the statements of Theorems 9.3 and 10.4 should be
corrected as in Theorem \ref{Main} above, i.e. the constant $C_s$
(respectively $C_{s,2}$) appearing in the statement (ii) of Theorem
9.3 (respectively 10.4) depends on $(\Re(V))_+$ (more precisely on its
distribution function) and not only on its
$L^{N/2}(\rho^{2^*},\Omega)$-norm (respectively
$L^s(\rho^{p},\Omega)$-norm) as incorrectly stated in
\cite{FFT}.

Anyway, the proofs of Theorems 9.3 and 9.4 contained in \cite{FFT}  
 are correct and lead to analogous conclusion as those stated in 
 Theorems \ref{t:BK} and \ref{Main} of the present paper. Moreover 
all the proofs and statements in the rest of the paper
\cite{FFT} are not affected by these corrections.

\end{remark}


\medskip\noindent
\begin{pfn}{Theorem \ref{t:BK}}
Let us define
$$
V(x):=\left\{
\begin{array}{ll}
\rho^{2-2^*}(x)\Big( \frac{f(x,u(x))}{u(x)}
-
\sum_{J\in{\mathcal
A}_k}\frac{\alpha_J^-}{|\theta_J|^2}-\sum_{(J_1,J_2)\in \mathcal
B_k}\frac{\alpha_{J_1\, J_2}^-}{|\theta_{J_1}-\theta_{J_2}|^2}\Big), &
 \quad {\rm if} \ u(x)\neq 0, \\
\rho^{2-2^*}(x)\Big(-
\sum_{J\in{\mathcal
A}_k}\frac{\alpha_J^-}{|\theta_J|^2}-\sum_{(J_1,J_2)\in \mathcal
B_k}\frac{\alpha_{J_1\, J_2}^-}{|\theta_{J_1}-\theta_{J_2}|^2}\Big), & \quad {\rm if} \ u(x)=0,
\end{array}
\right.
$$
where $\alpha_J^-=\max\{-\alpha_J,0\}$ and
$\alpha_{J_1J_2}^-=\max\{-\alpha_{J_1J_2},0\}$. By 
(\ref{F}) 
 and the Sobolev embedding 
$H^1(\Omega)\subset L^{2^*}(\Omega)$, 
we have that $V^+\in L^{N/2}(\rho^{2^*},\Omega)$ and  $u$ weakly solves 
$$
-\Delta
  u(x)-\frac{\hat a(\frac{x}{|x|})}{|x|^2} u(x)= \big(h(x)+\rho^{2^*-2}(x)
  V(x)\big) u(x) \quad \text{in } \Omega.
$$ 
From part ii) of Theorem \ref{Main}, it follows that $\rho^{-1}u\in
L^s(\rho^{2^*},\Omega')$ for any $\Omega'\Subset \Omega$ and for any
$s\geq 1$.  By (\ref{F}) we deduce that $V^+\in
L^{s}(\rho^{2^*},\Omega')$ for all $s\geq\frac{N-2}4$ and in
particular for some $s>N/2$.  The proof of the theorem follows now by
part i) of Theorem \ref{Main}.
\end{pfn}

\appendix
 \section*{Appendix}
 \setcounter{section}{1}
 \setcounter{Theorem}{0}

 To prove Theorem \ref{t:pohozaev} we used, for the approximating
 problems, a Pohozaev-type identity (see (\ref{eq:pohoapp})), whose
 proof is quite classical (see e.g. \cite{PS,STRU}) and requires just
 few adaptations due to the presence of a singularity. For the sake of
 completeness we give below a proof.

\begin{Proposition}\label{p:poho}
Let $\Omega\subset\R^N$, $N\geq 3$,
  be a bounded open set such that $0\in\Omega$. Let $b\in L^\infty(\SN)$,
  $h\in L^\infty(\Omega)$, and let $f$ satisfy \eqref{F}. Denote by
$\nu=\nu(x)$ the unit outer normal vector
  $\nu(x)=\frac{x}{|x|}$. If $u$ is a $H^1(\Omega)$-weak solution
  to $\mathcal L_{b} u=h(x)\, u+f(x,u)$ in $\Omega$ and $r_0>$ is such that
  $B_{r_0}\subseteq \Omega$, then for a.e. $r\in (0,r_0)$
\begin{multline}\label{eq:pohoapp-A}
  -\frac{N-2}2\int_{B_r}\bigg[ |\nabla u(x)|^2
  -\frac{b(\frac{x}{|x|})}
  {|x|^2}u^2(x)\bigg]\,dx
  +\frac{r}{2}\int_{\partial B_r}\bigg[ |\nabla u(x)|^2
  -\frac{b(\frac{x}{|x|})}{|x|^2}u^2(x)\bigg]\,dS\\
  =r\int_{\partial B_r}\bigg|\frac{\partial u}{\partial
    \nu}\bigg|^2\,dS+ \int_{B_r}h(x)u(x)\,(x\cdot {\nabla
    u(x)})\big)\,dx \\
+r\int_{\partial B_r} F(x,u(x))\, dS-\int_{B_r}
[\nabla_xF(x,u(x))\cdot x+NF(x,u(x))]\, dx \ .
\end{multline}

\end{Proposition}

\begin{pf}
  By classical Brezis-Kato \cite{BK} estimates, bootstrap,
  and elliptic regularity theory,  \eqref{F} and the boundedness of
the coefficients $b,h$ imply that
$u\in H^2_{\rm loc}(\Omega\setminus \{0\})\cap
C^{1,\alpha}_{\rm loc}(\Omega\setminus\{0\})$ for any $\alpha\in
(0,1)$.
Therefore by \eqref{F} and Hardy inequality, we have
\begin{multline*}
\int_0^r \left[\int_{\partial B_s} \left(|\nabla
u(x)|^2+\frac{u^2(x)}{|x|^2}+\left|\frac{\partial u}{\partial
\nu}(x)\right|^2+|F(x,u(x))|\right)  \, dS\right] \, ds \\
=\int_{B_r} \left(|\nabla
u(x)|^2+\frac{u^2(x)}{|x|^2}+\left|\frac{\partial u}{\partial
\nu}(x)\right|^2+|F(x,u(x))|\right) \, dx<+\infty
\end{multline*}
and hence there exists a decreasing sequence $\{\delta_n\}\subset
(0,r)$ such that $\lim_{n\to +\infty} \delta_n=0$ and
\begin{equation} \label{bordo-delta}
\delta_n \int_{\partial B_{\delta_n}} \left(|\nabla
u(x)|^2+\frac{u^2(x)}{|x|^2}+\left|\frac{\partial u}{\partial
\nu}(x)\right|^2+|F(x,u(x))|\right) \, dS \longrightarrow 0
\text{\quad as } n\to +\infty \, .
\end{equation}
Multiplying equation $\mathcal L_{b} u=h(x)\, u+f(x,u)$ 
 by $x\cdot \nabla u(x)$ and
integrating over $B_r\setminus B_{\delta_n}$, it follows that
\begin{align} \label{poho-A1}
\int_{B_r\setminus B_{\delta_n}} & \nabla u(x)\cdot \nabla(x\cdot
\nabla u(x)) \, dx-\int_{B_r\setminus B_{\delta_n}}
\frac{b(\frac{x}{|x|})}{|x|^2} u(x)(x\cdot \nabla u(x))  \, dx \\
\notag & =r\int_{\partial B_r} \left|\frac{\partial u}{\partial
\nu}\right|^2 \, dS-\delta_n \int_{\partial B_{\delta_n}}
\left|\frac{\partial u}{\partial \nu}\right|^2 \, dS +
\int_{B_r\setminus B_{\delta_n}} h(x)u(x)(x\cdot \nabla u(x)) \,
dx \\ \notag & \quad +\int_{B_r\setminus B_{\delta_n}}
f(x,u(x))(x\cdot \nabla u(x)) \, dx \, .
 \end{align}
Standard integration by parts shows that
\begin{multline} \label{poho-A2}
\int_{B_r\setminus B_{\delta_n}} \nabla u(x)\cdot \nabla(x\cdot
\nabla u(x)) \, dx \\
=-\frac{N-2}2 \int_{B_r\setminus
B_{\delta_n}} |\nabla u(x)|^2 \, dx+\frac{r}2 \int_{\partial B_r}
|\nabla u(x)|^2\, dS-\frac{\delta_n}2 \int_{\partial \delta_n}
|\nabla u(x)|^2 \, dS\, .
\end{multline}
Passing in radial coordinates $r=|x|$, $\theta=\frac{x}{|x|}$ and
observing that $\partial_r u(r,\theta)=\nabla u(r\theta)\cdot
\theta$, we obtain
\begin{align*}
\int_{B_r\setminus B_{\delta_n}} & \frac{b(\frac{x}{|x|})}{|x|^2}
u(x)(x\cdot \nabla u(x))  \, dx = \int_{\SN} b(\theta)\left(
\int_{\delta_n}^r s^{N-2} u(s\theta)\partial_s u(s\theta) \,
ds\right)\, dS(\theta) \\
&=\int_{{\mathbb S}^{N-1}}b(\theta)\bigg(
r^{N-2}u^2(r\theta)-\delta_n^{N-2}u^2(\delta_n\theta)\\
&\hskip3cm- (N-2)\int_{\delta_n}^rs^{N-3}u^2(s\theta)\,ds
-\int_{\delta_n}^rs^{N-2}u(s\theta)\partial_s u(s\theta)\,ds\bigg)dS(\theta)\\
&=r\int_{\partial B_r}\frac{b(\frac x{|x|})}{|x|^2}
   u^2(x)\,dS-\delta_n\int_{\partial B_{\delta_n}}\frac{b(\frac x{|x|})}{|x|^2}
   u^2(x)\,dS\\
&\hskip2cm-(N-2)\int_{B_r\setminus B_{\delta_n}}\frac{b(\frac
x{|x|})}{|x|^2}
    u^2(x)\,dx-  \int_{B_r\setminus B_{\delta_n}}\frac{b(\frac x{|x|})}{|x|^2}
    u(x)\,(x\cdot\nabla u(x))\,dx \ ,
\end{align*}
which yields
\begin{multline}\label{poho-A3}
  \int_{B_r\setminus B_{\delta_n}}\frac{b(\frac x{|x|})}{|x|^2}
   u(x)\,(x\cdot \nabla
    u(x))\,dx\\
\quad= -\frac{N-2}2\int_{B_r\setminus
    B_{\delta_n}}\frac{b(\frac x{|x|})}{|x|^2}
  u^2(x)\,dx+\frac r2\int_{\partial B_r}\frac{b(\frac x{|x|})}{|x|^2}
  u^2(x)\,dS
   -\frac{\delta_n}2\int_{\partial
    B_{\delta_n}}\frac{b(\frac x{|x|})}{|x|^2}
  u^2(x)\,dS \, .
\end{multline}
By \eqref{F} and the fact that $u\in C^{1,\alpha}_{\rm
loc}(\Omega\setminus\{0\})$ we obtain
\begin{align} \label{poho-A4}
r \int_{\partial B_r} & F(x,u(x))\, dS - \delta_n \int_{\partial
B_{\delta_n}} F(x,u(x))\, dS = \int_{B_r\setminus B_{\delta_n}}
{\rm div} \big (F(x,u(x))x\big) \, dx \\
&\notag = \int_{B_r\setminus B_{\delta_n}} \left[\nabla_x
F(x,u(x))\cdot x+NF(x,u(x))\right]\, dx +\int_{B_r\setminus B_{\delta_n}}
f(x,u(x))(\nabla u(x)\cdot x)\, dx
\end{align}
Letting $n\to+\infty$, \eqref{eq:pohoapp-A} follows by
(\ref{bordo-delta}--\ref{poho-A4}).
\end{pf}


\begin{thebibliography}{99}


\bibitem{AFP} B. Abdellaoui, V. Boumediene, I.  Peral, {\it Some
    remarks on systems of elliptic equations doubly critical in the
    whole $\R^N$},  Calc. Var. Partial Differential Equations 34
    (2009), no. 1, 97--137.

\bibitem{almgren} F. J. Jr. Almgren, {\it $Q$ valued functions
    minimizing Dirichlet's integral and the regularity of area
    minimizing rectifiable currents up to codimension two},
  Bull. Amer. Math. Soc. 8 (1983), no. 2, 327--328.

\bibitem{BBR} M. Badiale, V. Benci, S.  Rolando, {\it A nonlinear
    elliptic equation with singular potential and applications to
    nonlinear field equations}, J. Eur. Math. Soc.  9 (2007),
    no. 3, 355--381.

  \bibitem{BR} M. Badiale, S.  Rolando, {\it Elliptic problems with
      singular potential and double-power nonlinearity},
      Mediterr. J. Math.  2 (2005), no. 4, 417--436.

\bibitem{BT} M. Badiale, G. Tarantello, {\it A Sobolev-Hardy inequality
with applications to a nonlinear elliptic equation arising in astrophysics},
  Arch. Ration. Mech. Anal.  163  (2002),  no. 4, 259--293.

\bibitem{Baum-Juhl} H. Baum, A. Juhl, {\it Conformal differential
    geometry, Q-Curvature and conformal holonomy}, Oberwolfach
  Seminars, Birkh\"auser Verlag, 2010.

\bibitem{esteban} R. Bosi, J. Dolbeault, M. J. Esteban, {\it Estimates
 for the optimal constants in multipolar Hardy inequalities for
 Schr\"odinger and Dirac operators}, 
Commun. Pure Appl. Anal. 7 (2008), no. 3, 533--562. 

\bibitem{BK} H. Br{e}zis, T. Kato, {\it Remarks on the
    {S}chr\"odinger operator with singular complex potentials}, J.
  Math. Pures Appl. (9) 58 (1979), no.~2, 137--151.

\bibitem{buslaev_levin} V.S. Buslaev,S.B.  Levin, {\it Asymptotic
    behavior of the eigenfunctions of the many-particle Schr\"odinger
    operator. I. One-dimensional particles},  Spectral theory of
    differential operators, 55--71, Amer. Math. Soc. Transl. Ser. 2,
    225, Amer. Math. Soc., Providence, RI, 2008.

\bibitem{CKN} L. Caffarelli, R. Kohn, L. Nirenberg, {\it First order
    interpolation inequalities with weights}, Compositio Math. 53
  (1984), no. 3, 259--275.


\bibitem{CatrinaWang} F. Catrina, Z.-Q. Wang, {\it On the
    {C}affarelli-{K}ohn-{N}irenberg inequalities: sharp constants,
    existence (and nonexistence), and symmetry of extremal functions},
  Comm. Pure Appl. Math. 54 (2001), no. 2, 229--258.

\bibitem{CSW} J. Chabrowski, A. Szulkin, M. Willem, {\it
 Schr\"odinger equation with multiparticle potential and critical nonlinearity},
 Topol. Meth. Nonl. Anal. 34 (2009), 201-211.

\bibitem{chang} S.-Y. A. Chang, {\it Conformal Invariants and Partial Differential Equations}, 
Colloquium Lecture notes, BAMS, 42, (2005), no. 3, 365--393.

\bibitem{chen} J. Chen, {\it Multiple positive solutions for a
    semilinear equation with prescribed singularity,} J. Math. Anal.
  Appl.  305 (2005), no. 1, 140--157.

\bibitem{duyckaerts} T. Duyckaerts, {\it
In\'egalit\'es de r\'esolvante pour l'op\'erateur de
 Schr\"odinger avec potentiel multipolaire critique},
 Bulletin  Bull. Soc. Math. France 134 (2006), no. 2, 201--239.

\bibitem{egnell} H. Egnell, {\it Elliptic boundary value problems with
    singular coefficients and critical nonlinearities,} Indiana Univ.
  Math. J. 38 (1989), no. 2, 235--251.

\bibitem{FFT} V. Felli, A. Ferrero, S. Terracini, {\it Asymptotic
    behavior of solutions to Schr\"odinger equations near an isolated
    singularity of the electromagnetic potential}, Journal of the
  European Mathematical Society, 13 (2011), no. 1, 119--174.

\bibitem{addendumFFT} V. Felli, A. Ferrero, S. Terracini, {\it A note
    on local asymptotics of solutions to singular elliptic equations
    via monotonicity methods}, Preprint 2011, available at
  {\tt http://arxiv.org/abs/1007.4434}.

\bibitem{FMT1} V. Felli, E.M. Marchini, S. Terracini, {\it On
    {S}chr\"odinger operators with multipolar inverse-square
    potentials}, Journal of Functional Analysis 250 (2007), 265--316.


\bibitem{FMT2} V. Felli, E.M. Marchini, S. Terracini, {\it On the
    behavior of solutions to {S}chr\"odinger equations with
    dipole-type potentials near the singularity}, Discrete
  Contin. Dynam. Systems 21 (2008), 91--119.

\bibitem{FMT3} V. Felli, E.M. Marchini, S. Terracini, {\it On
    {S}chr\"odinger operators with multisingular inverse-square
    anisotropic potentials}, Indiana Univ. Math. Journal 58 (2009),
  no. 2, 617--676.


\bibitem{FS3} V. Felli, M. Schneider,
{\it A note on regularity of solutions to degenerate elliptic equations of
Caffarelli-Kohn-Nirenberg type,}
Adv. Nonlinear Stud. 3 (2003), no. 4, 431--443.


\bibitem{FT} V. Felli, S. Terracini, {\it Elliptic equations with
    multi-singular inverse-square potentials and critical
    nonlinearity,} Comm. Partial Differential Equations 31 (2006),
  no. 1-3, 469--495.

\bibitem{FG} A. Ferrero, F. Gazzola, {\it Existence of solutions for
    singular critical growth semilinear elliptic equations,}
  J. Differential Equations 177  (2001),  no. 2, 494--522.

\bibitem{GP} J. Garc\'{\i}a Azorero, I. Peral,
{\it Hardy Inequalities and some critical elliptic and parabolic problems,}
J. Diff. Equations 144 (1998), no. 2, 441--476.


\bibitem{GL} N. Garofalo, F.-H.  Lin, {\it Monotonicity properties of
    variational integrals, $A\sb p$ weights and unique continuation},
  Indiana Univ. Math. J.  35 (1986), no. 2, 245--268.

\bibitem{HHLT} M. Hoffmann-Ostenhof, T.  Hoffmann-Ostenhof, A. Laptev,
J. Tidblom, {\it Many-particle Hardy inequalities},  J. Lond. Math. Soc. (2)
  77  (2008),  no. 1, 99--114.

\bibitem{HS} W. Hunziker, I.  Sigal, {\it The quantum $N$-body problem},
  J. Math. Phys.  41  (2000),  no. 6, 3448--3510.


\bibitem{Jan} E. Jannelli,
{\it The role played by space dimension in elliptic critical problems,}
J. Differential Equations 156 (1999), no. 2, 407--426.


\bibitem{lesch} M. Lesch, {\it Operators of Fuchs type, conical
    singularities, and asymptotic methods}, Teubner Texts in
  Mathematics, 136. B. G. Teubner Verlagsgesellschaft mbH, Stuttgart,
  1997.

\bibitem{LT} E. H. Lieb, W. E. Thirring, {\it Gravitational collapse in
    quantum mechanics with relativistic kinetic energy},  Ann. Physics
    155 (1984), no. 2, 494--512.

\bibitem{MFS} G. Mancini, I.  Fabbri, K.  Sandeep, {\it
Classification of solutions of a critical Hardy-Sobolev operator},
  J. Differential Equations  224  (2006),  no. 2, 258--276.

\bibitem{Mazja} V. G. Maz'ja, {\it Sobolev Spaces},  Springer Series
 in Soviet Mathematics, Springer-Verlag, Berlin, 1985.

\bibitem{mazzeo91} R. Mazzeo, {\it Elliptic theory of differential
    edge operators. I}, Comm. Partial Differential Equations 16
  (1991), no. 10, 1615--1664.

\bibitem{mazzeo91-2} R. Mazzeo, {\it Regularity for the singular
    Yamabe problem}, Indiana Univ. Math. J. 40 (1991), no. 4,
  1277--1299.

\bibitem{Musina} R. Musina, {\it Ground state solutions of a critical
    problem involving cylindrical weights}, Nonlin. Anal. 68 (2008),
  3972--3986.

\bibitem{pinchover94} Y. Pinchover, {\it On positive Liouville
    theorems and asymptotic behavior of solutions of Fuchsian type
    elliptic operators,} Ann. Inst. H. Poincar\'e Anal. Non
  Lin\'eaire 11 (1994), no. 3, 313--341.

\bibitem{PS} P. Pucci, J. Serrin, {\it A General Variational Identity,}
Indiana Univ. Math. J. 35 (1986), 681--703.

\bibitem{SSW} S. Secchi, D. Smets, M. Willem,
{\it Remarks on a Hardy-Sobolev inequality,}
C. R. Math. Acad. Sci. Paris 336 (2003), no. 10, 811--815.

\bibitem{SM} D. Smets, {\it Nonlinear Schr\"{o}dinger equations with
Hardy potential and critical nonlinearities,}
Trans. AMS 357 (2005), 2909--2938.

\bibitem{STRU} M. Struwe, {\it Variational methods and applications to
    nonlinear partial differential equations and Hamiltonian systems,}
  Springer-Verlag, Berlin/New York (1990).

  \bibitem{terracini96} S. Terracini, {\it On positive entire solutions
      to a class of equations with singular coefficient and critical
      exponent,} Adv. Diff. Equa. 1 (1996), no. 2, 241--264.

\bibitem{wz} Z.-Q. Wang, M. Zhu, {\it {H}ardy inequalities with
    boundary terms}, Electron. J. Differential Equations 2003, No. 43.

\bibitem{wolff} T. H. Wolff, {\it A property of measures in $\R^ N$
    and an application to unique continuation},  Geom. Funct. Anal.  2
    (1992), no. 2, 225--284.


\end{thebibliography}
\end{document}